\documentclass[letterpaper,11pt,reqno]{amsart}

\usepackage[margin=0.75in,marginparwidth=2cm,marginparsep=0.5cm]{geometry}
\usepackage[T1]{fontenc}
\usepackage{amsmath,amssymb,amscd,amsthm,amsxtra,mathtools}
\usepackage{xcolor}
\usepackage[hidelinks,pagebackref]{hyperref}
\renewcommand*{\backrefalt}[4]{%
  \ifcase #1\relax
  \or\space{\footnotesize (Cited on page #2)}%
  \else\space{\footnotesize (Cited on pages #2)}%
  \fi
}
\AddToHook{cmd/appendix/before}{%
    \crefalias{section}{appendix}%
}
\usepackage{aliascnt}
\usepackage[capitalise]{cleveref}
\usepackage{enumitem}
\usepackage{graphicx}

\allowdisplaybreaks[2]
\numberwithin{equation}{section}

\newtheorem{theorem}{Theorem}[section]

\newaliascnt{lemma}{theorem}
\newtheorem{lemma}[lemma]{Lemma}
\aliascntresetthe{lemma}

\newaliascnt{proposition}{theorem}
\newtheorem{proposition}[proposition]{Proposition}
\aliascntresetthe{proposition}

\newaliascnt{remark}{theorem}
\newtheorem{remark}[remark]{Remark}
\aliascntresetthe{remark}

\newaliascnt{condition}{theorem}
\newtheorem{condition}[condition]{Condition}
\aliascntresetthe{condition}

\newaliascnt{assumption}{theorem}

\aliascntresetthe{assumption}

\newaliascnt{claim}{theorem}

\aliascntresetthe{claim}

\newaliascnt{example}{theorem}
\newtheorem{example}[example]{Example}
\aliascntresetthe{example}

\newaliascnt{definition}{theorem}
\newtheorem{definition}[definition]{Definition}
\aliascntresetthe{definition}

\newaliascnt{corollary}{theorem}
\newtheorem{corollary}[corollary]{Corollary}
\aliascntresetthe{corollary}

\crefname{theorem}{theorem}{theorems}
\Crefname{theorem}{Theorem}{Theorems}
\crefname{maintheorem}{theorem}{theorems}
\Crefname{maintheorem}{Theorem}{Theorems}
\crefname{lemma}{lemma}{lemmas}
\Crefname{lemma}{Lemma}{Lemmas}
\crefname{proposition}{proposition}{propositions}
\Crefname{proposition}{Proposition}{Propositions}
\crefname{remark}{remark}{remarks}
\Crefname{remark}{Remark}{Remarks}
\crefname{condition}{condition}{conditions}
\Crefname{condition}{Condition}{Conditions}
\crefname{assumption}{assumption}{assumptions}
\Crefname{assumption}{Assumption}{Assumptions}
\crefname{claim}{claim}{claims}
\Crefname{claim}{Claim}{Claims}
\crefname{example}{example}{examples}
\Crefname{example}{Example}{Examples}
\crefname{definition}{definition}{definitions}
\Crefname{definition}{Definition}{Definitions}
\crefname{corollary}{corollary}{corollaries}
\Crefname{corollary}{Corollary}{Corollaries}
\crefname{notation}{notation}{notations}
\Crefname{notation}{Notation}{Notations}

\DeclareMathOperator*{\supp}{supp}

\DeclareMathOperator{\Var}{Var}
\DeclareMathOperator{\Vol}{Vol}
\DeclareMathOperator{\dist}{dist}

\newcommand{\R}{\mathbb R}
\newcommand{\Z}{\mathbb Z}
\newcommand{\T}{\mathbb T}
\newcommand{\N}{\mathbb N}
\newcommand{\C}{\mathbb C}
\newcommand{\E}{\mathbb E}
\newcommand{\PP}{\mathbb P}

\newcommand{\fH}{\mathfrak H}

\newcommand{\eps}{\varepsilon}

\newcommand{\wpsi}{\widehat\psi}
\newcommand{\wE}{\widehat E}
\newcommand{\wmu}{\widehat\mu}

\makeatletter
\@namedef{subjclassname@2020}{\textup{2020} Mathematics Subject Classification}
\makeatother

\begin{document}
	
	\title[Spectral delocalization and central limit theorems]
	{Spectral delocalization and central limit theorems for\\ Fourier functionals of stationary Gaussian fields}
	
	\author[L. Maini]{Leonardo Maini}
	
	\address{
		Leonardo Maini,
		Dipartimento di Matematica\\
		Universit\`a di Roma Tor Vergata\\
		Via della Ricerca Scientifica, 1, 00133 Roma, Italy 
	}
	
	\email{maini@mat.uniroma2.it}
	
	

\begin{abstract}
We characterize Gaussian fluctuations of Fourier-Hermite
functionals of stationary Gaussian fields on locally compact
Abelian groups via spectral delocalization: conditional on
the total spectral sum being near the probed frequency,
no partial sum concentrates. Under suitable assumptions,
this reduces Gaussian convergence to small-ball estimates
for convolution powers of the spectral measure, without
model-specific contraction estimates.

On $\mathbb T^d$, the criterion yields CLTs at every fixed
frequency for band-limited fields and random trigonometric
polynomials, and extends the known failure of Gaussian
convergence for planar arithmetic random waves from the
fourth chaos to all even orders $q\ge4$. On $\mathbb R^d$
and $\mathbb Z^d$, we establish a Breuer--Major theorem
for Fourier functionals at every frequency with nonzero
limiting variance. Cyclical long-memory models exhibit
frequency-dependent Gaussian convergence.

We also establish CLTs and variance asymptotics for
singular continuous spectral measures, including those
of stationary Gaussian solutions of linear PDEs. At zero
frequency, weak Ahlfors regularity of the spectral support
and an $L^2$ density suffice for second-chaos CLTs,
covering fractal sets, the Airy surface of Gaussian sea
waves, and smooth compact manifolds. For a class of smooth
hypersurfaces, we obtain the same variance rates and CLTs known for classical nonlinear functionals of Berry's
random wave model. Finally, for random waves with spherical
spectral densities bounded above and below, we classify Gaussian convergence at every
frequency $\|z\|\le q$: it fails only for the planar second
chaos when $0<\|z\|<2$.
\end{abstract}
	\vspace{-10mm}
	\maketitle
    \vspace{-10mm}
	\tableofcontents
	
	\section{Introduction}
	
	Gaussian fields arise in a broad range of models and applications, see e.g.  \cite{AzaisLeonWschebor2011,EstradeFournier, LerarioMarinucciRossiStecconi2025,PistolatoStecconi2026,CaponeraMarinucciVidotto2026,CelliPeccati,DiLilloMainiMarinucci2026,DiLilloMarinucciSalviVigogna2026} and the references therein.
	
	One possible way to investigate their properties is to study sequences of their functionals. Examples of functionals include partial sums and integral functionals, \cite{BreuerMajor,DobrushinMajor,Taqqu1979}; geometric functionals such as numbers of critical points or volumes and components counts of excursion sets \cite{AnconaGassLetendreStecconi,BeliaevCammarotaWigman,BeliaevMcAuleyMuirhead,DiBernardinoShevchenkoTodino2026,NourdinPeccatiRossi}; Fourier coefficients and periodograms \cite{MarinucciPeccatiCompact,PW10,PZ19,FT86}. 
    
	In this work, we focus on stationary Gaussian fields on locally compact Abelian groups and on the asymptotic properties of a class of functionals defined below. Popular examples in the literature include integral functionals over growing domains in \(\mathbb R^d\), growing partial sums in \(\mathbb Z^d\), and Fourier coefficients of subordinated Gaussian fields on compact groups such as the flat torus $\T^d$. In all these situations, stationarity provides a spectral measure on the dual group, and the asymptotic behavior of the corresponding functionals can be investigated through spectral conditions.
	
	To make this more precise, let us fix a second countable, Hausdorff,
	locally compact Abelian group \(E=(E,+)\) and a Haar measure \(\mu\) on
	\(E\). Let \(\widehat E\) denote the corresponding Pontryagin dual group
	(see \Cref{sec:lca-spectral}). We consider families of stationary real-valued (jointly measurable $\Omega\times E \to \R$)
	Gaussian fields
	\[
	X_T=(X_T(x))_{x\in E},\qquad T>0.
	\]
	We assume that $X_T$ is centered, with unit variance, and
	has continuous covariance functions \(C_T:E\to\R\),
	\[
	C_T(x-y):=\E\left[X_T(x)X_T(y)\right],\qquad C_T(0)=1.
	\]
	By Bochner's theorem there is a unique symmetric Borel probability measure
	\(G_T\) on \(\widehat E\), called the \emph{spectral measure} of \(X_T\),
	such that\footnote{The exponential notation is literal for
		\(\R^d,\Z^d,\T^d\) and is only a convenient notation otherwise; see \Cref{sec:lca-spectral}.}
	\[
	C_T(x)=\widehat G_T(x)
	=\int_{\widehat E}e^{i\langle x,\lambda\rangle}\,G_T(d\lambda).
	\]
	A widely studied class of functionals\footnote{$X$ is stationary and jointly measurable, so $F_T^{(z)}$ is well defined at least in $L^2(\Omega)$, see e.g. \cite{MainiSPA}.} is
	\[
	F_T:=\int_E\varphi(X_T(x))\psi_T(x)\,\mu(dx),
	\]
	where $\mu$ is a fixed Haar measure on $E$, \(\psi_T\in L^1(E,\R,\mu)\) and \(\varphi\in L^2(\gamma):=L^2(\R,\R,\gamma)\), where $\gamma$ denotes the standard Gaussian measure. Regarding applications, we focus on
	two regimes:
	\begin{enumerate}[label=\textbf{\alph*)},leftmargin=2em]
		\item\label{high frequency}
		\(E=\T^d=\R^d/(2\pi\Z^d)\), $\widehat E=\Z^d$, $\mu(\T^d)=1$, \(\psi_T\equiv1\), where \(X_T\) and \(G_T\) vary with \(T\). In this case $F_T$ is
		\[
		F_T=\int_{\T^d}\varphi(X_T(x))\,\mu(dx).
		\]
		
		\item\label{large domain}
		\(E=\R^d\), with $\widehat E=\R^d$ and $\mu$ Lebesgue measure; or \(E=\Z^d\), with $\widehat E=\T^d$ and $\mu$ counting measure. Here \(X_T=X\), \(G_T=G\), and
		\(\psi_T(x)=\psi(x/T)\), where $\psi$ is in the Schwartz space $\mathcal S(\R^d)$ or \(\psi=\mathbf1_D\), $D\subseteq \R^d$. In the latter case, we recover growing-domain integrals and partial sums. In general, we have
        \[
        F_T=\int_{E}\varphi(X(x))\,\psi(x/T)\,\mu(dx)\,.
        \]
	\end{enumerate}
	Motivated by statistical applications, see e.g. \cite{FT86,FT87}, and by the study of Fourier transforms and Fourier coefficients of random fields, see e.g. \cite{MarinucciPeccatiCompact,PW10,PZ19,RosenblattFourier}, we also consider
	$$
	F_T^{(z)}
	:=\int_E\varphi(X_T(x))\psi_T(x)e^{-i\langle x,z\rangle}\,\mu(dx),
	\qquad z\in\widehat E.
	$$
	The functional is real-valued
	when \(z=-z\), and may be complex-valued if $z\neq -z$. 
	A very well-known and exploited fact in the literature is that the analysis of the functionals $F_T^{(z)}$ can be reduced to the study of their chaotic components. Up to constants depending on \(\varphi\), these components are 
	\[
	F_{T,q}^{(z)}
	=
	\int_E
	H_q(X_T(x))\,\psi_T(x)\,
	e^{-i\langle x,z\rangle}\,\mu(dx),
	\]
	where \(q\ge2\) and \(H_q\) denotes the \(q\)-th Hermite polynomial (see \Cref{sec:wiener-chaos}). We will refer to the functionals $F_{T,q}^{(z)}$ as \emph{Fourier-Hermite functionals}. Our main question is to understand their Gaussian fluctuations. More precisely, whenever
	\(\Var\left(F_{T,q}^{(z)}\right)>0\) for all large \(T\), we ask when
	\[
	\widetilde F_{T,q}^{(z)}
	:=
	\frac{
		F_{T,q}^{(z)}-\E\left[F_{T,q}^{(z)}\right]
	}{
		\sqrt{\Var\left(F_{T,q}^{(z)}\right)}
	},
	\qquad T\rightarrow\infty,
	\]
	converges in distribution to a standard Gaussian, complex-valued when \(z\neq-z\). As a byproduct, in some cases we will also determine the fluctuations for the general functionals $F_T^{(z)}$ under suitable assumptions. \\

    \noindent
	\textbf{Some notation.} 
	The \(q\)-fold convolution of a measure \(G\) will be denoted by \(G^{*q}\).
	For two non-negative quantities \(a_T\) and \(b_T\), we write
	\(a_T\lesssim b_T\) if \(a_T\le Cb_T\) for some constant \(C>0\) and every sufficiently large \(T\). We write \(a_T\asymp b_T\) if both
	\(a_T\lesssim b_T\) and \(b_T\lesssim a_T\) hold. For a possibly
	complex random variable $Z$, throughout the paper
	$\Var(Z):=\E|Z-\E Z|^2$.
	
	Finally, \(N\) denotes a standard real Gaussian random variable, whereas
	\(N_{\mathbb C}\) denotes a standard complex Gaussian random variable, that is,
	\[
	N_{\mathbb C}
	:=
	\frac{N_1+iN_2}{\sqrt2},
	\qquad
	N_1,N_2\overset{\mathrm{i.i.d.}}{\sim}\mathcal N(0,1).
	\]
	We write $N_z=N$ when $z=-z$ and $N_z=N_{\mathbb C}$ otherwise.
	\subsection{Motivations}\label{subsec:motivations}
	
	There are three main motivations for this work.
	
	The first one is to understand whether central limit theorems for functionals of
	stationary Gaussian fields that are usually proved by rather different arguments
	can be explained by a common mechanism. Classical examples are the
	Breuer--Major and Dobrushin--Major theorems, high-frequency results for
	Fourier coefficients of subordinated Gaussian fields on compact groups, and
	more recent results for arithmetic and Euclidean random waves. In these
	problems, the proofs may involve very different ingredients, such as
	integrability of the covariance, special
	properties of the spectral support, or oscillatory integral estimates, see e.g. \cite{BreuerMajor,CammarotaKlurmanWigman,GrottoMissingBerryCLT,MainiNourdinSpectralCLT}.
	One of the purposes of this paper is to show that, at a fixed chaos order,
	these phenomena can be described directly in terms of the way the spectral
	mass is distributed under the additive constraint \eqref{additive constraint}.
	In several examples, this allows one to recover known central limit theorems
	from a direct estimate on convolution powers of the spectral measures.
	
A second motivation comes from Gaussian fields with singular spectral measures. An important example is \textbf{Berry's random wave model}, the stationary isotropic Gaussian field on \(\R^d\) whose spectral measure is uniform on the unit sphere. The model was introduced in \cite{Berry1977} and also arises as a universal local scaling limit for high-frequency random waves on a class of Riemannian manifolds; see e.g. \cite{CanzaniHanin2015}. It has motivated an extensive literature on the geometry and fluctuations of random waves and their functionals; see e.g. \cite{BeliaevHegde,GassPeccatiStecconi,GrottoMissingBerryCLT,GrottoMainiTodino,MainiNourdinSpectralCLT, MainiRossiZheng,NourdinPeccatiRossi}.

More generally, \textbf{stationary Gaussian solutions of constant-coefficient linear PDEs} naturally have spectral measures with support contained in the zero set of the corresponding PDE; see e.g. \cite{EstradeFournier}. In the case of Berry's random wave, or more generally for anisotropic random waves, this zero set is precisely the unit sphere associated with the Helmholtz equation. Other examples arising in applications are \textbf{Gaussian sea waves}, Gaussian fields solving a PDE with zero set given by the Airy surface, see \cite{AzaisLeonOrtegaSeaWaves,EstradeFournier}. This suggests studying Gaussian fields with singular continuous spectral measures supported on lower-dimensional sets in $\R^d$. Such models lie outside the classical short memory setting and may have slowly decaying oscillating covariances. We will see that, for a variety of such Gaussian fields, the geometry of the spectral support controls the small-ball probabilities entering our criterion, yielding a \textbf{connection between spectral geometry and Gaussian fluctuations}. In particular, for both isotropic (i.e. Berry's random wave model) and anisotropic random waves with density bounded below (by a positive constant) and above, we will give a complete description of the CLTs for Fourier-Hermite functionals at every frequency $z$, highlighting that one does not have a CLT only if $d=2$, $q=2$ and $0<\|z\|<2$.

	The third motivation is to understand the role of the probed frequency \(z\). While \(z=0\) often corresponds to "spatial averages" with respect to the window $\psi_T$, the case \(z\neq0\) has also been extensively studied, in particular in connection with finite Fourier transforms, Fourier coefficients and periodograms of stationary processes and random fields; see e.g. \cite{MarinucciPeccatiCompact,PW10,PZ19,RosenblattFourier,FT86,FT87}. These results, obtained for different classes of processes and under rather different assumptions, motivate the search for a common spectral mechanism explaining the dependence of the limiting behavior on \(z\). 
	
	\subsection{Main result and interpretation in regimes \ref{high frequency} and \ref{large domain}}
	\label{subsec:main-nutshell}
	
	We now give an informal description of the main result of the paper. Fix \(q\geq 2\) and a frequency \(z\in\widehat E\). By Nualart-Peccati fourth moment theorem (see \Cref{sec:wiener-chaos}), and standard computations, one can show that the Gaussian convergence is equivalent to the convergence to $0$ of every normalized contraction norm involving the kernel
	\[
	\wpsi_T(\Lambda_1+\dots+\Lambda_q-z)\,.
	\]
	where $\Lambda_1,\dots,\Lambda_q$ i.i.d.$\sim G_T$. Under suitable localization and non-degeneracy assumptions given in \Cref{sec:mainresult}, we will show that the main contribution of these contractions is given by the spectral configurations satisfying the constraint
	\begin{equation}
    \label{additive constraint}
	    	S_q:=\Lambda_1+\cdots+\Lambda_q\in B_{1/T}(z)\,,
	\end{equation}
	where the scale \(1/T\) is the
	localization scale selected by the test functions under the assumptions
	introduced in Section~5, see in particular \Cref{cond:speclocalization}. Our main result shows that Gaussian convergence is
	equivalent to the fact that, under this constraint, no partial
	sum $S_r$, $ 1\leq r\leq q-1$
	can concentrate on a spectral region of the same scale. We call this property
	\emph{spectral delocalization}. More precisely, under \Cref{cond:doubling}, spectral delocalization admits two equivalent formulations.
	
	\medskip
	
	\noindent
	\textbf{Conditional formulation.}
	For some $M>0$ and every  \(r=1,\ldots,q-1\),
	\begin{equation}\label{eq:condform}
		\sup_{\lambda\in\widehat E}
		\mathbb P\left(
		S_{r}\in B_{M/T}(\lambda)
		\,\middle|\,
		S_{q}\in B_{M/T}(z)
		\right)
		\rightarrow0.
	\end{equation}
	
	\medskip
	
	\noindent
	\textbf{Convolution formulation.}
	For some $M>0$ and every \(r=1,\ldots,q-1\),
	\begin{equation}
		\label{eq:convform}
		\sup_{\lambda\in\widehat E}
		\frac{
			G_T^{*r}\!\left(B_{M/T}(\lambda)\right)
			G_T^{*(q-r)}\!\left(B_{M/T}(z-\lambda)\right)
		}{
			G_T^{*q}\!\left(B_{M/T}(z)\right)
		}
		\rightarrow0.
	\end{equation}
	The balls above are taken with respect to a translation invariant metric, that exists under our assumptions, see \Cref{sec:lca-spectral}. The first formulation is easier to interpret and says that, conditionally on the total spectral sum being
	localized around \(z\), no partial sum can concentrate around a
	single frequency. The second formulation is easier to use in applications, since it is expressed in terms of small ball probabilities with respect to the convolution powers of the varying spectral measure.
	
	Our main result, \Cref{thm:main}, can therefore be stated informally as follows.
	
	\begin{theorem}[\Cref{thm:main}-(iii)]
		\label{thm:informal-main}
		Under the assumptions in \Cref{thm:main}-(iii),
		\[
		\widetilde F_{T,q}^{(z)}
		\overset{\rm law}{\rightarrow} N_z
		\qquad\Leftrightarrow\qquad
		\text{spectral delocalization}\,.
		\]
	\end{theorem}
    \vspace{2mm}
	\noindent
	We now proceed formally stating and discussing our main result in the settings \ref{high frequency} and \ref{large domain}.
	
	\begin{theorem}\label{thm:mainTd}
		Consider regime~\ref{high frequency}. Fix $q\ge2$. Then
		\[
		\Var\left(F_{T,q}^{(z)}\right)
		=
		q!\,G_T^{*q}(\{z\})>0
		\quad\;\Leftrightarrow\;\quad
		z\in\supp\left(G_T^{*q}\right).
		\]
		Moreover, if \(z\in\supp(G_T^{*q})\) for every \(T\), then $\widetilde F_{T,q}^{(z)}\overset{\rm law}{\rightarrow}N_z$ if and only if for every \(r=1,\ldots,q-1\),
			\begin{equation}
				\label{eq:SACtoro}
				\frac{
					\sup_{\lambda\in\widehat E}
					G_T^{*r}(\{\lambda\})
					G_T^{*(q-r)}(\{z-\lambda\})
				}{
					G_T^{*q}(\{z\})
				}
				\rightarrow0.
			\end{equation}
	\end{theorem}
	
	Condition~\eqref{eq:SACtoro} is a weighted combinatorial condition; for trigonometric polynomials and arithmetic random waves it becomes purely combinatorial, see \Cref{sec:highfrequencyapplications}. As we will see, the thickness of the growing spectral window for trigonometric polynomials spreads the possible configurations, helps spectral delocalization and yields a CLT at every frequency. Conversely, the thin spectral windows of arithmetic random waves yield no CLT in many cases. 
	
	\begin{remark}
		In \cite[Theorem~2]{MarinucciPeccatiCompact}, Gaussian convergence
		of high-frequency Fourier coefficients on compact Abelian groups
		is characterized through a random-walk bridge condition.
		Conditionally on a \(q\)-step random walk on the dual group ending
		at the probed frequency, no intermediate position can retain
		positive probability uniformly over all positions. In particular,
		no ``privileged path'' survives asymptotically.
		
		This is the same discrete spectral delocalization condition as
		\eqref{eq:SACtoro}, but in a different asymptotic regime:
		in \cite{MarinucciPeccatiCompact} the spectral measure is fixed
		and the frequency tends to infinity, whereas here the frequency
		is fixed and \(G_T\) varies.
		Furthermore, \Cref{thm:mainRdZd} extends this principle to
		growing-domain functionals, replacing point masses by small-ball
		probabilities of convolution powers, also for singular continuous
		spectral measures.
	\end{remark}
	\vspace{2mm}
    
	In regime~\ref{large domain}, the criterion is a small-ball-probabilities
	condition. Here we need additional assumptions.
	\begin{condition}
		\label{cond:schwartz}
		The real-valued function \(\psi\) belongs to \(\mathcal S(\R^d)\) and
		\(\wpsi(0)>0\).
	\end{condition}
	\begin{condition}
		\label{cond:doublingintro}
		Fix $q\ge2$ and $z\in\supp(G^{*q})$. There exist $\varepsilon_0>0$ and
		$C_D>1$ such that
		\[
		G^{*q}(B_{2\varepsilon}(z))
		\le C_D\,G^{*q}(B_\varepsilon(z)),
		\qquad 0<\varepsilon<\varepsilon_0.
		\]
	\end{condition}
	\begin{remark}\label{rem:alwaysholdscases}
		\Cref{cond:doublingintro} holds whenever $G^{*q}$ has a continuous positive density at $z$, or more generally if $G^{*q}(B_\eps(z))\asymp\varepsilon^s L(1/\eps)^\ell$ for positive exponents $s,\ell>0$ and $L$ slowly varying at $\infty$. Moreover, it is always automatically satisfied if $z=0$ and $q$ is even, see \Cref{superem}. Indeed, in this case we have the stronger ball growth and rectangular growth conditions needed in \Cref{rem:indicator-windows}  to prove the analogous criterion \eqref{eq:SACgrowing} when $\psi=\mathbf{1}_D$, $D\subseteq \R^d$, see \Cref{prop:indicator-windows}.
	\end{remark}

	The Schwartz assumption is natural (see e.g. \cite{BeliaevHegde,GassPeccatiStecconi}) and allows to simplify the statement and the proof of \Cref{thm:mainRdZd}. Less regular test functions like indicators of cubes/balls can be treated under ball growth or rectangular growth assumptions
	on \(G^{*q}\), adapted to the Fourier decay of $\psi$; see
	\Cref{rem:indicator-windows} and \Cref{prop:indicator-windows}. On \(\widehat E=\mathbb R^d\) and \(\widehat E=\mathbb T^d\) we use, respectively, the standard Euclidean and toral distances/balls.
	
	\begin{theorem}\label{thm:mainRdZd}
		Consider regime~\ref{large domain}, fix $q\ge2$ and assume
		\Cref{cond:schwartz}. If
		\(z\notin\supp(G^{*q})\), then, $\forall$ \(N>0\),
		\[
		\Var\left(F_{T,q}^{(z)}\right)=O(T^{-N}).
		\]
		Conversely, if \(z\in\supp(G^{*q})\) and \Cref{cond:schwartz}, \Cref{cond:doublingintro} hold, then
		\[
		\Var\left(F_{T,q}^{(z)}\right)
		\asymp q!\,T^{2d}G^{*q}(B_{1/T}(z)),
		\]
		and we have that $\widetilde F_{T,q}^{(z)}\overset{\rm law}{\rightarrow}N_z$ if and only if for every \(r=1,\ldots,q-1\),
			\begin{equation}
				\label{eq:SACgrowing}
				\lim_{\varepsilon\rightarrow0}
				\sup_{\lambda\in\widehat E}
				\frac{G^{*r}(B_\varepsilon(\lambda))
					G^{*(q-r)}(B_\varepsilon(z-\lambda))}
				{G^{*q}(B_\varepsilon(z))}=0\,.
			\end{equation}
	\end{theorem}
	Combining the previous result for $q$ fixed with reduction lemmas such as \Cref{lem:finite-chaos-reduction}, one can also deduce CLTs for general functionals.
	Moreover, as a byproduct of \Cref{thm:mainRdZd}, we get the following corollary.
	
	\begin{corollary}[Equivalent spectral measures]\label{cor:equivalentspectralmeasures}
		\label{cor:equivalent-spectral-measures}
		Let \(G\) and \(G'\) be symmetric probability measures on \(\widehat E\) such that $cG'\le G\le CG'$
		for some \(0<c\le C<\infty\). Then, for every \(q\ge2\) and \(z\in\widehat E\), spectral delocalization holds for \(G\) if and only if it holds for \(G'\). Moreover, the local doubling condition at \(z\) holds for \(G^{*q}\) if and only if it holds for \(G'^{*q}\). Consequently, under the assumptions of \Cref{thm:mainRdZd}, the corresponding Fourier-Hermite functionals either both have Gaussian fluctuations or neither does.
	\end{corollary}
	One consequence, developed in \Cref{subsec:short-memory}, is that the classical Breuer-Major condition implies Gaussian fluctuations at every fixed frequency with non-zero limiting variance. Outside this regime, instead, the validity of the central limit theorem may depend on the probed frequency, see \Cref{sec:large-domain-applications}.
	
	Limit theorems in the regime \ref{large domain} under spectral assumptions have a long history; see, among others, \cite{AvramLeonenkoSakhno,IvanovLeonenkoRuizMedinaSavich,NualartZhengAveraging,RosenblattFourier,MainiNourdinSpectralCLT}. The equivalence between Gaussian convergence and the vanishing of contraction norms is given by the Nualart--Peccati Fourth Moment Theorem \cite{NualartPeccati} and its complex counterpart; see e.g. \cite{MarinucciPeccatiCompact} and \Cref{sec:wiener-chaos}. 
	
	The novelty of \Cref{thm:mainRdZd} is that, for Fourier--Hermite functionals, these contraction conditions are characterized by the explicit small-ball condition \eqref{eq:SACgrowing}, involving only convolution powers of the spectral measure. This characterization will allow us to directly recover known results by spectral estimates, as well as to obtain new applications.
	
	\subsection{Organization of the paper}
	
	In \Cref{sec:highfrequencyapplications,sec:large-domain-applications} we compare the criterion with previous results and discuss consequences and applications, respectively in the regimes \ref{high frequency} and \ref{large domain}. The preliminaries about topology and harmonic analysis on groups, isonormal representations and fourth moment theorems are collected in \Cref{sec:preliminaries}. The general version of the criterion, \Cref{thm:main}, is stated in \Cref{sec:mainresult}. Its proof, together with the proofs of \Cref{thm:mainTd} and \Cref{thm:mainRdZd}, is given in \Cref{sec:proofs}, where we also prove \Cref{prop:indicator-windows}, a version of \Cref{thm:mainRdZd} for indicator functions. Finally, in \Cref{sec:additional proofs} and \Cref{sec:ehresmann} we collect the proofs of some auxiliary technical results needed in \Cref{sec:large-domain-applications}.

	\section{Consequences in the regime~\ref{high frequency}: previous results and applications}\label{sec:highfrequencyapplications}
	We consider two
	standard toral examples: random trigonometric polynomials (or more generally band-limited fields) and arithmetic
	random waves. In both these cases the spectral measure $G_T$ is uniformly distributed on a spectral window $\Lambda_T\subseteq \Z^d$, that is
	\[
	G_T:=\frac{1}{\mathcal N_T}
	\sum_{\lambda\in\Lambda_T}\delta_\lambda\,,\qquad\mathcal{N}_T:=\left|\Lambda_T\right|.
	\]
	Thus, the spectral delocalization condition becomes purely combinatorial. Indeed, if we define
	\[
	D_{T,m}(w):=\left|
	\left\{
	(\lambda_1,\ldots,\lambda_m)\in\Lambda_T^m:
	\lambda_1+\cdots+\lambda_m=w
	\right\}
	\right|\,,
	\]
	then we have, $\forall $ $m\ge 1$ and $w\in\Z^d$,
	\[
	G_T^{*m}(w):=D_{T,m}(w)/\mathcal{N}_T^m
	\]
	and so by \Cref{thm:mainTd}, along every sequence such that \(D_{T,q}(z)>0\), we get
	\begin{equation}
		\label{eq:criterioARW}
		\widetilde F_{T,q}^{(z)}\overset{\rm law}{\rightarrow}N_z\quad \;\Leftrightarrow\; \quad \frac{\sup_{w\in\Z^d}\,D_{T,r}(w)D_{T,q-r}(z-w)}{D_{T,q}(z)}\rightarrow 0 \quad \forall \,r=1,\dots,q-1\,.
	\end{equation}
	
	\subsection{Band-limited toral fields and random trigonometric polynomials}
	\label{subsec:band-limited-toral-fields}
	
	Band-limited Gaussian fields on the torus correspond to the case where \(K\subseteq\mathbb R^d\) is a symmetric convex body and the spectral window is a dilation of $K$ in $\Z^d$, namely
	\[
	\Lambda_T:=TK\cap\mathbb Z^d,
	\]
	The classical random trigonometric polynomial is obtained in dimension
	\(d=1\), with \(K=[-1,1]\), after removing $0$. Namely, for
	\(T\in\mathbb N\),
	\[
	\Lambda_T=\{-T,\ldots,-1,1,\ldots,T\},
	\]
	and the corresponding field can be written as
	\[
	X_T(t)
	=
	\frac1{\sqrt T}
	\sum_{k=1}^T
	\bigl(a_k\cos(kt)+b_k\sin(kt)\bigr),
	\qquad t\in\mathbb T,
	\]
	where \((a_k)_{k\ge1}\) and \((b_k)_{k\ge1}\) are independent standard
	Gaussian random variables.
	Random trigonometric polynomials have a long literature on zeros, level
	crossings, universality and dependence; see, among others,
	\cite{AngstDalmaoPoly,AngstNguyenPoly,AngstPolySalem}.
	Using our criterion, we can prove the following result.
	\begin{proposition}
		\label{prop:band-limited-toral-fields}
		Fix \(q\ge2\) and \(z\in\mathbb Z^d\). Assume $\Lambda_T=TK\cap\Z^d$, or $\Lambda_T=TK\cap\Z^d\setminus\{0\}$. Then for $T$ large enough $D_{T,q}(z)>0$ and as $T\rightarrow\infty$
		\[
		\sup_{\alpha\in\mathbb Z^d}
		\frac{
			D_{T,r}(\alpha)D_{T,q-r}(z-\alpha)
		}{
			D_{T,q}(z)
		}
		=
		O(T^{-d})\,,\qquad\forall\,\,r=1,\ldots,q-1\,.
		\]
		Moreover, we have $\widetilde F_{T,q}^{(z)}
		\overset{\mathrm{law}}{\rightarrow}
		N_z$ and $\Var\left(F_{T,q}^{(z)}\right)=q!D_{T,q}(z)/\mathcal{N}_T^q\asymp T^{-d}$.
	\end{proposition}
	
	\begin{proof}
		Since \(K\) is a symmetric convex body, we have $\mathcal N_T\asymp T^d$.
		Moreover, for every \(m\ge1\) and \(w\in\mathbb Z^d\), the choices for the first \(m-1\)
		elements contributing to \(D_{T,m}(w)\) are bounded by $\mathcal{N}_T^{m-1}$, but the last choice must be $\lambda_m=w-\sum_{i=1}^{m-1}\lambda_i$,
		thus
		\[
		\sup_\alpha D_{T,m}(\alpha)\le \mathcal N_T^{m-1}\quad \;\Rightarrow\;\quad \sup_\alpha D_{T,r}(\alpha)D_{T,q-r}(z-\alpha)
		\le
		\mathcal N_T^{q-2}\asymp T^{d(q-2)}
		\]
		To lower bound the denominator, choose \(\rho>0\) such that $B_\rho(0)\subseteq K$ and consider the set
		\[
		A=\left\{(\lambda_1,\dots,\lambda_{q-1},z-\lambda_1-\cdots-\lambda_{q-1})\,:\,\lambda_i\in B_{\frac{\rho T}{2q}}(0)\cap \Z^d\setminus\{0\}\,\,,\,\forall\,i\le q-1\right\}\,.
		\]
		Then, for $T$ large enough we have $\lambda_q\in B_{\|z\|+\rho T/2}(0)\subseteq B_{\rho T}(0)\subseteq TK$. When $\Lambda_T=TK\cap\Z^d$ we get 
		\[
		D_{T,q}(z)\ge |A|
		\ge
		\left|B_{\frac{\rho T}{2q}}(0)\cap \Z^d\setminus\{0\}\right|^{q-1}\asymp T^{d(q-1)}\,.
		\]
        When $\Lambda_T=TK\cap\Z^d\setminus\{0\}$ we conclude observing that $z=\sum_{i=1}^{q-1}\lambda_i$, i.e. $\lambda_q=0\notin \Lambda_T$ must be excluded. Thus, since $|A\cap\{\lambda_q=0\}|\le|\{(\lambda_1,\dots,\lambda_{q-2},-\sum_{i=1}^{q-2}\lambda_i)\,:\,\lambda_i\in B_{\rho T/2q}\}|\lesssim T^{d(q-2)}$, we have
        \[
        D_{T,q}(z)\ge
\left|B_{\frac{\rho T}{2q}}(0)\cap \Z^d\setminus\{0\}\right|^{q-1}
-\,O(T^{d(q-2)})\,
\asymp\, T^{d(q-1)}\quad\,.
        \]
	\end{proof}
	As far as we know, \Cref{prop:band-limited-toral-fields} is a new result at this level of generality. By standard arguments, see e.g.
	\cite[Theorem~6.3.1]{NourdinPeccatiBook}, the conclusion can be extended
	to general functionals \(F_T^{(z)}\), provided the corresponding variance
	and chaotic-tail conditions hold. We do not pursue this extension here.
	
	Analogous results were proved in the case $d=1$ and $z=0$ in \cite{AngstPolyFluctuations} for random trigonometric polynomials with i.i.d. coefficients.

	\subsection{Arithmetic random waves}
	\label{subsec:arithmetic-random-waves}
	Fix \(d\ge2\). Arithmetic random waves (or toral random Laplace eigenfunctions) have been widely studied, see, among others,
	\cite{CammarotaKlurmanWigman,CammarotaMaffucciMarinucciRossi,KKW,PeccatiRossiARW}. For an admissible  \(T\in\N\) (i.e. such that $\mathcal N_T>0$), the spectral window in this case is
	\[
	\Lambda_T:=\{\lambda\in\Z^d:\|\lambda\|^2=T\},
	\qquad \mathcal N_T:=\left|\Lambda_T\right|\,.
	\]
	Arithmetic random waves are often expressed as
	\[
	X_T(x):=\frac1{\sqrt{\mathcal N_T}}
	\sum_{\lambda\in\Lambda_T}a_\lambda e^{i\langle x,\lambda\rangle},
	\qquad x\in\T^d,
	\]
	where $a_\lambda$ are standard complex Gaussian, \(\overline{a_\lambda}=a_{-\lambda}\) and $a_\lambda$, $a_\theta$ are independent if $\lambda\neq\pm\theta$.
	
	Along every sequence \(T_j\) such that
	\(\mathcal N_{T_j}\to\infty\) and \(D_{T_j,q}(z)>0\), the combinatorial quantities of interest $D_{T,m}(w)$, $m\ge1$, $w\in\widehat{E}$, are known as \textbf{spectral correlations} in the literature of arithmetic random waves and are harder to estimate. Indeed, the problem becomes a number-theoretic question, see e.g. \cite{BombieriBourgain}.
	
	If \(q=2\), the supremum in \eqref{eq:criterioARW} equals $1$ whenever
	\(D_{T,2}(z)>0\) and so we get the following.
	\begin{proposition}
		Let $q=2$. Then, if $\mathcal{N}_{T_j}\rightarrow\infty$ and $D_{T_j,2}(z)>0$, we have
		\[
		\widetilde F_{T_j,2}^{(z)}\overset{\rm law}{\rightarrow}N_z
		\quad\;\Leftrightarrow\;\quad
		D_{T_j,2}(z)\rightarrow\infty.
		\]
		In particular, since \(D_{T,2}(0)=\mathcal N_T\), $F_{T,2}^{(0)}$ satisfies a CLT. If $d=2$, for every
		\(z\ne0\) the intersection of the circle \(\Lambda_T\) with its
		translate \(z-\Lambda_T\) has at most two points, so
		\(D_{T,2}(z)\le2\) and $F_{T,2}^{(z)}$ does not satisfy a CLT.
	\end{proposition}
	For \(q=2\), \(d=2\) and \(z=0\), the central limit theorem is classical:
	\(F_{T,2}^{(0)}\) is a normalized sum of independent centered squared
	Gaussian coefficients; see e.g.
	\cite{CammarotaMarinucciRossi}.
	The proposition above gives the corresponding characterization in arbitrary
	dimension and at every fixed frequency. In particular, when $d=2$ the frequency $z=0$ is special and is the only one for which the functional $F_{T,2}^{(z)}$ satisfies a CLT. As we will see, this "no CLT away from $0$" phenomenon for second chaos functionals will be observed analogously for planar random waves in \Cref{prop:anisotropic-Berry-all-frequencies}.
	
	When $q>2$ things are more complicated. However, for planar arithmetic random waves at frequency $0$, i.e. when $d=2$ and $z=0$, sharp estimates for the zero-frequency
	spectral correlations are available. 
    
    For $d=2$ and $q$ odd, we have $D_{T,q}(0)=0$.
	Therefore, from now on let us fix $(d,z)=(2,0)$ and consider $q=2k$ even, $k\ge1$.
	The diagonal configurations cancelling in pairs already give the lower bound
	\[
	D_{T,2k}(0)\ge \mathcal N_{T}^k;
	\]
	see e.g. \cite[(1.26)-(1.27)]
	{CammarotaKlurmanWigman}. On the other hand,
	\cite[(1.31)]
	{CammarotaKlurmanWigman}, together with the inclusion of spectral
	correlations into spectral semi-correlations, implies that there exists a subsequence $T_{j}$ of relative density one among the admissible $T\in \N$ such that \(\mathcal N_{T_j}\to\infty\) and,
	simultaneously for every \(k\ge1  \),
	\begin{equation}
		\label{eq:optimal-ARW-correlations}
		D_{T_j,2k}(0)\asymp \mathcal N_{T_j}^k,
		\qquad \text{as } j\rightarrow\infty.
	\end{equation}
	
	Combining all these estimates and using our criterion, we immediately get the following.
	
	\begin{proposition}
		\label{prop:ARW-zero-frequency}
		Let \(d=2\), $q\ge4$ even and $z=0$. There exists a a subsequence $T_{j}$ of relative density one among the admissible $T\in \N$ such that \(\mathcal N_{T_j}\to\infty\) and $\widetilde F_{T,q}^{(z)}$ does not satisfy a CLT.
	\end{proposition}
	
	\begin{proof}
		Taking a sequence as in \eqref{eq:optimal-ARW-correlations}, the ratio \eqref{eq:criterioARW} for $r=2$ even is greater or equal than $D_{T,2}(0)D_{T,q-2}(0)/D_{T,q}(0)\asymp1$. We conclude by \eqref{eq:criterioARW} (i.e. \eqref{eq:SACtoro} in \Cref{thm:mainTd}).
	\end{proof}

	For \(q=4\), related non-Gaussian fourth-chaos phenomena for arithmetic
	random waves are well known; see, for instance,
	\cite{MPRW,PeccatiRossiARW}. To the best of our knowledge, the result above for arbitrary even order
	\(q\ge6\) has not been stated previously.
	
	\section{Consequences in the regime \ref{large domain}: previous results and applications}
	\label{sec:large-domain-applications}
	
	Throughout this section, \(E=\R^d\) or \(E=\Z^d\), and \(X\) is a fixed
	stationary centered Gaussian field with continuous covariance \(C\) and spectral measure
	\(G\). Unless otherwise stated, we assume the doubling \Cref{cond:doublingintro} for $G^{*q}$ at $z$ and we take
	\(\psi_T(x)=\psi(x/T)\), where \(\psi\in\mathcal S(\R^d)\) is real-valued and
	\(\wpsi(0)>0\). Indicator functions can be treated under ball growth and rectangular growth assumptions
	of \Cref{rem:indicator-windows} and \Cref{prop:indicator-windows}. When $q$ is even and $z=0$, all these assumptions on the growth of $G$ are satisfied, see \Cref{rem:alwaysholdscases} and \Cref{superem}.
	
	We shall mainly work in a fixed Wiener chaos, considering $F_{T,q}^{(z)}$. However, when considering a general
	nonlinear functional $F_T^{(z)}$, we assume that
	\(\varphi\in L^2(\gamma)\) and exploit its orthogonal decomposition in terms of the orthogonal basis of Hermite polynomials, see e.g. \cite{NourdinPeccatiBook}
\begin{equation}
		\label{eq:phiexp}
	\varphi(x)-\E[\varphi(N)]
	=
	\sum_{q\ge1}c_qH_q(x),
	\qquad
	\Var(\varphi(N))=\sum_{q\ge1}q!c_q^2<\infty\,,
\end{equation}
	where $c_q:=\frac{1}{q!}\E\left[H_q(N)\varphi(N)\right]$. If \(\varphi\) is non-constant, its Hermite rank $Q:=\min\{q\ge1:c_q\ne0\}<\infty$ is well defined.
	Moreover, we have a chaotic decomposition, i.e. an orthogonal decomposition in $L^2(\Omega)$ for the functional $F_T^{(z)}$ (see e.g. \cite{NourdinPeccatiBook} or \Cref{sec:wiener-chaos})
	\begin{equation}
	    F_T^{(z)}-\E\left[F_{T}^{(z)}\right]	=
	\sum_{q\ge Q}c_qF_{T,q}^{(z)},
	\qquad
	\Var\left(F_T^{(z)}\right)=\sum_{q\ge Q}c_q^2 \Var\left(F_{T,q}^{(z)}\right)\,.\label{ortovar}
	\end{equation}
	We refer to $c_qF_{T,q}^{(z)}$ as the $q$th chaotic component of $F_{T}^{(z)}$, see \Cref{sec:wiener-chaos} and \Cref{sec:mainresult}. By Hermite isometry, see \eqref{eq:variance-kernel}-\eqref{eq:pseudovariance-kernel}, one has the following expressions for the variance and pseudovariance of a chaotic component
    \[
    \Var\left(F_{T,q}^{(z)}\right)=q!\int_{\wE} \left|\wpsi_T(\lambda-z)\right|^2 G^{*q}(d\lambda)\,\,,\qquad\E\left[\left(F_{T,q}^{(z)}\right)^2\right]=q!\int_{\wE} \wpsi_T(\lambda-z) \overline{\wpsi_T(\lambda+z)} G^{*q}(d\lambda)\,,
    \]
    where $\wpsi_T$ is the Fourier transform of $\psi_T:E\rightarrow\R$, as defined in \eqref{eq:Fourierdef}.
    On $E=\R^d,\Z^d$ we consider respectively $\mu$ the Lebesgue and counting measure. As a consequence, the only Haar measure $\wmu$ on $\wE=\R^d,\T^d$ satisfying Plancherel identity \eqref{eq:plancherel} is given by
    \begin{equation}
        \wmu(d\lambda)=d\lambda/(2\pi)^d\,.\label{eq:muhat}
    \end{equation}
    Whenever $G^{*q}$ admits a density $g_q$, this density is taken with respect to $\wmu$, so that
    $$G^{*q}(d\lambda)=g_q(\lambda)\wmu(d\lambda)=g_q(\lambda)d\lambda/(2\pi)^d\,.$$
    Moreover, when $C^q\in L^1$, since $G^{*q}$ is its spectral measure, Fourier inversion and symmetry of $G^{*q}$ give
    \begin{equation}
        \label{eq:fourierinversion}
        g_q(z)=\int_E C(x)^q e^{-i\langle x,z\rangle}\,\mu(dx)=\int_E C(x)^q e^{i\langle x,z\rangle}\,\mu(dx)\,.
    \end{equation}
\subsection{Short memory and Breuer-Major CLT at every frequency}
\label{subsec:short-memory}

For nonlinear functionals with Hermite rank $Q$, the condition
$C\in L^Q(E)$ is the classical short-memory assumption of the
Breuer--Major theorem \cite{BreuerMajor}, which gives Gaussian fluctuations
for spatial averages, corresponding here to $z=0$. Several extensions
have been obtained, see e.g. \cite{NualartTilva,NualartZhengAveraging}.
Throughout \Cref{sec:large-domain-applications}, we sometimes also consider the indicator windows
in \Cref{prop:indicator-windows}:
    \begin{enumerate}[label=\textnormal{(\roman*)}]
			\item \label{psiquadro}If \(E=\mathbb Z^d, \R^d\), we consider also $\psi=\mathbf1_{Q_{\mathbf{r}}(a)}$, $Q_{\mathbf{r}}(a)$ rectangle,
			$a\in E$ and \(\mathbf r=(r_1,\ldots,r_d)\in(0,\infty)^d\);
			
			\item \label{psiliscio}If \(E=\mathbb R^d\), we consider also \(\psi=\mathbf1_D\), where \(D\) is a
			convex body with nonempty interior and \(C^\infty\) boundary with non-vanishing
			Gaussian curvature (e.g. Euclidean balls).
		\end{enumerate} 

The goal of this subsection is to show a Breuer-Major theorem at every fixed frequency $z$, i.e. to show that in the short memory regime we always have a CLT independently of the choice of the probed frequency.

In order to have the exact asymptotics of the variance as in the classical statement of the Breuer-Major theorem, we start determining the asymptotic variance of $F_{T,q}^{(z)}$ when $G^{*q}$ has bounded density, continuous at $z$. In particular, by Fourier inversion this is always the case when $C^q\in L^1(E)$. The proof of the following lemma is postponed to \Cref{sec:additional proofs}.

\begin{lemma}\label{lem:variance-density}
Let $q\ge1$ and suppose that $G^{*q}$ has a bounded density
$g_q$ with respect to $\widehat \mu$, continuous at $z$. If $\psi$ satisfies
\Cref{cond:schwartz}, or \ref{psiquadro}, or \ref{psiliscio}, then $a_T:=T^{-d}\|\psi_T\|^2_{L^2(E,\mu)}\rightarrow\|\psi\|_2^2:=\|\psi\|_{L^2(\R^d)}^2$ and
\begin{equation}
    T^{-d}\Var\left(F_{T,q}^{(z)}\right)
\rightarrow q!\,\|\psi\|_2^2\,g_q(z)\,\,,\qquad\qquad T^{-d}\Var\left(F_{T,q}^{(z)}\right)\le q!\,\left\|g_q\right\|_\infty\,a_T\,\,.\label{varasymp}
\end{equation}
Moreover, if $z\ne-z$, then
\[
T^{-d}\,\E\left[\left(F_{T,q}^{(z)}\right)^2\right]\rightarrow0\quad\;\Rightarrow\;\quad T^{-d}\,{\rm Cov}\left({\rm Re}\left(F_{T,q}^{(z)}\right),{\rm Im}\left(F_{T,q}^{(z)}\right)\right)\rightarrow\,0\,.
\]
\end{lemma}

We now show that \textbf{short-memory conditions yield spectral delocalization
at every frequency}. By Fourier inversion, if $C^q\in L^1$ then $G^{*q}$ has continuous bounded density given by \eqref{eq:fourierinversion}.

\begin{proposition}\label{prop:short-memory-criterion}
Let $q\ge2$. Assume $C\in L^q(E)$. Then, for every
$s=1,\ldots,q-1$,
\begin{equation}\label{topboundBM}
\sup_{\lambda\in\widehat E}
G^{*s}(B_\varepsilon(\lambda))
=o(\varepsilon^{ds/q}),\qquad \varepsilon\to0\,,
\end{equation}
If $\psi$ is as in \Cref{lem:variance-density} and $g_q(z)>0$, then $G^{*q}(B_\eps(z))\asymp\eps^d$, \eqref{eq:SACgrowing} holds, and
\[
\widetilde F_{T,q}^{(z)}\overset{\rm law}{\rightarrow}N_z,
\qquad
\Var\left(F_{T,q}^{(z)}\right)\sim q!T^d\|\psi\|_{L^2(\R^d)}^2g_q(z).
\]
\end{proposition}
\begin{proof}[Proof]
Let $\chi\in C_c^\infty(\R^d)$ be even, non-negative and equal to one
on the unit ball, and set
$\chi_\varepsilon(\lambda)=\chi(\lambda/\varepsilon)$,
$\lambda\in\widehat E$. If $E=\Z^d$, take $\varepsilon>0$ small so that
$\operatorname{supp}\chi_\varepsilon\subset\T^d=[-\pi,\pi)^d$.
Fix a compact set $K\subseteq\R^d$. Fourier inversion and
H\"older's inequality give, uniformly in $\lambda$,
\[
\begin{aligned}
G^{*s}(B_\varepsilon(\lambda))
\le\int_{\widehat E}\chi_\varepsilon(\lambda-\alpha)G^{*s}(d\alpha)\lesssim
\left\|C^s\mathbf1_K\right\|_{L^1(E)}\left\|\widehat\chi_\varepsilon\right\|_{L^\infty(E)}
+\left\|C^s\mathbf1_{K^c}\right\|_{L^{q/s}(E)}
\left\|\widehat\chi_\varepsilon\right\|_{L^{q/(q-s)}(E)}.
\end{aligned}
\]
Now $\|\widehat\chi_\varepsilon\|_{L^\infty(E)}\asymp\varepsilon^d$
and $\|\widehat\chi_\varepsilon\|_{L^{q/(q-s)}(E)}
\asymp\varepsilon^{ds/q}$. Indeed, in the Euclidean case
$\widehat\chi_\varepsilon(x)=\varepsilon^d\widehat\chi(\varepsilon x)$
and the estimate follows. If $E=\Z^d$,
then for $\eps$ small enough\footnote{Here the $\widehat\cdot$ to the left is on $\Z^d$, the $\widehat\cdot$ to the right is on $\R^d$.}
$\widehat\chi_\varepsilon(n)=\varepsilon^d\widehat\chi(\varepsilon n)$
for every $n\in\Z^d$, and the second estimate follows from
\Cref{lem:Riemann-sum} applied to
$|\widehat\chi|^{q/(q-s)}$.
Letting first $\varepsilon\to0$ and then $K\uparrow E$ proves
\eqref{topboundBM}.

The variance asymptotics follow by \Cref{lem:variance-density}. If $g_q(z)>0$, then $G^{*q}(B_\eps(z))\asymp\eps^d$ by continuity of $g_q$. By boundedness of $g_q$ we also have rectangular growth, ball growth of \Cref{prop:indicator-windows} and in particular the doubling  \Cref{cond:doublingintro} of $G^{*q}$ at $z$ (see \Cref{rem:alwaysholdscases}).
Thus, we can apply \Cref{thm:mainRdZd} and \Cref{prop:indicator-windows}, and
\eqref{topboundBM} yields \eqref{eq:SACgrowing}.
\end{proof}
A consequence of \Cref{prop:short-memory-criterion} is that if $C\in L^m$, then $\forall$ $q\ge\max\{m,2\}$ $C\in L^q$ $\implies \widetilde F_{T,q}^{(z)}\overset{\rm law}{\rightarrow}N_z$.
Exploiting this fact and using Peccati-Tudor \eqref{eq:peccatitudor}, we prove the following finite-chaos reduction lemma, which will also yield the desired Breuer–Major theorem for Fourier functionals. We will also use this result in the following subsections to extend Gaussian convergence from finite to infinite sums of chaotic components. Its proof is postponed to \Cref{sec:additional proofs}.

\begin{lemma}\label{lem:finite-chaos-reduction}
Let $\psi$ be as in \Cref{lem:variance-density}. Assume $C\in L^m(E)$ for some $m\ge1$ and
$\sigma_T^2:=\Var\left(F_T^{(z)}\right)\gtrsim T^d$. If, for every
$1\le q<m$ with $c_q\ne0$, either
$\Var\left(F_{T,q}^{(z)}\right)=o(\sigma_T^2)$ or
$\widetilde F_{T,q}^{(z)}\overset{\rm law}{\rightarrow}N_z$,
then
$\widetilde F_T^{(z)}\overset{\rm law}{\rightarrow}N_z$.
\end{lemma}

As a direct consequence, we obtain a Breuer--Major
theorem at every fixed frequency.

\begin{corollary}[Breuer--Major theorem for Fourier functionals]
\label{thm}
Let $\varphi\in L^2(\R,\gamma)$ have Hermite rank
$Q\ge1$. Assume $\psi$ as in \Cref{lem:variance-density} and $C\in L^Q(E)$, with $g_q$ density of $G^{*q}$ with respect to $\wmu$. Set
\[
\sigma_z^2:=\|\psi\|_2^2
\sum_{q\ge Q}q!c_q^2g_q(z).
\]
Then $T^{-d/2}\left(F_T^{(z)}-\E\left[F_T^{(z)}\right]\right)
\overset{\rm law}{\rightarrow}\sigma_zN_z$.
In particular, if $\sigma_z^2>0$, then
$\widetilde F_T^{(z)}
\overset{\rm law}{\rightarrow}N_z$.
\end{corollary}

\begin{proof}
Since $|C|\le1$, for $q\ge Q$ we have
$\|g_q\|_\infty\le\|C^Q\|_{L^1(E)}$.
By \eqref{varasymp} and \eqref{ortovar}, we have
\[
T^{-d}\Var\left(F_T^{(z)}\right)
\rightarrow
\|\psi\|_2^2\sum_{q\ge Q}q!c_q^2g_q(z)
=\sigma_z^2.
\]
If $\sigma_z^2=0$, $T^{-d/2}\left(F_T^{(z)}-\E\left[F_T^{(z)}\right]\right)
\overset{L^2}{\rightarrow}0$.
If $\sigma_z^2>0$ we conclude by \Cref{lem:finite-chaos-reduction}
with $m=Q$.
\end{proof}

For $z=0$, the preceding corollary is the classical Breuer--Major theorem.
It shows that the same condition $C\in L^Q(E)$ yields Gaussian
fluctuations at every fixed frequency $z$ with $\sigma_z^2>0$.
Our criterion is applied in \Cref{prop:short-memory-criterion} and gives a spectral interpretation:
in the short-memory regime, the frequency affects the limiting
variance and whether the limiting Gaussian is real or complex,
but not Gaussian convergence when the limiting variance is non-zero.
	
	\subsection{(Cyclical) long memory and the role of the frequency}
	\label{subsec:long-memory-frequency}
	
	Outside the short memory regime, the limiting behavior may depend on the
	probed frequency. The most classical long memory covariance functions are \textbf{radial long memory regularly varying} covariance functions, namely
	\begin{equation*}
		C(x)= L\left(\|x\|\right)\|x\|^{-\alpha}\,,\quad \|x\|\rightarrow\infty\,,\quad\alpha <d/Q\,,
	\end{equation*}
	with $L:[0,\infty)\rightarrow (0,\infty)$ slowly varying at $\infty$. Since $\alpha< d/Q$ we have long memory, i.e. $C\notin L^Q(E)$, and classical Dobrushin-Major-Taqqu limit theorems \cite{DobrushinMajor,Taqqu1979} provide non-Gaussian limits (the so-called Hermite distributions) for spatial averages (i.e. when $z=0$) when $\alpha<d/Q$ and $Q\ge2$ \footnote{If $Q\alpha>d$ we are in the short memory case. The case $Q\alpha=d$ is critical and one usually has a CLT. When $Q=1$, the first chaos dominates and so we always have a CLT.}. For these covariances, Abelian-Tauberian theorems \cite{LeonenkoOlenko} provide a control on the spectral small ball probabilities. More precisely, for $r\le Q$ we have (here $c(r\alpha,d)>0$ is explicit, see \cite[Theorem 3, p.723]{LeonenkoOlenko})
	\[
	C^r(x)= L^r\left(\|x\|\right) \|x\|^{-r\alpha}\,,\quad 0<r\alpha<d \quad\;\implies\; \quad G^{*r}(B_\eps(0))\sim c(r\alpha,d)\,\eps^{r\alpha}L^{r}(1/\eps)\,.
	\]
	This also implies that the ratio \eqref{eq:SACgrowing} is $\asymp 1$, which is perfectly consistent with the mentioned Dobrushin-Major-Taqqu non-central results for $F_{T,q}^{(0)}$.
	
	Less classical (but well studied) covariances are obtained taking stationary fields with
	spectral singularities at finitely many frequencies, usually referred
	to as \textbf{cyclical or seasonal long-memory models}; see e.g.\
	\cite{OlenkoCyclical,OuldHayeSeasonal}. These covariances have the form
\begin{equation}
		\label{eq:longmemoryxi}
		C(x)=
	R(x)
	\left(
	p_0
	+
	\sum_{j=1}^{m}
	p_j\cos(\langle x,\xi_j\rangle)
	\right)\,
\end{equation}
	where $R(x)= L\left(\|x\|\right)\|x\|^{-\alpha}$ is positive definite, $R(0)=1$, $0<\alpha<d/Q$ and $\xi_1,\dots,\xi_m\in\R^d\setminus\{0\}$ are the shifts, with $p_j\ge 0$  associated weights satisfying $p_0+\sum_{j=1}^{m}p_j=1$. 
	By construction, the spectral measure of $C$ is $G=H*\nu$, where $\widehat H=R$ and
	\[
\nu
	=
	p_0\delta_0
	+
	\frac{1}{2}
	\sum_{j=1}^{m}
	p_j
	\bigl(
	\delta_{\xi_j}+\delta_{-\xi_j}
	\bigr)\,,\qquad H^{*r}(B_\eps(0))\sim c(r\alpha,d)\,\eps^{r\alpha}L^{r}(1/\eps)\,,\quad r\le Q\,.
	\]
	Thus, a seasonal covariance function of the form \eqref{eq:longmemoryxi} has spectral measure which is a convolution of a $H$ singular in $0$ and an atomic measure $\nu$ of spectral shifts.
	We also assume that $H$ has a spectral density $h$ which is  bounded and continuous away from $0$, namely
	\begin{equation}
		\label{cond:hdensity}
		H(d\lambda)=h(\lambda)\,\wmu(d\lambda),\qquad
		h\in C(\R^d\setminus\{0\}),\qquad
		\sup_{\|\lambda\|\ge\delta}h(\lambda)<\infty
		\quad \forall\,\delta>0.
	\end{equation}
	 so that the only spectral singularity of $H$ is at $0$. This property is preserved under convolution: for every $r\ge1$,
	 the density $h_r$ of $H^{*r}$ with respect to $\wmu$ satisfies the conditions in \eqref{cond:hdensity}. This condition is satisfied in many situations, for example when $R$ is a Cauchy covariance function, see e.g. \cite[Example 3]{LeonenkoOlenko}. As a consequence, for every $r\le Q$, $G^{*r}$ is absolutely
	 continuous with density $g_r$, which is continuous on
	 $\R^d\setminus\supp(\nu^{*r})$, i.e. the possible spectral singularities
	 are contained in $\supp(\nu^{*r})$. 
	 
	 Under the above assumptions we have a control on the spectral small ball probabilities at every relevant $Q$,$z$ and we can apply \Cref{thm:mainRdZd} to characterize Gaussian fluctuations of $F_{T,q}^{(z)}$.  
	 
	 	\begin{proposition} 
	 	\label{prop:cyclical-long-memory}
	 	Assume $\psi$ satisfies \Cref{cond:schwartz}. Assume \eqref{eq:longmemoryxi} and \eqref{cond:hdensity}. Fix $2\le q\le Q$. Then
	 	\begin{itemize}
	 		\item If $z\notin\supp(\nu^{*q})$ and $g_q(z)>0$, then $\widetilde{F}_{T,q}^{(z)}
	 		\overset{\mathrm{law}}{\rightarrow}
	 		N_z$ and $\Var\left(F_{T,q}^{(z)}\right)
	 		\asymp
	 		T^d$.
	 		\item If $z\in\supp(\nu^{*q})$, then $\widetilde{F}_{T,q}^{(z)}
	 		\overset{\mathrm{law}}{\not\rightarrow}
	 		N_z$ and $\Var\left(F_{T,q}^{(z)}\right)
	 		\asymp
	 		T^{2d}T^{-q\alpha}L^q(T)$. 
	 	\end{itemize}
	 \end{proposition}
	 \begin{proof}
	 	Let us start estimating the small ball probabilities for every $r\le q$. Write \[ \nu^{*r} = \sum_{\eta\in\supp(\nu^{*r})} p_{r,\eta}\delta_\eta, \qquad p_{r,\eta}>0. \] 
	 	Then $g_r(x)=\sum_{\eta\in\supp(\nu^{*r})} p_{r,\eta}h_r(x-\eta)$, so for every \(z\in\supp(\nu^{*r})\), \[ G^{*r}\bigl(B_\varepsilon(z)\bigr) \sim p_{r,z}\,c(r\alpha,d)\, \varepsilon^{r\alpha}L^r(1/\varepsilon), \qquad \varepsilon\rightarrow0, \] whereas, for every $z\notin\supp(\nu^{*r})$, by continuity of $g_r$ at $z$
	 	\[
	 	G^{*r}(B_\varepsilon(z))=O(\varepsilon^d),
	 	\qquad
	 	G^{*r}(B_\varepsilon(z))
	 	\sim \wmu(B_1)\,g_r(z)\varepsilon^d
	 	\quad\text{if }g_r(z)>0,
	 	\]
	 	 Now take $q\le Q$ (with $Q$ as in \eqref{eq:longmemoryxi}) and let us use \Cref{thm:mainRdZd}. 
	 	Assume first that \(z\notin\supp(\nu^{*q})\) and \(g_q(z)>0\). Then $G^{*q}(B_\eps(z))\asymp\eps^d$. Fix \(r=1,\ldots,q-1\), and set \[ \rho:=\dist\bigl(z,\supp(\nu^{*q})\bigr)>0. \] Since $ \supp(\nu^{*r})+\supp(\nu^{*(q-r)}) = \supp(\nu^{*q})$, for every \(\lambda\in\R^d\) either $\dist\bigl(\lambda,\supp(\nu^{*r})\bigr)>\rho/2$ or $ \dist\bigl(z-\lambda,\supp(\nu^{*(q-r)})\bigr)>\rho/2$. Therefore, by \eqref{cond:hdensity}, one of the two factors in the numerator of \eqref{eq:SACgrowing} is uniformly \(O(\eps^d)\). The other one is uniformly \(o(1)\), since \(G^{*s}\) is atomless for every \(s\geq1\). Thus, we conclude applying \Cref{thm:mainRdZd}, after observing that \eqref{eq:SACgrowing} holds \[ \sup_{\lambda\in\R^d} G^{*r}(B_\eps(\lambda)) G^{*(q-r)}(B_\eps(z-\lambda)) = o(\eps^d). \]
	 	Assume now \(z\in\supp(\nu^{*q})\) and write $z=\eta_1+\cdots+\eta_q$, $\eta_j\in\supp(\nu)$. Taking \(\lambda=\eta_1\), we have \[ G(B_\eps(\eta_1)) G^{*(q-1)}(B_\eps(z-\eta_1)) \asymp \eps^{q\alpha}L^q(1/\eps) \asymp G^{*q}(B_\eps(z)). \] Thus \eqref{eq:SACgrowing} fails already for \(r=1\), and we conclude applying \Cref{thm:mainRdZd}.
	 \end{proof}
The proposition yields a frequency-dependent dichotomy in the
long-memory regime, which, to the best of our knowledge, is new
at this level of generality. At frequencies
\(z\in\supp(\nu^{*q})\), spectral singularities persist and
Gaussian convergence fails. Outside this set, provided
\(g_q(z)>0\), spectral delocalization yields a CLT with the
Breuer--Major variance rate.

In dimension one, Rosenblatt \cite{RosenblattFourier} obtained
non-Gaussian Hermite limits at singular frequencies. Our criterion
\eqref{eq:SACgrowing} recovers the failure of Gaussian convergence
and, under \eqref{cond:hdensity}, establishes the complementary
Gaussian regime. Related frequency-overlapping effects were studied
in \cite{IvanovLeonenkoRuizMedinaSavich}, where Gaussian limits
were established for a broad class of weighted nonlinear
functionals of one-dimensional Gaussian processes with multiple
spectral singularities, under the short memory assumption \(q\alpha>d\), see \cite[(A4)]{IvanovLeonenkoRuizMedinaSavich}.
Here we consider the more specific Fourier--Hermite functionals,
for which we characterize Gaussian convergence and its failure
also in the long-memory regime \(Q\alpha<d\).
	
	\subsection{Singular continuous spectra: PDEs, Ahlfors regularity, hypersurfaces and random waves}
	\label{subsec:random-waves-curved-surfaces}
	In this subsection, we only consider the case $E=\R^d$.
	In this setting, a natural source of singular spectral measures is provided by stationary
	Gaussian solutions of \textbf{constant-coefficient linear PDEs}. We follow the notation
	of \cite[Section~2.3]{EstradeFournier}. Let \(P\) be an even
	\(d\)-multivariate polynomial,
	\[
	P(\lambda)
	=
	\sum_{\substack{j\in\mathbb N_0^d\\ |j|\ {\rm even}}}
	\alpha_j\lambda^j,
	\]
	and associate to it the differential operator and \textbf{zero set}
	\[
	\mathcal L_P
	=
	\sum_{\substack{j\in\mathbb N_0^d\\ |j|\ {\rm even}}}
	(-1)^{|j|/2}\alpha_j\partial^j,
	\qquad
	\Sigma_P:=\{\lambda\in\mathbb R^d:P(\lambda)=0\}.
	\]
	It is known, see e.g. \cite[Theorem~2.5]{EstradeFournier}, that if
	\(X\) is regular enough, then
	\[
	\mathcal L_P(X)=0
	\qquad\implies\qquad
	\supp(G)\subseteq\Sigma_P\,.
	\]
	 A canonical example is \textbf{Berry's random wave model} (a.k.a. Euclidean isotropic random wave),
	already discussed in \Cref{subsec:motivations}, a Gaussian field with \(G\) uniformly distributed on
	\(\mathbb S^{d-1}\). This Gaussian field was widely studied in the literature and is the unique (in distribution) isotropic random wave, i.e. the unique isotropic Gaussian field satisfying the Helmholtz equation, with $P$ and $\Sigma_P$ as follows
	\[
	\Delta X +X=0\,,\qquad P(\lambda)=1-\|\lambda\|^2\,\qquad \Sigma_P=\mathbb S^{d-1}\,.
	\]
	More generally, arbitrary symmetric
	probability measures $G$ on \(\mathbb S^{d-1}\) give \textbf{anisotropic random waves}, this point of view is developed e.g. in \cite{EstradeFournier, GassSpectralCriteria}.
    
	 In the infinite depth regime of the \textbf{Gaussian sea-wave model}, see
	\cite{AzaisLeonOrtegaSeaWaves,EstradeFournier,LonguetHiggins1957}, the spectral
	measure support is contained in the Airy surface
	\begin{equation}
	    \label{eq:Airysurface}
        \Sigma_P=\left\{(\xi_1,\xi_2,\tau)\in\mathbb R^3: \xi_1^2+\xi_2^2=\tau^4\right\}\,.
	\end{equation}
	This model is very different from the spherical zero set of random waves, as the Airy surface is
	non-compact and singular at the origin.
	
	These examples motivate the study of spectral measures supported on
	hypersurfaces, or more generally on lower-dimensional subsets of \(\mathbb R^d\),
	without requiring the support to coincide with the whole zero set.
	We will apply our criterion \Cref{thm:mainRdZd} proceeding as follows. We start with general class of weakly Ahlfors regular spectral supports assuming $q=2$, $z=0$, then we strengthen the assumptions on the support and relax those on $q,z$, considering: first, smooth boundaries of convex bodies and general functionals at $z=0$; second, random waves and general Fourier functionals at frequency $z\in \supp(G^{*q})$.
	
\subsubsection{Second chaos, zero frequency and weakly Ahlfors regular spectral supports}

Let \(M\subset\mathbb R^d\) be a Borel set with Hausdorff dimension
$m\in(0,d)$, and set $\mu_M:=\mathcal H^m\!|_M$.
\begin{definition}
	\label{def:quasi-ahlfors}
	\(M\) is \textbf{weakly \(\textbf{m}\)-Ahlfors regular} if $\exists$
	constants \(c,C,\varepsilon_0>0\) such that
	\[
	\mu_M(B_\varepsilon(\lambda))
	\le C\varepsilon^m,
	\qquad
	\lambda\in\mathbb R^d,\quad
	0<\varepsilon<\varepsilon_0,
	\]
	and, for \(\mu_M\)-almost every \(\lambda\in M\), there exists
	\(\varepsilon_\lambda>0\) such that
	\[
	\mu_M(B_\varepsilon(\lambda))
	\ge c\varepsilon^m,
	\qquad
	0<\varepsilon<\varepsilon_\lambda.
	\]
\end{definition}

This is a weakening of the notion of Ahlfors regularity. This class includes also sets like the Airy surface, which is not \(2\)-Ahlfors regular, but is weakly
\(2\)-Ahlfors regular, see \Cref{ex:airy-second-chaos}. Moreover, \(m\) is not necessarily an
integer. Hence the class contains smooth compact manifolds and singular
rectifiable surfaces, but also fractal sets, allowing us to observe peculiar phenomena, see \Cref{ex:cantor}.

The next proposition shows that weak Ahlfors regularity completely
determines the second-chaos behavior at zero frequency, independently of
the particular \(L^2\) density on the spectral support.

\begin{proposition}
	\label{prop:second-chaos-spectral-set}
	Fix \(q=2\) and \(z=0\). Assume that \(\psi\) satisfies
	\Cref{cond:schwartz}, or that \(\psi\) is as in \Cref{lem:variance-density}. Let \(M\subset\mathbb R^d\) be Borel, weakly
	\(m\)-Ahlfors regular for some \(m\in(0,d)\), and suppose that \(X\) has
	symmetric spectral measure
	\[
	dG=f\,d\mu_M,
	\qquad
	f\in L^2(M,\mu_M).
	\]
	Then
	\[
	G^{*2}(B_\varepsilon(0))\asymp\varepsilon^m,
	\qquad
	\widetilde F_{T,2}^{(0)}
	\overset{\mathrm{law}}{\rightarrow}N,
	\qquad
	\Var\left(F_{T,2}^{(0)}\right)\asymp T^{2d-m}.
	\]
\end{proposition}

\begin{proof}
	Write \(\mu=\mu_M\). Since \(M\) is Borel, \(\mu\) is Borel regular. Moreover, by weak Ahlfors regularity $\mu$ is locally finite. Hence \(\mu\) is
	Radon. By the differentiation theorem
	\cite[Corollary~2.14(2)]{Mattila},
	\[
	\frac{G(B_\varepsilon(\lambda))}
	{\mu(B_\varepsilon(\lambda))}
	\rightarrow f(\lambda)
	\qquad\text{for \(\mu\)-a.e. \(\lambda\).}
	\]
	Therefore, by symmetry of \(G\), weak Ahlfors regularity and Fatou's
	lemma,
	\[
	\liminf_{\varepsilon\rightarrow0}
	\frac{G^{*2}(B_\varepsilon(0))}{\varepsilon^m}
	\ge c\int f^2\,d\mu>0.
	\]
	For the upper bound, symmetry, \(2ab\le a^2+b^2\) and weak Ahlfors regularity give
	\[
	G^{*2}(B_\varepsilon(0))
	\le
	\int f(\lambda)^2\mu(B_\varepsilon(\lambda))\,d\mu(\lambda)
	\lesssim
	\varepsilon^m\|f\|_{L^2(\mu)}^2.
	\]
	Finally, Cauchy--Schwarz and the uniform upper growth yield
	\[
	\sup_\lambda G(B_\varepsilon(\lambda))^2
	\lesssim
	\varepsilon^m
	\sup_\lambda
	\int_{B_\varepsilon(\lambda)}f^2\,d\mu
	=o(\varepsilon^m),
	\]
	where the last estimate follows from \(f^2\in L^1(\mu)\) and
	\(\sup_\lambda\mu(B_\varepsilon(\lambda))\to0\). Thus \eqref{eq:SACgrowing} holds for
	\(q=2\) and \(z=0\). Moreover, both \Cref{cond:doublingintro} and the
	rectangular growth required for indicator functions in \Cref{prop:indicator-windows} hold automatically if $q=2$ and $z=0$, see \Cref{rem:alwaysholdscases} and \Cref{superem}.
	Thus, we conclude by \Cref{thm:mainRdZd} and
	\Cref{prop:indicator-windows}.
\end{proof}

As a corollary, we get a CLT for general functionals with non-zero second chaos, under additional minimal integrability assumptions on  $C$. In particular, note that this gives rise to a variety of CLTs in the long memory regime, a direction explored in recent works, see e.g. \cite{GrottoMissingBerryCLT,LeonenkoMainiNourdinPistolato,LeonenkoMainiNourdinPistolatoGneiting,MainiNourdinSpectralCLT}. We postpone the proof to \Cref{sec:additional proofs}.
\begin{corollary}\label{cor:Ahlfors general}
Let $\varphi$ be as in \eqref{eq:phiexp}, with Hermite coefficients $(c_q)_{q\ge1}$, and assume that $c_2\neq0$. If $G$ and $\psi$ are as in \Cref{prop:second-chaos-spectral-set}, and $C\in L^p(\R^d)$ for some $p>0$, then 
\[
\widetilde F_{T}^{(0)}
\overset{\mathrm{law}}{\rightarrow}N,
\qquad 
\Var\left(F_{T}^{(0)}\right)\asymp c_1^2\,\Var\left(F_{T,1}^{(0)}\right) +  T^{2d-m}.
\]
\end{corollary}

Including also indicators $\psi=\mathbf{1}_D$, $D\subseteq \R^d$, in the previous results gives access also to geometric functionals such as  excursion volumes at non $0$ level, namely
\begin{equation*}
	F_T^{(0)}=\Vol\left(\left\{x\in TD\,:\,X(x)\ge u\right\}\right)\,,\qquad\varphi=\mathbf{1}_{[u,\infty)}\,,\quad u\neq0\,,
\end{equation*}
Indeed, in this case $c_2\neq0$, see e.g. \cite[Proof of Lemma 9]{MainiMalliavin}. When $u=0$ we have only odd chaoses and we need different arguments, see the next sections.

Note that the proposition works for smooth compact $M$, but requires in general neither compactness nor smoothness.
We illustrate this with rather different examples: a fractal spectral
measure, the sea-wave model associated with the
Airy surface and anisotropic random waves.

\begin{example}[Cantor spectral measure]
	\label{ex:cantor}
	Let \(M\subseteq \R \) be the centered middle-third Cantor set, $M
	=
	\left\{
	\sum_{k\ge1}\varepsilon_k3^{-k}:
	\varepsilon_k\in\{-1,1\}
	\right\}$. 
	Since \(M\) is \(m\)-Ahlfors regular with $
	m=\frac{\log2}{\log3},$
\footnote{At the $n$th step of its construction, the Cantor set is covered by $2^n$ disjoint intervals $I_{n,j}$ of length $\ell_n=3^{-n}$, separated by gaps of length at least $\ell_n$. $M$ is self-similar and its Hausdorff dimension is $\log(2)/\log(3)$, see \cite{Falconer}. By self-similarity,
$\mathcal H^m(M\cap I_{n,j})=\ell_n^m\mathcal H^m(M)=2^{-n}\mathcal H^m(M)$. Given $3^{-(n+1)}<r\le 3^{-n}$, any ball $B_r(y)$ meets a uniformly finite number of $n$-step intervals and contains a $n+1$-step interval, so
$\mathcal H^m(M\cap B_r(y))\asymp 2^{-n}\asymp 3^{-nm}\asymp r^m$,
uniformly in $y\in M$.}

	\[
	\widetilde F_{T,2}^{(0)}
	\overset{\mathrm{law}}{\rightarrow}N,
	\qquad
	\Var\left(F_{T,2}^{(0)}\right)
	\asymp T^{\,2-\log2/\log3}.
	\]
	The stationary Gaussian process $X$ here has covariance
	\[
	C(t)=\prod_{k\ge1}\cos(3^{-k}t);
	\]
	In particular,
	\(C(3^n\pi)=(-1)^nC(\pi)\), so \(C(t)\not\rightarrow 0\). Related processes have been considered e.g. in \cite{LachiezeRey}.
\end{example}

\begin{example}[Gaussian sea waves]
	\label{ex:airy-second-chaos}
	Consider the space-time Gaussian sea-wave model discussed above, with spectral support contained in the Airy surface $M$, see \eqref{eq:Airysurface}.
	Here \(d=3\) and \(m=2\). The surface is smooth away from the origin,
	while \(\{0\}\) is \(\mathcal H^2\)-negligible. Parametrizing the surface, one can check that 
it is weakly 2-Ahlfors regular
\footnote{
The uniform upper bound follows since $M$ is covered by finitely
many graphs with uniformly bounded gradients.
Indeed, for $|t|\ge1$ we write
$t=\pm(x^2+y^2)^{1/4}$, whose gradient has norm
$1/(2|t|)\le1/2$. For $|t|\le1$, on the region
$|x|\ge|y|$ we write
$x=\pm\sqrt{t^4-y^2}$, with
$|\partial_t x|\le2\sqrt2$ and
$|\partial_y x|\le1$; the region $|y|\ge|x|$
is treated symmetrically.
For the lower bound, since $M\setminus\{0\}$ is smooth,
$\mathcal H^2(M\cap B_r(p))\sim\pi r^2$
for every $p\in M\setminus\{0\}$ as $r\downarrow0$.
}
, even if it is not \(2\)-Ahlfors regular in the usual sense, because
	\(\mathcal H^2(M\cap B_\varepsilon(0))
	\asymp\varepsilon^3\).
	Thus, by \Cref{prop:second-chaos-spectral-set}, every symmetric probability measure
	\(dG=f\,d\mu_M\) with
	\(f\in L^2(\mu_M)\) satisfies
	\[
	\widetilde F_{T,2}^{(0)}
	\overset{\mathrm{law}}{\rightarrow}N,
	\qquad
	\Var\left(F_{T,2}^{(0)}\right)\asymp T^4.
	\]
    CLTs for other geometric
	characteristics of Gaussian sea-wave models were obtained previously in
	\cite{AzaisLeonOrtegaSeaWaves}. 
\end{example}

\begin{example}[Random waves with $L^2$ density] Let $M=\mathbb S^{d-1}$ and $X$ be a random wave with $L^2$ spherical spectral density, that is, assume $X$ has spectral measure $dG=f\,d\mu_M$ with $f\in L^2(\mu_M)$. In this case, one has also $C\in L^p(\R^d)$ for $p$ large enough (see e.g. \cite[Chapter VIII, Theorem 3]{SteinHarmonicAnalysis}) and so by \Cref{cor:Ahlfors general}, if $c_2\neq0$  (since $0\not\in \supp(G)\implies \Var(F_{T,1}^{(0)})=o(T^d)$ by reasoning as in \eqref{eq:negligiblefirst}) we have 
\[
\widetilde F_{T}^{(0)}
\overset{\mathrm{law}}{\rightarrow}N,
\qquad
\Var\left(F_{T}^{(0)}\right)\asymp T^{2d-(d-1)}=T^{d+1}.
\]
This is consistent with the results in \cite{BeliaevHegde,GassPeccatiStecconi,MainiNourdinSpectralCLT} for Berry's random wave model, and generalizes them to every situation in which one has a generic $L^2$ spherical spectral density. The result holds for example for the excursion volume at level $u\neq0$.
\end{example}

\subsubsection{General chaos at frequency $0$ and smooth spectral curved hypersurfaces}

For $q=2$ we already proved a CLT under general assumptions on the local behavior of $\mu_M$ on $M$, i.e. weakly Ahlfors regularity. For $q\ge3$, we need stronger assumptions on $M$. Here local regularity is not enough and spectral delocalization depends more deeply on the geometry of $M$.

We start with a technical geometric lemma that gives us a more precise control on $G^{*2}$. The proof of the lemma, together with some geometric preliminaries about co-area formula and Ehresmann theorem, are postponed to \Cref{sec:ehresmann}.

\begin{lemma}
	\label{lem:twofold-convolution}
	Let \(K\subset\mathbb R^d\), \(K=-K\), \(d\ge2\), be compact and
	strictly convex, with \({\rm int}(K)\neq\emptyset\) and smooth boundary
	\(M:=\partial K\) with non-vanishing Gaussian curvature. Then \(G:=\mu_M/\mu_M(M)\) is
	\((d-1)\)-Ahlfors regular: there exist \(r_0,c,C>0\) such that
	\[
	c r^{d-1}\le G(B_r(x))\le C r^{d-1},
	\qquad x\in M,\quad 0<r<r_0.
	\]
	Moreover \(G^{*2}(dx)=g_2(x)\,\wmu(dx)\), with $g_2$ continuous and positive on ${\rm int}(2K)\setminus\{0\}$. Finally, if $d=2$\footnote{This asymptotic in $0$ is true for every $d\ge2$, however we will need it only for $d=2$ and the proof is simpler in this case.}
	\[
	g_2(z)\asymp\|z\|^{-1},
	\qquad
	0<\|z\|<\delta\,\qquad\text{ for some }\delta>0\,.
	\]
\end{lemma}

Combining \Cref{thm:mainRdZd}, \Cref{prop:indicator-windows}, \Cref{prop:second-chaos-spectral-set} and \Cref{lem:twofold-convolution}, we get a CLT for every Fourier-Hermite functional $F_{T,q}$ at frequency $0$ when $dG=f\,d\mu_M$ and $0<c\le f\le C<\infty$. Then, \Cref{lem:finite-chaos-reduction} yields the same result for every functional $F_T$ with $\varphi-\E[\varphi(N)]$ non-linear. 

\begin{proposition}\label{third chaos general}
	Let \(K\subset\mathbb R^d\), \(K=-K\), \(d\ge2\), be compact,
	strictly convex with \({\rm int}(K)\neq\emptyset\) and smooth boundary
	\(M:=\partial K\) with non-vanishing Gaussian curvature.
	Set \(dG:=f\,d\mu_M\), where $0<c\le f\le C<\infty$ on $M$.
	Assume that $\psi$ satisfies \Cref{cond:schwartz}, or is the indicator of a convex body as in \ref{psiliscio}. Then \(\Var\left(F_{T,1}^{(0)}\right)=o(T^d)\) and for \(q\ge2\)
	\[
	\widetilde F_{T,q}^{(0)}
	\overset{\mathrm{law}}{\rightarrow}N,
	\qquad
	\Var\left(F_{T,q}^{(0)}\right)\asymp
	\begin{cases}
		T^{d+1}, & q=2,\\
		T^2\log T, & d=2,\ q=4,\\
		T^d, & \text{otherwise}.
	\end{cases}
	\]
	Moreover, if \(\varphi\in L^2(\mathbb R,\gamma)\) is as in
	\eqref{eq:phiexp} and \(\varphi-\E[\varphi(N)]\) is not linear, then
	\[
	\widetilde F_T^{(0)}
	\overset{\mathrm{law}}{\rightarrow}N,
	\qquad
	\Var\left(F_T^{(0)}\right)\asymp
	\begin{cases}
		T^{d+1}, & c_2\ne0,\\
		T^2\log T, & d=2,\ c_2=0,\ c_4\ne0,\\
		T^d, & \text{otherwise}.
	\end{cases}
	\]
\end{proposition}

\begin{proof}
	Let us first deal with $F_{T,q}^{(z)}$. By
	\Cref{cor:equivalentspectralmeasures}, WLOG
	\(f\equiv1/\mu_M(M)\). Since \(0\notin M=\supp(G)\) and $|\wpsi(\lambda)|^2\lesssim \left(1+\|\lambda\|\right)^{-d-1}$ \footnote{This follows by the assumptions on $D$ in \ref{psiliscio}, see \cite{BHI03}.}, by \eqref{eq:variance-kernel} we have
    \begin{equation}
        \label{eq:negligiblefirst}
\Var\left(F_{T,1}^{(0)}\right)=\int_{\R^d}\left|\wpsi_T(\lambda)\right|^2G(d\lambda)\lesssim T^{2d}\int_{\R^d}\frac{G(d\lambda)}{\left(T\|\lambda\|\right)^{d+1}}=o(T^d)\,.
    \end{equation}
	By \Cref{lem:twofold-convolution}, \(M\) is \((d-1)\)-Ahlfors regular,
	so the result for \(q=2\) follows by \Cref{prop:second-chaos-spectral-set}. 
    
    Setting
	\(M_\delta:=\{x\in\R^d:{\rm dist}(x,M)<\delta\}\), for $\eps$ small enough we have by Ahlfors regularity 
	\[
	0<
	\inf_{x\in M_{\eps/2}}
	\frac{G(B_\eps(x))}{\eps^{d-1}}
	\le
	\sup_{x\in\R^d}
	\frac{G(B_\eps(x))}{\eps^{d-1}}
	<\infty.
	\]
	Since \(G(B_\eps(x))=0\) if \(x\notin M_\eps\) and, by
	\Cref{lem:twofold-convolution}, \(g_2\) is continuous and positive close to \(M\), we get
	\[
	\eps^{d-1}\mathcal H^d(M_{\eps/2})
	\lesssim
	G^{*3}(B_\eps(0))
	=
	\int_{\R^d}g_2(x)G(B_\eps(x))\,\wmu(dx)
	\lesssim
	\eps^{d-1}\mathcal H^d(M_\eps).
	\]
	Since \(\mathcal H^d(M_\eps)\asymp\eps\), we get
	\(G^{*3}(B_\eps(0))\asymp\eps^d\). Moreover, again by continuity of $g_2$ on $M_\eps$,
	\[
	\sup_{x\in\R^d}
	G(B_\eps(x))G^{*2}(B_\eps(x))
	\lesssim
	\eps^{d-1}
	\sup_{x\in M_\eps}G^{*2}(B_\eps(x))
	\lesssim
	\eps^{2d-1}.
	\]
	Thus \eqref{eq:SACgrowing} holds and we conclude the case \(q=3\) by
	\Cref{thm:mainRdZd} and \Cref{prop:indicator-windows}. Indeed, the ball growth condition in \Cref{prop:indicator-windows} (and in particular the doubling \Cref{cond:doublingintro} assumed in \Cref{thm:mainRdZd}) hold, since $G^{*3}(B_\eps(0))\asymp \eps^d$.
	
	It remains to treat \(q\ge4\). Since \(M\) is smooth with non-vanishing
	Gaussian curvature, by \cite[Chapter VIII, Theorem 1]{SteinHarmonicAnalysis}
	we have
	\[
	\left|\widehat G(x)\right|=|C(x)|
	\lesssim
	\|x\|^{-(d-1)/2},
	\]
	so that \(C\in L^q(\R^d)\) if \(d\ge3,\ q\ge4\), or if \(q\ge5\).
	Thus the conclusion follows from the Breuer--Major theorem \Cref{thm} if
	\(g_q(0)>0\). For even \(q\), $g_q(0)=\int C^q >0$, since $C$ continuous and $C(0)=1$. For \(q=3\), for some
	\(c,\eta>0\),
	\[
	g_3(z)
	=
	\int_M g_2(z-x)\,G(dx)
	\ge
	\inf_{y\in M_\eta}g_2(y)
	\ge c,
	\qquad
	\|z\|<\eta.
	\]
	Hence, for $q\ge5$ odd, since $G$ symmetric $\implies$
$0\in\supp(G^{*(q-3)})$, we have
\[
g_q(0)
=
\int_{\R^d}g_3(x)\,G^{*(q-3)}(dx)
\ge
c\,G^{*(q-3)}(B_\eta(0))
>0.
\]
	It remains the case \(d=2,\ q=4\). Here the covariance decay gives\footnote{Since $\Var\left(F_{T,4}^{(0)}\right)=4!\int_{(\R^2)^2}C^4(x-y)\psi_T(x)\psi_T(y)\,dx\,dy=4!T^2\int_{\R^2}C^4(z)\psi*\psi_-(z/T)\,dz$, $\psi_-:=\psi(-\cdot)$, see e.g. \cite{MainiSPA}, by covariance decay and $\psi\in L^1\cap L^\infty$ we get $\Var\left(F_{T,4}^{(0)}\right)\lesssim T^2\int_{\|z\|\le T}\left(1+\|z\|^2\right)^{-1}\,dz+T^2\|\psi*\psi_-\|_1\lesssim T^2\log(T)+T^2$.}
	\(\Var\left(F_{T,4}^{(0)}\right)\lesssim T^2\log (T)\). Since $\sup_{\lambda\in\R^d} G^{*r}(B_\eps(\lambda))\lesssim \sup_{\lambda\in\R^d} G(B_\eps(\lambda))\lesssim \eps$, the numerator of \eqref{eq:SACgrowing} is $O(\eps^2)$. By \Cref{lem:twofold-convolution}
	\[
	G^{*4}(B_\eps(0))
	\ge
	\int_{\|z\|<\eps}
	\int_{2\eps<\|x\|<\delta}
	g_2(x)g_2(z-x)\,\wmu(dx)\,\wmu(dz)
	\gtrsim
	\eps^2
	\int_{2\eps<\|x\|<\delta}
	\frac{dx}{\|x\|^2}
	\asymp
	\eps^2\log(1/\eps),
	\]
	where we used \(\|z-x\|\asymp\|x\|\) on the integration domain. Thus \eqref{eq:SACgrowing} vanishes and we conclude by \Cref{thm:mainRdZd} and \Cref{prop:indicator-windows}, noting that the ball growth condition in \Cref{prop:indicator-windows} holds because $q$ is even and $z=0$, see  \Cref{superem}.
    
    Finally, for general nonlinear functionals, we have to consider again the general density $f$. Since $f\in L^2(\mu_M)$, by
	\cite[Chapter VIII, Theorem 3]{SteinHarmonicAnalysis} we have
	\(C\in L^m(\R^d)\) for some integer \(m\) large enough. Thus, for $q\ge m$ by \Cref{lem:variance-density}, \eqref{ortovar} and \eqref{eq:phiexp} we have $\Var(\sum_{q\ge m}c_qF_{T,q}^{(z)})\lesssim T^d$.
	Thus, the fixed-chaos estimates above, the fact that $\exists$ $q\ge2$ with $c_q\neq0$ and again \eqref{ortovar} give
	\[
	\Var\left(F_T^{(0)}\right)\asymp
	\begin{cases}
		T^{d+1}, & c_2\ne0,\\
		T^2\log T, & d=2,\ c_2=0,\ c_4\ne0,\\
		T^d, & \text{otherwise},
	\end{cases}\,.
	\]
	Since \(\Var\left(F_T^{(0)}\right)\gtrsim T^d\), $\Var(F_{T,1})=o(T^d)$ and $\widetilde F_{T,q}\overset{\rm law}{\rightarrow}N$ for $q\ge2$, and we conclude by \Cref{lem:finite-chaos-reduction}.
\end{proof}
The previous result recovers and generalizes many results on Berry's random wave, which
corresponds to \(M=\mathbb S^{d-1}\) and \(f\equiv 1/\mu_M(M)\) constant on the sphere; see
e.g.\
\cite{BeliaevHegde,GassPeccatiStecconi,GrottoMissingBerryCLT,
	GrottoMainiTodino,MainiNourdinSpectralCLT}.
In particular, fluctuations and variances rates do not change if we take $M$ as above.
\begin{remark}
    In the works \cite{BeliaevHegde,MainiSPA,GassPeccatiStecconi}, the authors also proved stronger multivariate or functional CLTs for functionals controlling the asymptotic covariances of the field $(F_T(\psi))_{\psi}$. In particular, \cite{GassPeccatiStecconi} highlighted a ${\rm FGF}_{1/2}$ asymptotic structure and related it to the presence of "random scars". We believe that these results and interpretations may be extended (with different interesting asymptotic structure and different kind of "scars") in this context. We leave this for future research.
\end{remark}

\subsubsection{Fourier functionals of random waves}

Let \(\sigma\) be the normalized surface measure on
\(\mathbb S^{d-1}\). For \(G=\sigma\), the convolution \(\sigma^{*q}\) is
the law of a \(q\)-step Pearson random walk; we use the following
density estimates, see e.g. \cite{BorweinStraubVignat,BorweinStraubWanZudilin}. 
We shall use the following facts: by \cite[Lemma~4.1]{BorweinStraubVignat} we have
	\begin{equation}
	    g_2(x)=a_d\|x\|^{-1}\left(4-\|x\|^2\right)^{(d-3)/2}
	\mathbf 1_{\{0<\|x\|<2\}},\qquad a_d>0.\label{eq:g2berry}
	\end{equation}
If \(d=2\), then \(g_3\) is continuous, positive on
	\(B_3(0)\setminus\mathbb S^1\), and
	\(g_3(x)\asymp\log\left(1/\left|\|x\|-1\right|\right)\) as \(\|x\|\to1\). This follows
	from \cite[(3.1)--(3.3)]{BorweinStraubWanZudilin}. If \(d=2\), then \(g_4\) is continuous and positive on
	\(B_4(0)\setminus\{0\}\); see
	\cite[ Example~2.6]{BorweinStraubVignat},
	together with $g_4(z)=g_2*g_2(z)$.
In particular, if $d=2$ these estimates imply
\begin{equation*}
    \label{eq:berryestimates}
    \sup_x\sigma^{*2}(B_\eps(x))\lesssim\eps\,,\qquad\sup_{\|x\|\ge\delta}\sigma^{*2}(B_\eps(x))
\lesssim\eps^{3/2}\,,\qquad \sup_x\sigma^{*3}(B_\eps(x))
\lesssim\eps^2\log(1/\eps)\,.
\end{equation*}

 Together with
\Cref{thm:mainRdZd} and \Cref{cor:equivalentspectralmeasures}, these bounds give
the following classification. 

\begin{proposition}[Random waves at every fixed frequency]
	\label{prop:anisotropic-Berry-all-frequencies}
	Let \(dG=f\,d\sigma\), where \(0<c\le f\le C<\infty\).
	Assume $\psi$ satisfies \Cref{cond:schwartz}. Fix \(d,q\ge2\). 	If \(\|z\|>q\), then
	\(\Var\left(F_{T,q}^{(z)}\right)=O(T^{-N})\) $\forall \,\,N>0$. 
	If \(\|z\|\le q\), we have
	\(\widetilde F_{T,q}^{(z)}\overset{\rm law}{\not\rightarrow}N_z\)
	$\;\Leftrightarrow\;$ \(d=q=2\) and \(0<\|z\|<2\). Moreover, the variance rates are
	\[
	\Var\left(F_{T,q}^{(z)}\right)\asymp
	\begin{cases}
		T^{d+1}, & q=2,\ z=0,\\
		T^2\log T, & d=2,\ q=3,\ \|z\|=1,\\
		T^2\log T, & d=2,\ q=4,\ z=0,\\
		T^{2d-(q+1)(d-1)/2}, & \|z\|=q,\\
		T^d, & \text{otherwise}.
	\end{cases}
	\]
\end{proposition}

\begin{proof}
	By \Cref{cor:equivalentspectralmeasures}, WLOG we can take
	\(G=\sigma\). Moreover $|C(x)|\lesssim \|x\|^{-\frac{d-1}{2}}$, see \cite[Chapter VIII]{SteinHarmonicAnalysis}.
	Throughout the proof, \Cref{cond:doublingintro} will follow from our estimates, so \Cref{thm:mainRdZd} applies.
	
We first treat the case \(\|z\|=q\). Geometric arguments
	\footnote{For \(j=1\), we conclude by $(d-1)$-Ahlfors regularity of $\mathbb S^{d-1}$.
			Assume the claim at order \(j\). By convolution,
			\[
			\sigma^{*(j+1)}(B_\eps((j+1)e))
			=
			\int_{\mathbb S^{d-1}}
			\sigma^{*j}(B_\eps((j+1)e-x))\,\sigma(dx).
			\]
			Set \(v=(j+1)e-x\). Since
			\(\|v\|^2=j^2+(j+1)\|x-e\|^2\) and
			\(\supp(\sigma^{*j})\subseteq\overline B_j(0)\),
			the integrand vanishes unless \(\|v\|-j<\eps\),
			hence \(\|x-e\|\le\sqrt{2\eps}\).
			Moreover, \(B_\eps(v)\subset B_{2\eps}\left(jv/\|v\|\right)\).
			The induction hypothesis and the cap estimate give the upper bound.
			Analogously, we lower bound restricting to \(\|x-e\|\le c\sqrt\eps\).
			For \(c\) small, \(\|v\|-j\le\eps/2\), so
			\(B_{\eps/2}\left(jv/\|v\|\right)\subset B_\eps(v)\).} 
		give
	\(\sigma^{*j}(B_\eps(je))\asymp
	\eps^{(j+1)(d-1)/2}\) $\forall j\ge1$, uniformly for $e\in\mathbb S^{d-1}$. Thus, we have a CLT by \Cref{thm:mainRdZd}, since \eqref{eq:SACgrowing} is bounded by
	\[
	\eps^{-(q+1)(d-1)/2}\sup_{{\rm dist}(\lambda,r\mathbb S^{d-1})\le \eps}\sigma^{*r}(B_\eps(\lambda))\sup_{{\rm dist}(\lambda',(q-r)\mathbb S^{d-1})\le \eps}\sigma^{*(q-r)}(B_\eps(\lambda'))
	\,\,\lesssim\,\,
	\eps^{(d-1)/2}.
	\]
	If $\|z\|>q$, then $z\not\in \supp \sigma^{*q}$ and we apply \Cref{thm:mainRdZd}.
	Assume now \(0<\|z\|<q\), since $z=0$ follows from \Cref{third chaos general}. We have 4 cases.
	
	\emph{Case \(d\ge3\).}
	For \(q=2\),
	when \(z\ne0\) \eqref{eq:g2berry} gives
	\(\sigma^{*2}(B_\eps(z))\asymp\eps^d\). Since
	\(\sup_x\sigma(B_\eps(x))\lesssim\eps^{d-1}\), the ratio in
	\eqref{eq:SACgrowing} is \(O(\eps^{d-2})=o(1)\). Since $C\in L^q$ for $q\ge4$ when $d\ge3$, the only remaining case not
	covered by \Cref{thm} is \((d,q)=(3,3)\). Here
	\(g_3=g_2*\sigma\) is continuous and positive on \(B_3(0)\), since the
	singularity \(\|z-\cdot\|^{-1}\) is locally integrable on
	\(\mathbb S^2\), so $G^{*3}(B_\eps(z))\asymp\eps^d=\eps^3$, whereas
	the numerator is \(O(\eps^4)\). 
	
	\emph{Case \(d=2\), \(q=2\) or \(q=4\).}
	If \(q=2\) and \(0<\|z\|<2\), then
	\(\sigma^{*2}(B_\eps(z))\asymp\eps^2\) and there exist \(x,y\in\mathbb S^1\) with \(z=x+y\), so the numerator in \eqref{eq:SACgrowing} is also of order \(\eps^2\) and by \Cref{thm:mainRdZd} we do not have a CLT. For \(q=4\) and \(0<\|z\|<4\), continuity and
	positivity of \(g_4\) give a denominator of order \(\eps^2\). Thus, if $r=1$ in \eqref{eq:SACgrowing} the numerator has order
	\(O(\eps^3\log(1/\eps))\), while for $r=2$ at least one center stays away
	from the origin, giving \(O(\eps^{5/2})\). Thus spectral delocalization
	holds and by \Cref{thm:mainRdZd} we have a CLT.
	
	\emph{Case \(d=2\), \(q=3\).}
	The mentioned estimates on \(g_3\) give
	\[
	\sigma^{*3}(B_\eps(z))\asymp
	\begin{cases}
		\eps^2\log(1/\eps), & \|z\|=1,\\
		\eps^2, & \|z\|<3,\ \|z\|\ne1.
	\end{cases}
	\]
	If \(\|z\|=1\), the numerator in \eqref{eq:SACgrowing} is
	\(O(\eps^2)\). If \(\|z\|\ne1\), $G(B_\eps(x))>0$ only if $|\|x\|-1|<\eps$ $\implies$ $\|z-x\|\ge |\|z\|-1|-\eps$, so for $\eps$ sufficiently small $\sup_{|\|x\|-1|<\eps}G^{*2}(B_\eps(z-x))=O(\eps^{3/2})$ and the numerator in \eqref{eq:SACgrowing} is \(O(\eps^{5/2})\). We conclude by \Cref{thm:mainRdZd}.
	
	\emph{Case \(d=2\), \(q\ge5\).} $C\in L^q$ and so we conclude by positivity of $g_q(z)$ and \Cref{thm}. The positivity of $g_q(z)$ for $0<\|z\|<q$ can be deduced by that of $g_2$, $g_3$ respectively on $B_2(0), B_3(0)$ and induction.
\end{proof}

\begin{figure}[htbp]
    \centering
    \includegraphics[scale=0.8]{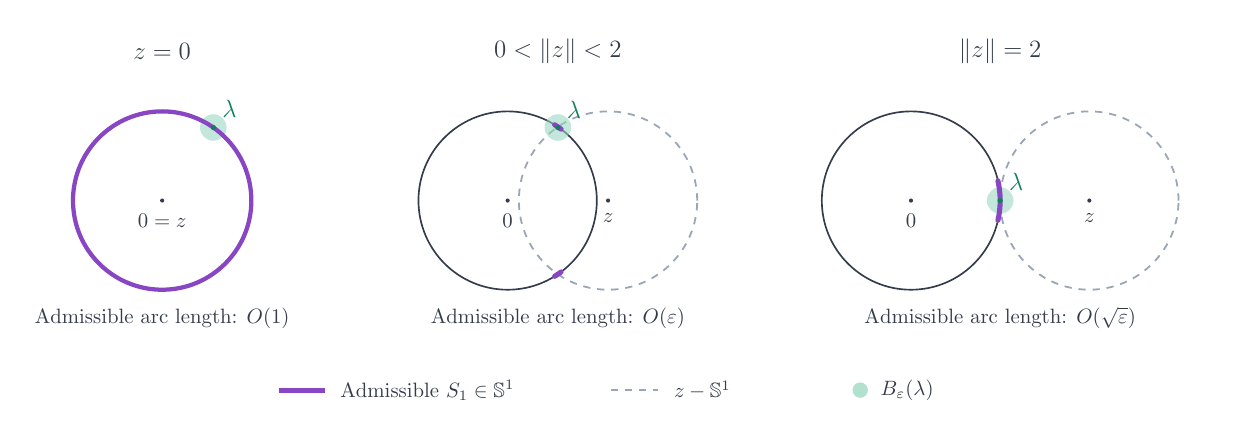}
    \caption{The purple arcs represent the admissible positions of
$S_1$ given $S_2\in B_\varepsilon(z)$, while the green disks the admissible positions for $S_1\in B_\varepsilon(\lambda)$ given $S_2\in B_\varepsilon(z)$. If $z=0$ every position is admissible for $S_1$, so the event $S_1\in B_\eps(\lambda)$ is negligible as $\eps \rightarrow 0$. If $0<\|z\|<2$, $S_1$ concentrates on two arcs of length $\eps$, so the event $S_1\in B_\eps(\lambda)$ is not negligible as $\eps \rightarrow 0$. If $\|z\|=2$, $S_1$ concentrates on an arc of length $\sqrt{\eps}$, so $S_1\in B_\eps(\lambda)$ is negligible as $\eps \rightarrow 0$.}
\label{fig:spectral-delocalization}
\end{figure}

While \Cref{third chaos general} extends existing zero-frequency results for Berry's random wave to more general curved spectral supports, \Cref{prop:anisotropic-Berry-all-frequencies} extends them by varying the probed frequency. For second chaos Fourier functionals of planar random waves the frequency $z=0$ is special, as for planar arithmetic random waves on the torus (see \Cref{subsec:arithmetic-random-waves}). An interpretation of this phenomenon with the conditional formulation \eqref{eq:condform} of spectral delocalization is given in \Cref{fig:spectral-delocalization}.

	\section{Preliminaries}\label{sec:preliminaries}
	
\subsection{Locally compact Abelian groups and spectral representations}
	\label{sec:lca-spectral}
	
	We collect here preliminaries related to topology and harmonic analysis on LCA groups, for missing details
	see e.g. 
	\cite[Chapter~4]{FollandAHA}.

    A topological group is a group $(E,+)$ endowed
with a topology for which addition and inversion are
continuous. In particular, translations are
homeomorphisms. For a Hausdorff locally compact
Abelian (LCA) group $(E,+)$, its Pontryagin dual
$\widehat E$ is the group of continuous homomorphisms
from $E$ into the unit circle $\mathbb S^1\subseteq\C$,
called characters, endowed with pointwise multiplication
and the topology of uniform convergence on compact
subsets of $E$. The dual $\widehat E$ is again a
Hausdorff LCA group, see \cite[Section 4.1]{FollandAHA}, and we have the Pontryagin duality  $E\simeq \widehat \wE$, see  \cite[Theorem 4.32]{FollandAHA}, where $\simeq$ denotes the existence of an isomorphism of topological groups: a continuous group isomorphism $E\to \widehat\wE$ with continuous inverse. 

For $x\in E$ and a character $\lambda\in\widehat E$, we adopt
the exponential notation
\begin{equation*}
e^{i\langle x,\lambda\rangle}:=\lambda(x).
\end{equation*}
In general, $\langle x,\lambda\rangle$ is the argument/phase in $[0,2\pi)$ of the complex number $\lambda(x)$. Consistently with the exponential notation, we use additive notation on $\widehat E$, although
its group operation is pointwise multiplication of
characters. In other words, for $\lambda,\eta\in\widehat E$, we write
\[
\begin{aligned}
e^{i\langle x,\lambda+\eta\rangle}
&=(\lambda+\eta)(x)
=\lambda(x)\eta(x)
=e^{i\langle x,\lambda\rangle}
 e^{i\langle x,\eta\rangle},\\
e^{i\langle x,-\lambda\rangle}
&=(-\lambda)(x)
=\overline{\lambda(x)}=\lambda(x)^{-1}
=e^{-i\langle x,\lambda\rangle}.
\end{aligned}
\]
In particular, the identity $0\in\widehat E$
is the trivial character, satisfying
$e^{i\langle x,0\rangle}=1$ for every $x\in E$.

\subsubsection{Some topological properties of LCA groups}
We have the following properties.

\begin{lemma}\label{lem:lca-topology}
Let $(E,+)$ be a Hausdorff LCA group.
\begin{enumerate}[label=(\alph*),leftmargin=2em]
\item
$E$ is first countable if and only if it admits
a compatible translation-invariant metric.
\item
$E$ is second countable if and only if it is
metrizable and $\sigma$-compact.
\item
$E$ is metrizable if and only if  $\widehat E$ is $\sigma$-compact. $E$ is $\sigma$-compact if and only if  $\widehat E$ is metrizable.
\item
$E$ is second countable if and only if
$\widehat E$ is second countable.
\end{enumerate}
\end{lemma}

Part (a) is the Birkhoff--Kakutani theorem
\cite{Birkhoff1936,Kakutani1936}, since $E$ is
Hausdorff and Abelian. Part (b) is the
characterization of second countability for
locally compact groups, see
\cite[Section~2.1]{Reid2021}, while (c) is
\cite[Lemma~2.8(i)]{Reid2021} and Pontryagin duality. (d) follows immediately by (b) and (c).

Throughout the paper, we will always assume that our Gaussian field $X$ is indexed by a Hausdorff, second
countable, LCA group $E$. By \Cref{lem:lca-topology}, $\wE$ is second countable by (d) and metrizable with translation invariant metric by (a). Thus, we fix a translation-invariant compatible metric \(d_{\widehat E}\), satisfying
	\[
	d_{\widehat E}(\lambda+\eta,\lambda'+\eta)
	=
	d_{\widehat E}(\lambda,\lambda'),
	\qquad \lambda,\lambda',\eta\in\widehat E.
	\]
and we denote by \(B_\varepsilon(\lambda)\) the corresponding open balls, for which we have $B_\varepsilon(\lambda)=\lambda+B_\varepsilon(0)$, $\lambda\in\wE$.

\subsubsection{Harmonic analysis and Bochner's theorem on LCA groups}   

Fix now a Haar measure $\mu$ on $E$. By Plancherel's
theorem for LCA groups, see e.g. \cite[Theorem~4.26]{FollandAHA},
there is a unique Haar measure $\widehat\mu$
on $\widehat E$ such that the Fourier transform
\begin{equation}
    \widehat f(\lambda)
:=\int_E f(x)e^{i\langle x,\lambda\rangle}\,\mu(dx),
\qquad f\in L^1(E,\mu):=L^1(E,\C,\mu),\label{eq:Fourierdef}
\end{equation}
extends to a unitary operator from $L^2(E,\mu)$ to $L^2(\widehat E,\widehat\mu)$, so that Plancherel identity holds
\begin{equation}
    \langle f,g\rangle_{L^2(E,\mu)}=\left\langle\, \widehat f\,,\,\widehat g\, \right\rangle_{L^2\left(\wE,\wmu\right)}\,.\label{eq:plancherel}
\end{equation}
If \(G\) is a finite Borel measure on \(\widehat E\), we also set
	\[
	\widehat G(x)
	=
	\int_{\widehat E}e^{i\langle x,\lambda\rangle}\,G(d\lambda),
	\qquad x\in E.
	\]
	The usual relation between products and convolutions remains valid on LCA
	groups. If \(G,G'\) are Borel probability measures on \(\widehat E\), their
	convolution \(G*G'\) is defined in the standard way and is the law
	of \(\Lambda+\Lambda'\) when \(\Lambda\) and \(\Lambda'\) are independent
	with laws \(G\) and \(G'\). In particular,
	\(\widehat{G*G'}=\widehat G\,\widehat G'\).
	The analogous identity
	\(\widehat{f*g}=\widehat f\,\widehat g\) holds for integrable functions.
    
	We shall repeatedly use \textbf{Bochner's theorem}, see e.g. \cite[Theorem 4.19]{FollandAHA}: a continuous function
	\(C:E\to\mathbb R\), with \(C(0)=1\), is positive definite if and only if
	there is a unique symmetric Borel probability measure \(G\) on \(\widehat E\),
	called the spectral measure of \(C\), such that
	\[
	C(x)
	=\widehat G(x)=
	\int_{\widehat E}e^{i\langle x,\lambda\rangle}\,G(d\lambda),
	\qquad x\in E.
	\]
	This also explains why convolution powers of \(G\) appear naturally
	throughout the paper. Indeed, if \(C=\widehat G\), then
	\(C^q=\widehat{G^{*q}}\). Hence, for a centered stationary Gaussian field
	\(X\) with covariance \(C\), the Hermite isometry identity (see e.g. \cite{NourdinPeccatiBook} or \Cref{sec:wiener-chaos})
	\begin{equation}\label{eq:hermiso}
		\E\left[H_q(X(x))H_q(X(y))\right]
		=
		q!\,C(x-y)^q
	\end{equation}
	shows that the spectral measure of \(H_q(X)\), which is the random factor of the integrand in $F_{T,q}^{(z)}$, is \(q!G^{*q}\).

    Before passing to the next subsection, we also need to show the separability of $L^2$ spaces on $\wE$. Since $\wE$ is second countable, its Borel
$\sigma$-algebra is countably generated. Hence,
for every finite symmetric\footnote{Symmetric means $G(A)=G(-A)$ for every $A\subseteq \wE$ Borel set.} Borel measure $G$ on $\widehat E$,
the Hilbert space $L^2(G):=L^2(\widehat E,\mathbb C,G)$ is a separable complex Hilbert space. In particular, the closed real subspace $L^2_\R(G)$, defined as
	\[
	L^2_{\mathbb R}(G)
	=
	\left\{
	f\in L^2(G):
	f(-\lambda)=\overline{f(\lambda)}
	\right\}\,,
	\]
	is a real separable Hilbert space. Indeed, since $G$ is symmetric, for \(f,g\in L^2_{\mathbb R}(G)\) we have
	\[
	\langle f,g\rangle_{L^2_{\mathbb R}(G)}=\langle f,g\rangle_{L^2(G)}
	=
	\int_{\widehat E}f(\lambda)\overline{g(\lambda)}\,G(d\lambda)
	=
	\int_{\widehat E}f(\lambda)g(-\lambda)\,G(d\lambda)\in\R\,.
	\]

    We conclude with a remark related to the important examples $E=\R^d,\Z^d,\T^d$.

     \begin{remark}\label{rem:Bochner classic}
The standard identifications
\[
\widehat{\mathbb R^d}\simeq\mathbb R^d,\qquad
\widehat{\mathbb Z^d}\simeq\mathbb T^d,\qquad
\widehat{\mathbb T^d}\simeq\mathbb Z^d
\]
can be realized by the isomorphisms of topological groups
\[
\Phi:F\to\widehat E,\qquad
\Phi(\xi)(x):=e^{ix\cdot\xi},
\]
where $x\cdot\xi$ is the Euclidean
scalar product and $F=\mathbb R^d,\mathbb T^d,\mathbb Z^d$
when $E=\mathbb R^d,\mathbb Z^d,\mathbb T^d$
respectively, see e.g. \cite[Theorem~4.6 and Proposition~4.7]{FollandAHA}.
Thus, a Borel measure $G$ on $\widehat E$ can
be identified with the pushforward measure
$\widetilde G:=(\Phi^{-1})_*G$ on $F$, satisfying $\widetilde G(A)=G(\Phi(A))$ $\forall$ $A\subseteq F$ Borel and
\[
\int_{\widehat E}e^{i\langle \lambda,x\rangle}\,G(d\lambda)=\int_{\widehat E}\lambda(x)\,G(d\lambda)
=
\int_F e^{i x\cdot\xi}
\,\widetilde G(d\xi)\,,\qquad\forall\,\,x\in E\,,
\]
which gives the classical Bochner's theorems when $E=\R^d,\Z^d,\T^d$.
Thus, when working with
$E=\R^d,\Z^d,\T^d$ in
\Cref{sec:highfrequencyapplications,sec:large-domain-applications,sec:proofs},
we identify $\widehat E$ with $F$ and the
spectral measure $G$ with
$\widetilde G$. For simplicity, we
continue to denote them by $\widehat E$ and
$G$, calling them Pontryagin dual of $E$ and spectral measure of $C$, respectively. Moreover, in these cases we use $x\cdot \xi$ and $\langle x,\xi\rangle$ interchangeably to denote the Euclidean scalar product.
\end{remark}

	\subsection{Wiener chaoses and multiple Wiener-It\^o integrals}
	\label{sec:wiener-chaos}
	
	We now introduce Wiener chaoses, multiple Wiener-It\^o integrals and fourth moment theorems, see e.g.
	\cite{NourdinPeccatiBook} for more details. 
	
	Let \(\fH\) be a real separable Hilbert space and
	\(I=\{I(h):h\in\fH\}\) an isonormal Gaussian process over \(\fH\), that is,
	a centered Gaussian field indexed by \(\fH\) satisfying
	\[
	\E[I(h)I(g)]
	=
	\langle h,g\rangle_{\fH},
	\qquad h,g\in\fH.
	\]
	A classical result is the orthogonal chaotic decomposition
	\[
	L^2(\Omega)
	:=
	L^2(\Omega,\sigma(I),\PP)
	=
	\bigoplus_{q\ge0}\mathcal H_q,
	\]
	where \(\mathcal H_0=\mathbb R\) and, for \(q\ge1\), the \(q\)-th Wiener
	chaos is defined by
	\[
	\mathcal H_q
	:=
	\overline{
		\operatorname{span}
		\left\{
		H_q(I(h)):
		h\in\fH,\ \|h\|_{\fH}=1
		\right\}
	}^{\,L^2(\Omega)}.
	\]
	Here the Hermite polynomials $H_q$ can be defined recursively
	\[
	H_0(x)=1,\qquad H_1(x)=x,\qquad
	H_{q+1}(x)=xH_q(x)-qH_{q-1}(x),
	\]
	and form a complete orthogonal system in \(L^2(\mathbb R,\gamma)\) satisfying the isometry relation \eqref{eq:hermiso}.
	
	We now recall the tensor representation of \(\mathcal H_q\) in terms of Wiener-It\^o integrals $I_q$. Let
	\(\fH^{\otimes q}\) be the \(q\)-fold Hilbert tensor product of \(\fH\),
	whose inner product is characterized on elementary tensors by
	\[
	\left\langle
	h_1\otimes\cdots\otimes h_q,
	g_1\otimes\cdots\otimes g_q
	\right\rangle_{\fH^{\otimes q}}
	=
	\prod_{j=1}^q\langle h_j,g_j\rangle_{\fH},
	\]
	and let \(\fH^{\odot q}\) denote its closed subspace of symmetric tensors,
	that is, tensors invariant under permutations of their \(q\) coordinates.
	By polarization identities
	\({\rm span}\left\{h^{\otimes q}:h\in\fH\,,\,\|h\|_{\fH}=1\right\}\) is dense in \(\fH^{\odot q}\). For \(h\in\fH\), $\|h\|_\fH=1$, we set $I_q\left(h^{\otimes q}\right):=H_q(I(h))$.
	We first extend this definition linearly to ${\rm span}\left\{h^{\otimes q}:h\in\fH\,,\,\|h\|_{\fH}=1\right\}$, then by density \(I_q\) extends uniquely to a continuous linear isomorphism
	\[
	I_q:\fH^{\odot q}\rightarrow\mathcal H_q,
	\]
	called the $q$th multiple Wiener--It\^o integral, satisfying the Hermite isometry relation
	\begin{equation*}
		\E\left[I_q(f)I_q(g)\right]
		=
		q!\,\langle f,g\rangle_{\fH^{\otimes q}},
		\qquad f,g\in\fH^{\odot q}.
	\end{equation*}
	Thus, every element in
	\(\mathcal H_q\) admits a unique representation \(I_q(f)\), with
	\(f\in\fH^{\odot q}\).
	
	Combining this representation with the orthogonal decomposition of
	\(L^2(\Omega)\), every \(F\in L^2(\Omega)\) admits a unique chaotic
	decomposition in $L^2(\Omega)$
	\[
	F-\E[F]
	=
	\sum_{q\ge1}I_q(f_q),
	\qquad f_q\in\fH^{\odot q}\,.
	\]
	
	We now specialize this construction to the stationary Gaussian fields of
	this paper. Let \(X_T\) be a sequence of stationary centered Gaussian fields
	on \(E\), with covariance function \(C_T\), $C_T(0)=1$, and spectral measure \(G_T\).
	Consider the real separable Hilbert space
	\(\fH_T=L^2_{\mathbb R}(G_T)\), with
	inner product
	\[
	\langle f,g\rangle_{\fH_T}
	=
	\int_{\widehat E}f(\lambda)\overline{g(\lambda)}\,G_T(d\lambda)
	=
	\int_{\widehat E}f(\lambda)g(-\lambda)\,G_T(d\lambda).
	\]
	Let \(I_T=\{I_T(h):h\in\fH_T\}\) be an isonormal Gaussian process over
	\(\fH_T\). For \(x\in E\), define
	\(e_x(\lambda):=e^{i\langle x,\lambda\rangle}\). Notice that
	\(e_x\in L^2_{\mathbb R}(G_T)\) and \(\|e_x\|_{\fH_T}=1\). Moreover,
	\[
	\E\left[I_T(e_x)I_T(e_y)\right]
	=
	\int_{\widehat E}e^{i\langle x-y,\lambda\rangle}\,G_T(d\lambda)
	=
	C_T(x-y)
	=
	\E\left[X_T(x)X_T(y)\right].
	\]
	Thus \(X_T\) and \((I_T(e_x))_{x\in E}\) are centered Gaussian fields with
	the same covariance function, hence have the same finite-dimensional
	distributions. We may therefore work with the realization
	\[
	X_T(x)=I_T(e_x),\qquad x\in E.
	\]
	
	We denote by \(I_{T,q}\) the corresponding multiple Wiener--It\^o integrals.
	The tensor powers of \(\fH_T\) can then be identified with the corresponding
	\(L^2\)-spaces on the spectral side:
	\[
	\fH_T^{\otimes q}
	\cong
	L^2_{\mathbb R}\left(G_T^{\otimes q}\right),
	\qquad
	\fH_T^{\odot q}
	\cong
	L^2_{\mathbb R,\mathrm{sym}}\left(G_T^{\otimes q}\right),
	\]
	where \(G_T^{\otimes q}\) is the \(q\)-fold product measure and
	\({\rm sym}\) indicates invariance under permutations of the \(q\)
	variables.
	
	We shall also need contractions. 
	Under the spectral identification above, the $r$th contraction of
	\(f\in\fH_T^{\odot p}\) and \(g\in\fH_T^{\odot q}\) is the tensor
	\(f\otimes_r g\in\fH_T^{\otimes(p+q-2r)}\) given  by
	\begin{equation}
		\label{eq:defcontractions}
		\left(f\otimes_r g\right)(u,u')
		=
		\int_{\left(\widehat E\right)^r}
		f(\lambda,u)\,
		g(-\lambda,u')\,
		G_T^{\otimes r}(d\lambda),
	\end{equation}
	where \(u\in(\widehat E)^{p-r}\), \(u'\in(\widehat E)^{q-r}\). 
	We shall repeatedly use the Nualart--Peccati Fourth Moment Theorem \cite{NualartPeccati}: if \(q\ge2\), \(f_T\in\fH_T^{\odot q}\), and $\Var(I_{T,q}(f_T))
	=
	q!\left\|f_T\right\|_{\fH_T^{\otimes q}}^2
	\rightarrow1$, then $I_{T,q}(f_T)\overset{\rm law}{\rightarrow}N(0,1)$ if and only if
	\begin{equation}
		\label{eq:fourth-moment-real}
		\E\left[I_{T,q}(f_T)^4\right]\rightarrow3
		\quad\Leftrightarrow\quad
		\left\|f_T\otimes_r f_T\right\|_{\fH_T^{\otimes2(q-r)}}\rightarrow0\quad \text{for every \(r=1,\ldots,q-1\)}\,.
	\end{equation}
Another result we will use is the Peccati-Tudor multi-variate Fourth-Moment theorem
\cite{PeccatiTudor}: for fixed integers $q_1,\ldots,q_m\ge1$,
if $F_T=(I_{T,q_i}(f_{i,T}))_{i=1}^m$ has covariance
matrices converging to $\Sigma\in\R^{m\times m}$, then
\begin{equation}
    \label{eq:peccatitudor}
    F_T\xrightarrow{\mathrm{law}}\mathcal N_m(0,\Sigma)
\quad\Leftrightarrow\quad
I_{T,q_i}(f_{i,T})
\xrightarrow{\mathrm{law}}\mathcal N(0,\Sigma_{ii})
\quad\text{for every }i.
\end{equation}

	For Fourier functionals at frequencies \(z\neq-z\), we shall also need a complex-valued version of the Fourt Moment theorem. We therefore introduce the
	complexification of $\fH_T$, namely
	\[
	\fH_{T,\mathbb C}
	=
	\fH_T+i\fH_T\,,
	\]
	endowed with the Hermitian inner
	product, linear in the first argument and antilinear in the second,
	\[
	\langle a+ib,u+iv\rangle_{\fH_{T,\mathbb C}}
	=
	\langle a,u\rangle_{\fH_T}
	+
	\langle b,v\rangle_{\fH_T}
	+
	i\big(
	\langle b,u\rangle_{\fH_T}
	-
	\langle a,v\rangle_{\fH_T}
	\big).
	\]
	Since $L^2(G_T)$, as defined in \Cref{sec:lca-spectral}, is the complexification
of $L^2_{\mathbb R}(G_T)$, we also have
	\[
	\fH_{T,\mathbb C}\cong L^2(G_T),
	\qquad
	\langle h,g\rangle_{\fH_{T,\mathbb C}}
	=
	\int_{\widehat E}
	h(\lambda)\overline{g(\lambda)}\,G_T(d\lambda).
	\]
	Similarly,
	\(\fH_{T,\mathbb C}^{\otimes q}\cong L^2\left(G_T^{\otimes q}\right)\) and
	\(\fH_{T,\mathbb C}^{\odot q}\cong
	L^2_{\mathrm{sym}}\left(G_T^{\otimes q}\right)\). It is useful here to distinguish two notions of conjugation. If
	\(h=a+ib\), with \(a,b\in\fH_T^{\odot q}\), we denote by
	\[
	h^\star:=a-ib
	\]
	the conjugation associated with the complexification. Under the spectral
	identification, this is simply
	\[
	h^\star(\lambda)
	=
	\overline{h(-\lambda)}\,,
	\]
	Thus \(h^\star\) should not be confused with the pointwise complex conjugate
	\(\overline h\).
	
	The multiple Wiener--It\^o integral is extended to
	\(\fH_{T,\mathbb C}^{\odot q}\) by complex linearity:
	\(I_{T,q}(a+ib):=I_{T,q}(a)+iI_{T,q}(b)\).
	With this convention,
	\[
	\overline{I_{T,q}(h)}
	=
	I_{T,q}(h^\star).
	\]
	The Wiener--It\^o isometry becomes, for
	\(h,g\in\fH_{T,\mathbb C}^{\odot q}\),
	\[
	\E\!\left[
	I_{T,q}(h)\overline{I_{T,q}(g)}
	\right]
	=
	q!\,
	\langle h,g\rangle_{\fH_{T,\mathbb C}^{\otimes q}}.
	\]
	In particular, its variance and pseudo-variance are
	\begin{align*}
		\E\left[\left|I_{T,q}(h)\right|^2\right]
		&=
		q!\int_{\left(\widehat E\right)^q}|h(\lambda)|^2\,
		G_T^{\otimes q}(d\lambda), \\
		\E\left[I_{T,q}(h)^2\right]
		&=
		q!\int_{\left(\widehat E\right)^q}h(\lambda)h(-\lambda)\,
		G_T^{\otimes q}(d\lambda)=q!\int_{\left(\widehat E\right)^q}h(\lambda)\overline{h^\star(\lambda)}\,
		G_T^{\otimes q}(d\lambda).
	\end{align*}
	The contractions are extended to the complexified tensor spaces by
	complex bilinearity. With this convention,
	the spectral formula \eqref{eq:defcontractions} remains unchanged.
	
	We can finally state the complex version of the Fourth Moment criterion, see, e.g. \cite{MarinucciPeccatiCompact}.
	
	\begin{proposition}[Complex Fourth Moment criterion]
		\label{prop:complex-fmt}
		Let \(q\ge2\), \(f_T\in\fH_{T,\mathbb C}^{\odot q}\), and set
		\(F_T=I_{T,q}(f_T)\). Assume that
		\[
		\E\left|F_T\right|^2\rightarrow1,
		\qquad
		\E\left[F_T^2\right]\rightarrow0.
		\]
		Let
		\[
		N_{\mathbb C}
		=
		\frac{N_1+iN_2}{\sqrt2},
		\qquad
		N_1,N_2\overset{\mathrm{i.i.d.}}{\sim}\mathcal N(0,1).
		\]
		Then
		\[
		F_T\overset{\rm law}{\rightarrow}N_{\mathbb C}
		\quad\Leftrightarrow\quad
		\left\|f_T\otimes_r f_T^\star\right\|_
		{\fH_{T,\mathbb C}^{\otimes2(q-r)}}
		\rightarrow0,
		\qquad r=1,\ldots,q-1.
		\]
	\end{proposition}
	\begin{proof}
		If \(f_T=a_T+ib_T\), with \(a_T,b_T\in\fH_T^{\odot q}\), then $A_T:={\rm Re}(F_T)=I_{T,q}(a_T)$, $B_T:={\rm Im}(F_T)=I_{T,q}(b_T)$. The proof strategy can be found e.g. in \cite{MarinucciPeccatiCompact}, here we provide it in our notation for completeness and because we use a different isonormal representation and notation. Since by assumption
		\[
		\E\left[F_T^2\right]
		=
		\E\left[A_T^2\right]-\E\left[B_T^2\right]+2i\E\left[A_TB_T\right]\rightarrow 0,
		\]
		we get
		\[
		\E\left[A_T^2\right]\to\frac12,\qquad
		\E\left[B_T^2\right]\to\frac12,\qquad
		\E\left[A_TB_T\right]\to0 
		\]
		Hence, by \eqref{eq:peccatitudor} with $\Sigma=\frac{1}{2}I_2\in \R^{2\times 2}$ and by \eqref{eq:fourth-moment-real}, we have \(F_T\Rightarrow N_{\mathbb C}\)
		if and only if, for every \(r=1,\ldots,q-1\),
		\[
		\left\|a_T\otimes_r a_T\right\|\to0,
		\qquad
		\left\|b_T\otimes_r b_T\right\|\to0.
		\]
		Moreover, by bilinearity,
		\[
		f_T\otimes_r f_T^\star
		=
		a_T\otimes_r a_T+b_T\otimes_r b_T
		+i\left(b_T\otimes_r a_T-a_T\otimes_r b_T\right).
		\]
Using the standard identity for symmetric tensors
\[
\left\langle a_T\otimes_r a_T,
b_T\otimes_r b_T\right\rangle
=\left\|a_T\otimes_{q-r}b_T\right\|^2\ge0,
\]
together with Cauchy--Schwarz and
$\left\|a_T\otimes_r b_T\right\|
=\left\|b_T\otimes_r a_T\right\|$, we obtain
\[
\begin{aligned}
\left\|a_T\otimes_r a_T\right\|^2
+\left\|b_T\otimes_r b_T\right\|^2\le &\left\|f_T\otimes_r f_T^\star\right\|^2\le 
\\
&\le 2\left(\left\|a_T\otimes_r a_T\right\|^2
+\left\|b_T\otimes_r b_T\right\|^2\right)+
4\left\|a_T\otimes_{q-r}a_T\right\|
\left\|b_T\otimes_{q-r}b_T\right\|.
\end{aligned}
\]
	\end{proof}
	
	\section{Main result in the general setting}\label{sec:mainresult}
	
	Let \(E\) be a second countable Hausdorff LCA group, fix a Haar measure
	\(\mu\), and endow \(\widehat E\) with a compatible translation-invariant metric $d_{\wE}$. Let
	\(B_\varepsilon(x)\) denotes the corresponding open ball centered in $x\in\wE$, satisfying $B_\eps(x)=x+B_\eps(0)$, see \Cref{sec:lca-spectral} for the details.
    
    Let
	\((X_T)_T\) be stationary Gaussian fields with spectral measures \(G_T\),
	and let \(\psi_T\in L^1(E,\mu)\) be real-valued. For \(q\ge2\), define
	\[
	F_{T,q}^{(z)}
	:=\int_EH_q(X_T(x))\psi_T(x)e^{-i\langle x,z\rangle}\,\mu(dx).
	\]
	The definition and the continuity of Wiener-It\^o integrals give
	\begin{equation*}
		F_{T,q}^{(z)}
		=\int_E I_{T,q}\left(e^{i\langle \lambda_1+\cdots+\lambda_q, x\rangle}\right)\psi_T(x)e^{-i\langle x,z\rangle}\,\mu(dx) =I_{T,q}\left(\wpsi_T(\lambda_1+\cdots+\lambda_q-z)\right).
	\end{equation*}
	Consequently,
	\begin{align}
		\Var\left(F_{T,q}^{(z)}\right)
		&=q!V_{T,q}(z),\label{eq:variance-kernel}\\
		V_{T,q}(z)
		&:=\int_{\widehat E}\left|\wpsi_T(\lambda-z)\right|^2G_T^{*q}(d\lambda),\notag
	\end{align}
	and
	\begin{equation}
		\label{eq:pseudovariance-kernel}
		\E\left[\left(F_{T,q}^{(z)}\right)^2\right]
		=q!\int_{\widehat E}\wpsi_T(\lambda-z)
		\wpsi_T(-\lambda-z)\,G_T^{*q}(d\lambda).
	\end{equation}
	For \(\varepsilon>0\), set
	\begin{align*}
		\mathcal D_{T,q,z}(\varepsilon)
		&:=G_T^{*q}(B_\varepsilon(z)),\\
		\mathcal A_{T,r,q,z}(\varepsilon)
		&:=\sup_{\alpha\in\widehat E}
		G_T^{*r}(B_\varepsilon(\alpha))
		G_T^{*(q-r)}(B_\varepsilon(z-\alpha)),
		\quad 1\le r<q.
	\end{align*}
	Let \(\Lambda_{T,1},\ldots,\Lambda_{T,q}\) i.i.d. \(\sim G_T\)
	and set \(S_{T,r}=\Lambda_{T,1}+\cdots+\Lambda_{T,r}\). Whenever
	\(\mathcal D_{T,q,z}(\varepsilon)>0\), define
	\begin{equation}
		\eta_{T,r,z}(\varepsilon)
		:=\sup_{\alpha\in\widehat E}
		\PP\!\left(
		S_{T,r}\in B_\varepsilon(\alpha)
		\,\middle|\,
		S_{T,q}\in B_\varepsilon(z)
		\right).\label{eq:condprobproof}
	\end{equation}
	
	The main goal of this section is to prove that under suitable assumptions the Gaussian convergence $\widetilde F_{T,q}^{(z)}\overset{\rm law}{\rightarrow}N_z$ is equivalent to the spectral delocalization property discussed in \Cref{subsec:main-nutshell}, see \Cref{thm:main}-(iii). Let us give the formal definition.
	
	\begin{definition}[Spectral delocalization]
		\label{def:resonant-delocalization}
		Assume \(z\in\supp(G_T^{*q})\) for all sufficiently large \(T\). We say
		that the family \((G_T)_T\) satisfies \emph{spectral delocalization} at
		order \(q\) and frequency \(z\) if for some \(M>0\),
		\[
		\eta_{T,r,z}(M/T)\rightarrow0\,,\quad \forall r=1,\dots,q-1.
		\]
	\end{definition}
	
	The previous definition corresponds to the conditional formulation \eqref{eq:condform}. The fact that \eqref{eq:condform} is equivalent to the convolution formulation \eqref{eq:convform} is proved in 
	\Cref{prop:conditional-product}.\\
	
	We will impose three assumptions. The first one is a spectral localization assumption given simultaneously on the sequence of test functions and spectral measures. It allows to control variances and pseudovariances in terms of small probabilities of $G_T^{*q}$.
	
	\begin{condition}[Localization at \(z\)]\label{cond:speclocalization}
		We have \(V_{T,q}(z)>0\) for all sufficiently large \(T\), and
		\[
		\lim_{M\to\infty}\limsup_{T\to\infty}
		\frac{\displaystyle
			\int_{B_{M/T}^c(z)}\left|\wpsi_T(\lambda-z)\right|^2G_T^{*q}(d\lambda)}
		{V_{T,q}(z)}=0.
		\]
	\end{condition}
	The second condition is given only on test functions and allows to control their Fourier transform from above and from below.
	\begin{condition}[Global and local Fourier control]
		\label{cond:upperelowerpsi}
		There exist \(a,c_\psi>0\) and \(C_\psi<\infty\) such that, for all
		sufficiently large \(T\),
		\[
		\wpsi_T(0)>0,
		\qquad
		\left|\wpsi_T(\theta)\right|\le C_\psi\wpsi_T(0),
		\quad \theta\in\widehat E,
		\]
		and
		\[
		\inf_{\theta\in B_{a/T}(0)}
		\operatorname{Re}\frac{\wpsi_T(\theta)}{\wpsi_T(0)}
		\ge c_\psi.
		\]
	\end{condition}
	
	The third, final assumption is purely spectral and imposes a control of $G_T^{*q}$ near $z$.
	
	\begin{condition}[Local doubling at the probed frequency]
		\label{cond:doubling}
		We have $z\in\supp(G_T^{*q})$ for all sufficiently large $T$. Moreover, there exist
		$\varepsilon_0>0$ and $C_D>1$, independent of $T$, such that
		\[
		G_T^{*q}(B_{2\varepsilon}(z))
		\le C_D\,G_T^{*q}(B_\varepsilon(z)),
		\qquad 0<\varepsilon<\varepsilon_0,
		\]
		for all sufficiently large $T$.
	\end{condition}
	
	\begin{remark}\label{rem:profile-doubling}
		Under \Cref{cond:doubling}, for every fixed $b,c>0$, and uniformly for $T$ large enough,
		\[
		\mathcal D_{T,q,z}(b/T)\asymp_{b,c}\mathcal D_{T,q,z}(c/T)
		\]
	\end{remark}
	
	Under \Cref{cond:doubling}, the conditions \eqref{eq:condform} and \eqref{eq:convform} are equivalent. 
	\begin{proposition}[Equivalence of \eqref{eq:condform} and \eqref{eq:convform}]
		\label{prop:conditional-product}
		For every \(\varepsilon>0\) such that
		\(\mathcal D_{T,q,z}(\varepsilon)>0\),
		\begin{equation*}
			\frac{\mathcal A_{T,r,q,z}(\varepsilon/2)}
			{\mathcal D_{T,q,z}(\varepsilon)}
			\le
			\eta_{T,r,z}(\varepsilon)
			\le
			\frac{\mathcal A_{T,r,q,z}(2\varepsilon)}
			{\mathcal D_{T,q,z}(\varepsilon)}.
		\end{equation*}
		If \Cref{cond:doubling} holds, spectral delocalization is
		equivalent to the following: for some $M>0$ 
		\[
		\frac{\mathcal A_{T,r,q,z}(M/T)}
		{\mathcal D_{T,q,z}(M/T)}\rightarrow0\,,\qquad\forall\, r=1,\dots,q-1\,.
		\]
	\end{proposition}
	\begin{proof}
		For fixed \(\alpha\), by \eqref{eq:condprobproof} and translation invariance, the numerator of the conditional probability is
		\[
		\int_{B_\varepsilon(\alpha)}
		G_T^{*(q-r)}(B_\varepsilon(z-x))\,G_T^{*r}(dx).
		\]
		If \(x\in B_\varepsilon(\alpha)\), then
		\(B_\varepsilon(z-x)\subseteq B_{2\varepsilon}(z-\alpha)\), which gives
		the upper bound. Restricting the integral to
		\(B_{\varepsilon/2}(\alpha)\) and using
		\(B_{\varepsilon/2}(z-\alpha)\subseteq B_\varepsilon(z-x)\) gives the
		lower bound. The equivalence follows from \Cref{cond:doubling}.
	\end{proof}

	The next result allows to control precisely the variance rates and the pseudovariance vanishing for our family of Fourier-Hermite functionals.
	\begin{proposition}\label{prop:varcontrol}
		Assume that \Cref{cond:speclocalization}, \Cref{cond:upperelowerpsi} and \Cref{cond:doubling} hold. Then
		$\exists$ \(C>1\) such that
		\[
		C^{-1}\wpsi_T(0)^2\mathcal D_{T,q,z}(1/T)
		\le V_{T,q}(z)
		\le C\,\wpsi_T(0)^2\mathcal D_{T,q,z}(1/T).
		\]
		Moreover, if \(z\ne-z\), then
		\[
		\frac{\E\left[\left(F_{T,q}^{(z)}\right)^2\right]}{\Var\left(F_{T,q}^{(z)}\right)}\rightarrow0.
		\]
	\end{proposition}
	Finally, to apply the classical and complex Fourth Moment theorems, we need a control on contractions norms. For
	\[
	f_T(\lambda_1,\ldots,\lambda_q)
	:=\wpsi_T(\lambda_1+\cdots+\lambda_q-z),
	\]
	the relevant complex contraction in \Cref{prop:complex-fmt} is
	\[
	\left(f_T\otimes_rf_T^\star\right)(u,u')
	=\int_{\left(\widehat E\right)^r}
	f_T(\lambda,u)\overline{f_T(\lambda,-u')}
	\,G_T^{\otimes r}(d\lambda)\,,
	\]
	with squared norm
	\begin{equation}\label{eq:conorm}
		\left\|f_T\otimes_rf_T^\star\right\|^2_{\fH_{T,\C}^{\otimes 2(q-r)}}=\int_{\left(\widehat E\right)^{2q}}
		f_T(\lambda,u)\overline{f_T(\lambda,-u')}\,\overline{f_T(\lambda',u)}f_T(\lambda',-u')
		\,dG_T^{\otimes 2q}(u,u',\lambda,\lambda')\,.
	\end{equation}
	For $M>0$ and $1\le r<q$, write
	\[
	\beta_{T,r}(M)
	:=\frac{\mathcal A_{T,r,q,z}(M/T)}{\mathcal D_{T,q,z}(M/T)},
	\qquad
	\kappa_{T,r}
	:=\frac{\left\|f_T\otimes_rf_T^\star\right\|^2}{V_{T,q}(z)^2},
	\]
	and set
	\[
	L_{T,M}:=\int_{B_{M/T}^c(z)}
	\left|\wpsi_T(\lambda-z)\right|^2G_T^{*q}(d\lambda).
	\]
	
	\begin{theorem}[Contractions and localization scales]
		\label{thm:contractions-delocalization}
		Assume that \Cref{cond:speclocalization}, \Cref{cond:upperelowerpsi} and \Cref{cond:doubling} hold. Let $a$ be the radius in \Cref{cond:upperelowerpsi}.
		For every fixed $M>0$ and $0<\delta<a/2$, there are constants
		$c_\delta>0$ and $C_M<\infty$, independent of $T$, such that, for all sufficiently large $T$,
		\begin{equation}
			\label{eq:scale-explicit-bounds}
			c_\delta\beta_{T,r}(\delta)^2
			\le\kappa_{T,r}
			\le C_M\beta_{T,r}(3M)
			+4\left(\frac{L_{T,M}}{V_{T,q}(z)}\right)^{1/2},
			\qquad 1\le r<q\,.
		\end{equation}
	\end{theorem}
	Combining the previous result with the real and complex Fourth Moment theorems in \Cref{sec:wiener-chaos}, we can finally obtain the main result of the paper.
	\begin{theorem}[Spectral criterion]
		\label{thm:main}
		Fix $q\ge2$, $z\in\widehat E$, and assume that \Cref{cond:speclocalization}, \Cref{cond:upperelowerpsi} and \Cref{cond:doubling} hold. Let $a$ be the radius in \Cref{cond:upperelowerpsi}. Then we have:
		\begin{enumerate}[label=\textnormal{(\roman*)}]
			\item If $\beta_{T,r}(M)\rightarrow 0$ $\forall\, r=1,\dots,q-1$ and for every fixed $M>0$, then
			$\widetilde F_{T,q}^{(z)}\overset{\rm law}{\rightarrow} N_z$.
			\item If $\widetilde F_{T,q}^{(z)}\overset{\rm law}{\rightarrow} N_z$, then
			$\beta_{T,r}(\delta)\rightarrow 0$ for every
			$0<\delta<a/2$ and every $r=1,\dots,q-1$.
			\item If
			$\widehat E$ is discrete, or $G_T=G$ is fixed, then $\beta_{T,r}(M)\rightarrow 0$ $\forall M>0$ $\;\Leftrightarrow\;$ $\beta_{T,r}(\delta)\rightarrow 0$ for some $\delta>0$. Thus, $\widetilde F_{T,q}^{(z)}\overset{\rm law}{\rightarrow} N_z$ $\;\Leftrightarrow\;$ $\beta_{T,r}(\delta)\rightarrow 0$ $\forall r=1,\dots,q-1$, for some $\delta>0$.
		\end{enumerate}
	\end{theorem}
	
	\section{Proofs}\label{sec:proofs}
	
	\subsection{Proof of the main result}
	\begin{proof}[Proof of \Cref{prop:varcontrol}]
		By \Cref{cond:upperelowerpsi} and
		\Cref{cond:doubling}-\Cref{rem:profile-doubling},
		\[
		V_{T,q}(z)
		\ge c_\psi^2\wpsi_T(0)^2\mathcal D_{T,q,z}(a/T)
		\gtrsim \wpsi_T(0)^2\mathcal D_{T,q,z}(1/T).
		\]
		By \Cref{cond:speclocalization}, one can choose a fixed \(M\) such that, for all sufficiently
		large \(T\)
		\[
		L_{T,M}=\int_{B_{M/T}^c(z)}
		\left|\wpsi_T(\lambda-z)\right|^2G_T^{*q}(d\lambda)
		\le\frac12V_{T,q}(z).
		\]
		\Cref{cond:upperelowerpsi} and \Cref{cond:doubling} then give
		\[
		V_{T,q}(z)
		\le2C_\psi^2\wpsi_T(0)^2\mathcal D_{T,q,z}(M/T)
		\lesssim\wpsi_T(0)^2\mathcal D_{T,q,z}(1/T).
		\]
		Suppose now \(z\ne-z\). For fixed \(M\), the balls \(B_{M/T}(z)\) and
		\(B_{M/T}(-z)\) are disjoint for large \(T\). Splitting the integral in
		\eqref{eq:pseudovariance-kernel} over \(B_{M/T}(z)\) and its complement,
		then using symmetry of \(G_T^{*q}\) and Cauchy--Schwarz, yields
		\[
		\left|\E\left[\left(F_{T,q}^{(z)}\right)^2\right]\right|
		\le2q!V_{T,q}(z)^{1/2}L_{T,M}^{1/2}.
		\]
		We conclude dividing by \(q!V_{T,q}(z)\), letting first \(T\to\infty\) and then
		\(M\to\infty\), using \Cref{cond:speclocalization}.
	\end{proof}
	
	Before proceeding, we need a technical lemma that allows to control spectral overlap probabilities in terms of $\mathcal{D}_{T,q,z}$ and $\mathcal{A}_{T,r,q,z}$.
	\begin{lemma}[Contractions and overlap probabilities]
		\label{lem:fourfold-overlap}
		Fix $1\le r<q$. Let $\Lambda,\Lambda'\sim G_T^{*r}$ and
		$U,U'\sim G_T^{*(q-r)}$ be all independent, and let
		$\Gamma_T(\varepsilon)$ be the overlap event 
		\begin{equation}
		    \label{eq:overlapevent}
            \Gamma_T(\eps):=\{\Lambda+U,\Lambda'+U',\Lambda'+U,\Lambda +U'\in B_\eps(z)\}\,.
		\end{equation}
		For every $\varepsilon>0$,
		\begin{equation*}
			\mathcal A_{T,r,q,z}(\varepsilon/2)^2
			\le\PP(\Gamma_T(\varepsilon))
			\le\mathcal D_{T,q,z}(\varepsilon)
			\mathcal A_{T,r,q,z}(3\varepsilon).
		\end{equation*}
	\end{lemma}
	
	\begin{proof}
		On the overlap event,
		$\Lambda'\in B_{2\varepsilon}(\Lambda)$ and
		$U'\in B_{2\varepsilon}(U)$. Moreover, if
		$\Lambda,\Lambda'\in B_{\varepsilon/2}(\alpha)$ and
		$U,U'\in B_{\varepsilon/2}(z-\alpha)$, all four sums are in
		$B_\varepsilon(z)$. Thus, we get
        \[
        \{\Lambda,\Lambda'\in B_{\varepsilon/2}(\alpha),U,U'\in B_{\varepsilon/2}(z-\alpha)\}\subseteq \Gamma_T(\eps)\subseteq\{\Lambda+U\in B_\eps(z), \Lambda\in B_{2\eps}(\Lambda'), U\in B_{2\eps}(U')\}\,.
        \]
        Thus, the lower bound follows immediately by independence, taking the supremum over $\alpha\in\wE$. Moreover, conditioning on $(\Lambda,U)$,
		\[
		\PP(\Gamma_T(\varepsilon))
		\le\E\left[\mathbf1_{\{\Lambda+U\in B_\varepsilon(z)\}}
		G_T^{*r}(B_{2\varepsilon}(\Lambda))G_T^{*(q-r)}(B_{2\varepsilon}(U))\right].
		\]
		On the event inside the expectation,
		$B_{2\varepsilon}(U)\subset B_{3\varepsilon}(z-\Lambda)$, so the product
		is at most $\mathcal A_{T,r,q,z}(3\varepsilon)$, which proves the upper bound.
	\end{proof}
	
	\begin{proof}[Proof of \Cref{thm:contractions-delocalization}]
		For the independent variables of \Cref{lem:fourfold-overlap}, set
		\[
		X_1=\Lambda+U,\quad X_2=\Lambda'+U',\quad
		X_3=\Lambda'+U,\quad X_4=\Lambda+U'.
		\]
		The contraction identity \eqref{eq:conorm} and the change $u'\leftrightarrow-u'$ (by symmetry of $G_T$) yields
		\[
		\left\|f_T\otimes_rf_T^\star\right\|^2
		=\E\left[\wpsi_T(X_1-z)\wpsi_T(X_2-z)
		\overline{\wpsi_T(X_3-z)\wpsi_T(X_4-z)}\right].
		\]
		Let $\Gamma(M/T)$ be defined as in \eqref{eq:overlapevent}.
		By Cauchy-Schwarz, using that \(X_1\), \(X_2\) are independent, \(X_3\), \(X_4\) are independent, and the fact $\Gamma(M/T)^c= \cup_{i=1}^4\{X_i\in B_{M/T}(z)^c\}$, we have the bound
		\[
		\E\left[\left|\wpsi_T(X_1-z)\wpsi_T(X_2-z)
		\overline{\wpsi_T(X_3-z)\wpsi_T(X_4-z)}\mathbf{1}_{\Gamma(M/T)^c}\right|\right]\le 4V_{T,q}(z)^{3/2}L_{T,M}^{1/2}\,.
		\]
		On the other hand, by \Cref{cond:upperelowerpsi} on $\Gamma(M/T)$ the integrand is at most $C_\psi^4\wpsi_T(0)^4$.
		Thus \Cref{lem:fourfold-overlap} gives the bound
		\begin{equation*}
			\left\|f_T\otimes_rf_T^\star\right\|^2
			\le C_\psi^4\wpsi_T(0)^4\mathcal D_{T,q,z}(M/T)
			\mathcal A_{T,r,q,z}(3M/T)
			+4V_{T,q}(z)^{3/2}L_{T,M}^{1/2}.
		\end{equation*}
		The proof of the upper bound in \eqref{eq:scale-explicit-bounds} is therefore concluded dividing by $V_{T,q}(z)^2$ and using
		\Cref{prop:varcontrol} with $1/T$, $M/T$, and $3M/T$. 
For the lower bound, fix $0<\delta<a/2$ and
$\alpha\in\widehat E$. If
$G_T^{*(q-r)}(B_{\delta/T}(z-\alpha))=0$,
there is nothing to prove. Otherwise,
conditioning on $(U,U')$, by Cauchy--Schwarz
\[
\begin{aligned}
\left\|f_T\otimes_r f_T^\star\right\|^2
&=\E\left|
\E\left[
\wpsi_T(X_1-z)\overline{\wpsi_T(X_4-z)}
\mid U,U'
\right]\right|^2\\
&\ge
\frac{
\left|
\E\left[
\wpsi_T(X_1-z)\overline{\wpsi_T(X_4-z)}
\mathbf1_{\{U,U'\in B_{\delta/T}(z-\alpha)\}}
\right]\right|^2
}{
G_T^{*(q-r)}(B_{\delta/T}(z-\alpha))^2
}.
\end{aligned}
\]
By independence, the expectation in the numerator
is
\[
\int_{\widehat E}
\left|
\int_{B_{\delta/T}(z-\alpha)}
\wpsi_T(\lambda+u-z)\,
G_T^{*(q-r)}(du)
\right|^2
G_T^{*r}(d\lambda).
\]
For $\lambda\in B_{\delta/T}(\alpha)$ and
$u\in B_{\delta/T}(z-\alpha)$, we have
\[
\lambda+u-z\in B_{2\delta/T}(0)
\subseteq B_{a/T}(0).
\]
Hence, by \Cref{cond:upperelowerpsi},
the last integral is bounded below by
\[
c_\psi^2\wpsi_T(0)^2
G_T^{*r}(B_{\delta/T}(\alpha))
G_T^{*(q-r)}(B_{\delta/T}(z-\alpha))^2.
\]
Consequently,
\[
\begin{aligned}
\left\|f_T\otimes_r f_T^\star\right\|^2
\ge c_\psi^4\wpsi_T(0)^4
\bigl[
G_T^{*r}(B_{\delta/T}(\alpha))
G_T^{*(q-r)}(B_{\delta/T}(z-\alpha))
\bigr]^2.
\end{aligned}
\]
Taking the supremum over $\alpha\in\widehat E$ yields
\begin{equation*}
\left\|f_T\otimes_r f_T^\star\right\|^2
\ge c_\psi^4\wpsi_T(0)^4
\mathcal A_{T,r,q,z}(\delta/T)^2.
\end{equation*}
Dividing by
$V_{T,q}(z)^2$ and applying
\Cref{prop:varcontrol,cond:doubling}
completes the proof.
	\end{proof}
We are now ready to prove our main result.
	\begin{proof}[Proof of \Cref{thm:main}]
    
Since $\widetilde F_{T,q}^{(z)}
=I_{T,q}\left(\widetilde f_T\right)$, where
$\widetilde f_T=f_T/\sqrt{q!V_{T,q}(z)}$,
by \Cref{prop:varcontrol} we have $\E\left[\left|\widetilde F_{T,q}^{(z)}\right|^2\right]=1$ and $\E\left[\left(\widetilde F_{T,q}^{(z)}\right)^2\right]\rightarrow 0$
when $z\ne-z$.
Moreover, $(q!)^2\left\|\widetilde f_T\otimes_r\widetilde f_T^\star\right\|^2=\kappa_{T,r}$.
Thus, the real Fourth Moment Theorem
\eqref{eq:fourth-moment-real} and its complex version
\Cref{prop:complex-fmt} yield
\[
\widetilde F_{T,q}^{(z)}
\xrightarrow{\mathrm{law}}N_z
\quad\Leftrightarrow\quad
\kappa_{T,r}\to0
\quad\text{for every }r=1,\dots,q-1\,.
\]
		If $\beta_{T,r}(M)\to0$ for every fixed $M$, taking first
		$T\to\infty$ at fixed $M$ in the upper bound of
		\eqref{eq:scale-explicit-bounds} and then $M\to\infty$, by \Cref{cond:speclocalization} we get \textnormal{(i)}.
		
		Now fix \(0<\delta<a/2\), with $a$ as in \Cref{cond:upperelowerpsi}. If
		\(\widetilde F_{T,q}^{(z)}\overset{\rm law}{\rightarrow}N_z\), the lower bound in \eqref{eq:scale-explicit-bounds} yields (ii).
        
		Finally, let us prove (iii).
		If $\widehat E$ is discrete, i.e. its singletons are open, then for $\eps_0>0$ sufficiently small $B_{\eps_0}(0)=\{0\}$ where $0=0_{\widehat{E}}$.
		Translation invariance of the considered metric on $\wE$ then gives $B_{M/T}(\alpha)=\{\alpha\}$ for every
		fixed $M>0$ and for all $T\ge M/\eps_0$. Therefore, we have
		\begin{equation*}
			\beta_{T,r}(M)\rightarrow 0\quad \;\Leftrightarrow\; \quad \sup_{\alpha \in \widehat E} \frac{G_T^{*r}(\alpha)G_T^{*(q-r)}(z-\alpha)}{G_T^{*q}(z)}\rightarrow 0\,\,\,\quad \forall\,\, 1\le r<q\,.
		\end{equation*}
		Suppose now $G_T=G$. Then, for every fixed $\delta>0$,
		\begin{equation*}
			\beta_{T,r}(\delta)\rightarrow 0\qquad \;\Leftrightarrow\; \qquad \lim_{\eps \rightarrow 0}\sup_{\alpha \in \widehat E} \frac{G^{*r}(B_\eps(\alpha))G^{*(q-r)}(B_\eps(z-\alpha))}{G^{*q}(B_{\eps}(z))}= 0\,\,\,\quad \forall\,\, 1\le r<q\,.
		\end{equation*}
		Thus, in both cases all
		fixed scales give the same condition, proving
		\textnormal{(iii)}.
	\end{proof}

\subsection{Proof of the criterion in the specific regimes}
In this subsection we prove \Cref{thm:mainTd} and \Cref{thm:mainRdZd}. To do this, we only need to prove variance estimates when $z\not\in\supp(G^{*q})$ and check that the conditions of \Cref{thm:main} hold.
Let us start with the regime \ref{high frequency}, where the proof is immediate.
\begin{proof}[Proof of \Cref{thm:mainTd}]
Here $E=\T^d$, $\widehat E=\Z^d$, and
$\psi_T\equiv1$, so
$\wpsi_T=\mathbf1_{\{0\}}$ and
$V_{T,q}(z)=G_T^{*q}(\{z\})$.
\Cref{cond:speclocalization} and \Cref{cond:upperelowerpsi} are immediate,
while \Cref{cond:doubling} holds since sufficiently small
balls are singletons. We conclude by
\Cref{thm:main}\textnormal{(iii)}.
\end{proof}

Let us move to the regime \ref{large domain}.
We start with a technical lemma that is needed to pass from $\Z^d$ to $\R^d$ by Riemann-sum approximation.

\begin{lemma}
\label{lem:Riemann-sum}
Let $f:\R^d\to\R$ be continuous almost everywhere
and assume $|f(x)|\lesssim\left(1+\|x\|\right)^{-\alpha}$
for some $\alpha>d$. Then, as $\varepsilon\to0$, uniformly for $u$ on compact sets of $\R^d$,
\[
\varepsilon^d\sum_{n\in\Z^d}
f(\varepsilon n)e^{i\langle\varepsilon n,u\rangle}
\rightarrow\widehat f(u)\,.
\]
\end{lemma}

\begin{proof}
For $0<\varepsilon<1$, define the step function $f_\eps:\R^d\rightarrow\R$ by 
$f_\varepsilon(x)=f(\varepsilon n)$ whenever
$x\in\varepsilon(n+[0,1)^d)$.
Since $\|x-\varepsilon n\|\le\sqrt d\,\varepsilon$ for $x\in\varepsilon(n+[0,1)^d)$,
the decay assumption gives
$|f_\varepsilon(x)|\lesssim\left(1+\|x\|\right)^{-\alpha}$, uniformly in  $\varepsilon$. Moreover, since $f_\varepsilon(x)\to f(x)$ at every
continuity point of $f$ and $f$ is a.e. continuous, dominated convergence
gives $\|f_\varepsilon-f\|_1\to0$. Finally, for $M>0$, adding and subtracting $\widehat f_\eps(u)$,
\[
\begin{aligned}
&\sup_{\|u\|\le M}\left|
\varepsilon^d\sum_n
f(\varepsilon n)e^{i\langle\varepsilon n,u\rangle}
-\widehat f(u)\right|\le
\sqrt d\,M\varepsilon\left\|f_\varepsilon\right\|_1+\left\|f_\varepsilon-f\right\|_1
\rightarrow0,
\end{aligned}
\]
where we used
$|e^{i\langle\varepsilon n,u\rangle}
-e^{i\langle x,u\rangle}|
\le \|u\|\|x-\varepsilon n\|\le M\sqrt d\,\varepsilon$
on each cube. This completes the proof.
\end{proof}

We also prove the following lemma, that allows to check \Cref{cond:upperelowerpsi} in all the cases of interest.

\begin{lemma}
\label{lem:window-conditions}
Let $E=\R^d$ or $\Z^d$ and
$\psi_T=\psi(\cdot/T)$, where $\psi:\R^d\rightarrow\R$ is continuous a.e. and satisfies $|\psi(x)|\lesssim \left(1+\|x\|\right)^{-\alpha}$, $\alpha>d$, and $\wpsi(0)>0$.
Then \Cref{cond:upperelowerpsi} holds.
\end{lemma}
\begin{proof}
By change of variable if $E=\R^d$, and by
\Cref{lem:Riemann-sum} applied to $\psi$
and $|\psi|$ if $E=\Z^d$, we have
\[
T^{-d}\wpsi_T(u/T)\to\wpsi(u)
\quad\text{locally uniformly in }u,
\qquad
T^{-d}\left\|\psi_T\right\|_1\to\|\psi\|_1.
\]
In particular,
$\wpsi_T(0)\sim T^d\wpsi(0)>0$ and
$|\wpsi_T(\theta)|\le\|\psi_T\|_1
\lesssim\wpsi_T(0)$ uniformly in $\theta$.
Finally, since $\wpsi$ is continuous and
$\wpsi(0)>0$, the local uniform convergence and the continuity of $\wpsi:\R^d\rightarrow \C$
imply that, for some $a>0$, the lower bound in \Cref{cond:upperelowerpsi} holds.
\end{proof}

Finally, when $E=\Z^d$ we use Poisson's summation formula to derive a decay for $\wpsi_T$.

\begin{lemma}
\label{lem:lattice-poisson}
Let $E=\R^d$ or $\Z^d$, $\psi\in \mathcal S(\R^d)$ and
$\psi_T(x)=\psi(x/T)$, $x\in E$.
Then, for every $N>0$,
\[
\left|\wpsi_T(\theta)\right|
\lesssim_N T^d
\left(1+T\,d_{\wE}(\theta,0)\right)^{-N},
\qquad \theta\in\wE\,.
\]
\end{lemma}
\begin{proof}If $E=\R^d$ then the result immediately follows by a change of variable. If $E=\Z^d$, recall $\wE=\T^d$ and representing $\T^d$ as $[-\pi,\pi)^d$ modulo $2\pi$ one has $d_{\wE}(0,\theta)=\inf_{k\in\Z^d}\|\theta+2\pi k\|$.
By Poisson's summation formula (see e.g. \cite[p.244]{Lang1993}), recalling our definition \eqref{eq:Fourierdef} for the Fourier transform, we have $\wpsi_T(\theta)
=T^d\sum_{k\in\Z^d}
\wpsi(T(\theta+2\pi k))$.
Moreover, since $\wpsi$ is Schwartz, we get
\[
\left|\wpsi_T(\theta)\right|
\lesssim_L T^d
\sum_{k\in\Z^d}
\left(1+T\|\theta+2\pi k\|\right)^{-L}\lesssim_L T^d
\left(1+T\,d_{\T^d}(\theta,0)\right)^{-N}\sum_{k\in\Z^d}
\left(1+\|\theta+2\pi k\|\right)^{-(L-N)}.
\]
where we used $d_{\T^d}(\theta,0)\le
\|\theta+2\pi k\|$, $\forall$ $k$. Taking $L>N+d$, the last series is uniformly bounded in $\theta$.
\end{proof}

\begin{proof}[Proof of \Cref{thm:mainRdZd}]
First of all, if $z\notin\supp(G^{*q})$, setting
$\rho=d_{\widehat E}(z,\supp(G^{*q}))>0$, by \Cref{cond:schwartz}, \eqref{eq:variance-kernel} and
\Cref{lem:lattice-poisson} we get $V_{T,q}(z)\lesssim_N
T^{2d}(1+T\rho)^{-2N}$, $\forall$ $N\in \N$. Let us now assume $z\in \supp(G^{*q})$ and check \Cref{cond:speclocalization}, \Cref{cond:upperelowerpsi} and \Cref{cond:doubling}.
\Cref{cond:doubling} is immediate by \Cref{cond:doublingintro} in \Cref{thm:mainRdZd}. By \Cref{lem:window-conditions}, also
\Cref{cond:upperelowerpsi} holds by \Cref{cond:schwartz}. 

It remains to prove \Cref{cond:speclocalization}.
Let $C_D$ be the doubling constant in
\Cref{cond:doublingintro}, and set
$\gamma=\log_2 C_D$.
Iterating \Cref{cond:doublingintro} gives
$G^{*q}(B_{1/T}(z))\gtrsim T^{-\gamma}$
and, for every $R\ge1$,
\begin{equation}
    \label{carino}
    G^{*q}(B_{R/T}(z))
\lesssim R^\gamma G^{*q}(B_{1/T}(z)).
\end{equation}
The lower bound from
\Cref{cond:upperelowerpsi}, \eqref{eq:variance-kernel} and \Cref{cond:doublingintro}
also give
$V_{T,q}(z)\gtrsim
T^{2d}G^{*q}(B_{1/T}(z))$.
Decomposing
$B_{M/T}^c(z)$ into annuli $2^{j}/T<d_{\wE}(\lambda,z)<2^{j+1}/T$,
by \Cref{lem:lattice-poisson} and \eqref{carino} we get for $N>\gamma/2$
\[
\frac{
\displaystyle\int_{B_{M/T}^c(z)}
\left|\wpsi_T(\lambda-z)\right|^2G^{*q}(d\lambda)}
{V_{T,q}(z)}
\lesssim_N \frac{T^{2d}G^{*q}(B_{1/T}(z))\sum_{2^j\gtrsim M}
2^{-(2N-\gamma)j}}{T^{2d}G^{*q}(B_{1/T}(z))}\lesssim \sum_{2^j\gtrsim M}
2^{-(2N-\gamma)j}
\]
Since the tail vanishes as $M\rightarrow\infty$,
\Cref{cond:speclocalization} holds and the proof of \Cref{thm:mainRdZd} is concluded by applying
\Cref{prop:varcontrol} and \Cref{thm:main}\textnormal{(iii)},
using that $\wpsi_T(0)\asymp T^d$ by \Cref{lem:Riemann-sum}.
\end{proof}

	\subsection{Proof of the criterion for indicator test functions}
	
		\label{rem:indicator-windows}
		Let $E=\R^d$ or $E=\Z^d$ and consider the regime \ref{large domain}, where $G=G_T$ and $\psi_T(x)=\psi(x/T)$. In \Cref{thm:mainRdZd} we assumed the doubling \Cref{cond:doublingintro} of $G^{*q}$ at $z$  and that $\psi$ is a Schwartz function with $\wpsi(0)>0$, see \Cref{cond:schwartz}. In this subsection, we show that the same criterion holds taking $\psi=\mathbf{1}_{D}$, with $D\subseteq\R^d$ suitable domain, under stronger assumptions on $G$. For \(\mathbf r=(r_1,\ldots,r_d)\in(0,\infty)^d\) and $z\in\wE$, let
		\[
		Q_{\mathbf r}(z)
		:=
		\left\{\lambda\in \widehat E:
		d_j(\lambda_j,z_j)\le r_j,\ j=1,\ldots,d\right\},
		\]
		where \(d_j\) is the Euclidean distance when \(\widehat E=\mathbb R^d\)
		and the toral distance when \(\widehat E=\mathbb T^d\). We say that
		\(G^{*q}\) has \emph{rectangular growth at \(z\)} if there exist
		\(C>0\) and \(\varepsilon_0>0\) such that
		\[
		G^{*q}\bigl(Q_{\varepsilon\mathbf R}(z)\bigr)
		\le
		C\left(\prod_{j=1}^dR_j\right)
		G^{*q}\bigl(Q_{\varepsilon\mathbf 1}(z)\bigr)
		\]
		for every \(\mathbf R=(R_1,\ldots,R_d)\in[1,\infty)^d\) and every
		\(\varepsilon\in(0,\varepsilon_0)\), where
		\(\mathbf1=(1,\ldots,1)\). A possible weakening of the rectangular growth condition is the \emph{ball growth} condition: there exists $C>0$ and $\eps_0>0$ such that
        \[
        G^{*q}\bigl(B_{\varepsilon R}(z)\bigr)
		\le
		C R^d
		G^{*q}\bigl(B_{\varepsilon}(z)\bigr)\,\qquad \forall \,\,R\ge 1\,,\,\,\eps\in(0,\eps_0)\,.
        \]
        We observe that rectangular growth $\;\Rightarrow\;$ ball growth $\;\Rightarrow\;$ \Cref{cond:doublingintro}. The first implication follows by inclusions of balls and cubes, while the second follows by definition, taking $R=2$.
        \begin{remark}\label{superem}
Let $\nu$ be a probability measure on
$\widehat E=\mathbb R^d$ or $\mathbb T^d$, and suppose that
$\widehat\nu\ge0$ on $E$. Then $\nu$ has rectangular
growth at the origin. In particular, $\widehat \nu=C^q\ge0$ $\implies$ $\nu=G^{*q}$ has rectangular growth at $z=0$.

To prove it, let $g_\eps$ be the normalized covariogram of $Q_{(\varepsilon/2)\mathbf1}(0)$ (see e.g. \cite{MainiRossiZheng})
\[
g_\eps(x)
=\frac{\wmu(Q_{(\varepsilon/2)\mathbf1}(0)\cap (Q_{(\varepsilon/2)\mathbf1}(0)+x)))}
       {\widehat\mu(Q_{(\varepsilon/2)\mathbf1}(0))}=\frac{\mathbf1_{Q_{(\varepsilon/2)\mathbf1}(0)}*
        \mathbf1_{Q_{(\varepsilon/2)\mathbf1}(0)}}
       {\Vol(Q_{(\varepsilon/2)\mathbf1}(0))}\,.
\]
For $\varepsilon$ sufficiently small when $\wE=\T^d$,
$g_\varepsilon\ge2^{-d}$ on
$Q_{(\varepsilon/2)\mathbf1}(0)$.
Since $\widehat g_\varepsilon\ge0$,
the function $h_\varepsilon=\nu*g_\varepsilon$ has non-negative Fourier transform and so is positive definite. Hence we have
\[
2^{-d}\nu(Q_{(\varepsilon/2)\mathbf1}(a))
\le h_\varepsilon(a)
\le h_\varepsilon(0)
\le\nu(Q_{\varepsilon\mathbf1}(0))\quad \implies \quad \nu(Q_{\varepsilon\mathbf R}(0))\lesssim \left(\prod_{j=1}^d R_j\right)\,\nu(Q_{\varepsilon\mathbf 1}(0))\,.
\]
\end{remark}
		We will assume that $z\in\supp(G^{*q})$ and that \(G^{*q}\) has ball growth or rectangular growth at \(z\). Then the conclusions of \Cref{thm:mainRdZd} when  \(z\in\supp(G^{*q})\) remain valid in each of the following cases:
		\begin{enumerate}[label=\textnormal{(\roman*)}]
			\item \(E=\mathbb Z^d, \R^d\), $G^{*q}$ has rectangular growth at $z$ and \(\psi=\mathbf1_{Q_{\mathbf{r}}(a)}\),
			$a\in E$ and \(\mathbf r=(r_1,\ldots,r_d)\in(0,\infty)^d\);
			
			\item \(E=\mathbb R^d\), $G^{*q}$ has ball growth at $z$ and \(\psi=\mathbf1_D\), where \(D\) is a
			convex body with nonempty interior and \(C^\infty\) boundary with non-vanishing
			Gaussian curvature (e.g. Euclidean balls).
		\end{enumerate} 
	
	\begin{theorem}[\Cref{thm:mainRdZd} for indicator functions]
		\label{prop:indicator-windows}
		Consider regime~\ref{large domain},
		fix \(q\ge2\), and suppose that $z\in\supp(G^{*q})$.
		Assume that $E,\psi,G^{*q}$ are as in (i) or (ii). Then \[
		\Var\left(F_{T,q}^{(z)}\right)
		\asymp q!\,T^{2d}G^{*q}(B_{1/T}(z)),
		\]
		and we have that $\widetilde F_{T,q}^{(z)}\overset{\rm law}{\rightarrow}N_z$ if and only if for every \(r=1,\ldots,q-1\),
			\begin{equation*}
				\lim_{\varepsilon\rightarrow0}
				\sup_{\lambda\in\widehat E}
				\frac{G^{*r}(B_\varepsilon(\lambda))
					G^{*(q-r)}(B_\varepsilon(z-\lambda))}
				{G^{*q}(B_\varepsilon(z))}=0\,.
			\end{equation*}
	\end{theorem}
	
	\begin{proof}
		We verify the assumptions of \Cref{thm:main}.
		Rectangular growth or ball growth imply \Cref{cond:doublingintro} and so \Cref{cond:doubling} holds.
		Thus, we can and will use
		\(G^{*q}(Q_{T^{-1}\mathbf1}(z))\) and
		\(G^{*q}(B_{1/T}(z))\) interchangeably up to constants. Moreover, in both cases $\psi$ is compactly supported, a.e. continuous and $\wpsi(0)>0$, so \Cref{cond:upperelowerpsi} holds by \Cref{lem:window-conditions}. 
        In particular, the lower bound in \Cref{cond:upperelowerpsi} and \Cref{cond:doubling} give $V_{T,q}(z)\gtrsim
T^{2d}G^{*q}(B_{1/T}(z))$.
It remains to verify
\Cref{cond:speclocalization}.

\emph{Case (i):} 
Here $\wpsi_T$ factors into one-dimensional Fourier transforms on $E$ of indicators of intervals of length $2Tr_j$. Therefore, by classical identities and bounds for Fourier transforms and Dirichlet kernels, we have
$|\wpsi_T(\theta_j)|\lesssim T
(1+Td_j(\theta_j,0))^{-1}$ and
\[
\left|\wpsi_T(\theta)\right|
\lesssim T^d\prod_{j=1}^d
(1+Td_j(\theta_j,0))^{-1}.
\]
Now we reason similarly as done in the proof of \Cref{thm:mainRdZd}, but with dyadic rectangles instead of dyadic annuli.
Decompose $\widehat E$ into
rectangles, indexed by
$\mathbf k\in\mathbb N_0^d$, on which
$2^{k_j}/T<d_j(\lambda_j,z_j)<2^{k_j+1}/T$
On the $\mathbf{k}$th rectangle $|\wpsi_T(\lambda-z)|^2\lesssim\prod_j2^{-2k_j}$,
while by rectangular growth the
$G^{*q}$-mass of the $\mathbf{k}$th rectangle is bounded by
$G^{*q}(B_{1/T}(z))\prod_j2^{k_j}$.
Consequently, as $M\rightarrow\infty$, uniformly for sufficiently large $T$
\[
\frac{\displaystyle
\int_{B_{M/T}^c(z)}
\left|\wpsi_T(\lambda-z)\right|^2G^{*q}(d\lambda)}
{V_{T,q}(z)}
\lesssim \frac{T^{2d}G^{*q}(B_{1/T}(z))\sum_{\max_j2^{k_j}\gtrsim M}
\prod_{j=1}^d2^{-k_j}}{T^{2d}G^{*q}(B_{1/T}(z))}\lesssim
\sum_{\max_j2^{k_j}\gtrsim M}
\prod_{j=1}^d2^{-k_j}
\rightarrow0
\]

\emph{Case (ii):}
By Fourier decay estimates for indicators of domains, see e.g.
\cite{BHI03}, we have
\[
\left|\wpsi_T(\theta)\right|
\lesssim T^d\left(1+T\|\theta\|\right)^{-(d+1)/2}.
\]
By ball growth, the same argument with dyadic annuli used in the proof of
\Cref{thm:mainRdZd} yields
\[
\frac{\displaystyle
\int_{B_{M/T}^c(z)}
\left|\wpsi_T(\lambda-z)\right|^2G^{*q}(d\lambda)}
{V_{T,q}(z)}
\lesssim
\sum_{2^j\gtrsim M}2^{-(d+1)j}2^{dj}
=
\sum_{2^j\gtrsim M}2^{-j}
\rightarrow0.
\]
Thus \Cref{cond:speclocalization} holds
in both cases. The conclusions follow from
\Cref{prop:varcontrol,thm:main}\textnormal{(iii)}.
	\end{proof}

\appendix

    \section{Variance asymptotics and chaotic reductions}\label{sec:additional proofs}
    In this section we collect the proof of \Cref{lem:variance-density} and \Cref{lem:finite-chaos-reduction},\Cref{cor:Ahlfors general}, which provide variance asymptotics and chaotic reductions used in \Cref{sec:large-domain-applications}.
    \begin{proof}[Proof of \Cref{lem:variance-density}]
Let $D_T=\R^d$ if $E=\R^d$, and
$D_T=T[-\pi,\pi)^d$ if $E=\Z^d$.
By change of variables in \eqref{eq:variance-kernel}, since $\wmu(du)=c_d \,du$, $c_d=(2\pi)^{-d}$ by \eqref{eq:muhat}, we get
\[
\frac{\Var\left(F_{T,q}^{(z)}\right)}{q!T^d}=c_d\,
\int_{D_T}
\left|T^{-d}\wpsi_T(u/T)\right|^2g_q(z+u/T)\,du.
\]
If $E=\R^d$,
$T^{-d}\wpsi_T(u/T)=\wpsi(u)$.
If $E=\Z^d$, by \Cref{lem:Riemann-sum} $T^{-d}\wpsi_T(u/T)\rightarrow\wpsi(u)$.
Moreover, by Plancherel
\[
c_d\int_{D_T}
\left|T^{-d}\wpsi_T(u/T)\right|^2\,du
=\int_{\wE}
T^{-d}\left|\wpsi_T(u)\right|^2\,\wmu(du)=T^{-d}\int_E\left|\psi_T(x)\right|^2\,\mu(dx)
\rightarrow\|\psi\|_2^2
=c_d\int_{\R^d}\left|\wpsi(u)\right|^2\,du.
\]
where the convergence follows by \Cref{lem:Riemann-sum} applied to $|\psi|^2$ if $E=\Z^d$. 
Since $\mathbf1_{D_T}
|T^{-d}\wpsi_T(\cdot/T)|^2$
converges pointwise to $|\wpsi|^2$
and the integrals converge, Scheffé's lemma
gives convergence in $L^1$.
 We conclude \eqref{varasymp} by
boundedness of $g_q$, continuity
of $g_q$ at $z$ and dominated convergence.
Moreover, the $L^1$ convergence also implies for every $\delta>0$
\[
R_{\delta,T}:=T^{-d}\int_{\|\lambda\|>\delta}\left|\wpsi_T(\lambda)\right|^2\wmu(d\lambda)=c_d\int_{D_T\cap\{\|u\|>\delta T\}}
\left|T^{-d}\wpsi_T(u/T)\right|^2\,du\rightarrow 0\,.
\]
Thus, if $z\ne-z$, choosing $\delta>0$ such that $\widehat E= \{d_{\widehat E}(\lambda,z)>\delta\} \cup \{d_{\widehat E}(\lambda,-z)>\delta\}$, Cauchy--Schwarz
gives
\[
\begin{aligned}
T^{-d}\left|\E\left[\left(F_{T,q}^{(z)}\right)^2\right]\right|
&\le q!\left\|g_q\right\|_\infty T^{-d}
\int_{\widehat E}
\left|\wpsi_T(\lambda-z)\wpsi_T(-\lambda-z)\right|
\,\wmu(d\lambda)\lesssim R_{\delta,T}^{1/2}\rightarrow 0\,.
\end{aligned}
\]
\end{proof}

\begin{proof}[Proof of \Cref{lem:finite-chaos-reduction}]
For $q\ge\max(m,2)$, we have $C\in L^q(E)$, so
\Cref{lem:variance-density} and \Cref{prop:short-memory-criterion}
show that either
$\Var\left(F_{T,q}^{(z)}\right)=o(\sigma_T^2)$ (if $g_q(z)=0$) or
$\widetilde F_{T,q}^{(z)}\Rightarrow N_z$.
If $m=1$, the same alternative holds for $q=1$ by
Gaussianity of the first chaos (see \Cref{sec:wiener-chaos}) and \Cref{lem:variance-density}.
Thus, by assumption the alternative holds for every $q$ with $c_q\neq0$.
Let $I:=\left\{q\ge1\,:\,c_q=0\text{ or }\Var\left(F_{T,q}^{(z)}\right)=o(\sigma_T^2)\right\}$. Define $H_{T}^{(z) }=\sum_{q\not\in I}c_qF_{T,q}^{(z)}$ and $\widetilde H_{T}^{(z)}:=H_{T}^{(z)}/\sqrt{\Var(H_T^{(z)})}$. Then by orthogonality \eqref{ortovar} and standard arguments (see e.g. \cite[Proof of Lemma 3.1]{MainiRossiZheng}) we have 
\[
\E\left[\left|\widetilde F_T^{(z)}-\widetilde H_T^{(z)}\right|^2\right]\lesssim \frac{\Var\left(F_T^{(z)}-H_T^{(z)}\right)}{\sigma_T^2}=\sum_{q\in I}c_q^2\frac{\Var\left(F_{T,q}^{(z)}\right)}{\sigma_T^2}\,.
\]
For $q\ge m$, by \eqref{varasymp} and
$\sigma_T^2\gtrsim T^d$, uniformly in $q$ and for sufficiently large $T$
\[
\frac{\Var\left(F_{T,q}^{(z)}\right)}{\sigma_T^2}
\lesssim\frac{\Var\left(F_{T,q}^{(z)}\right)}{T^d}\lesssim q!\|C^m\|_{L^1(E)}\,.
\]
Thus, breaking the sum in $q\in I\,,\,q<m$ and $q\in I\,,\,q\ge m$, since $\sum_q q!c_q^2<\infty$ by
\eqref{eq:phiexp} and $C\in L^m(E)$, dominated convergence gives \[
\E\left[\left|\widetilde F_T^{(z)}-\widetilde H_T^{(z)}\right|^2\right]\rightarrow 0\,.
\]
Thus, WLOG from now on we can and we will assume that for every $q$  we have $c_q\neq0$ and we do not have $\Var\left(F_{T,q}^{(z)}\right)=o(\sigma_T^2)$. In particular, by the alternative above $\widetilde F_{T,q}^{(z)}\overset{\rm law}{\rightarrow}N_z$ for every $q$. Set now
\[
S_{T,L}:=\sum_{q=1}^L c_qF_{T,q}^{(z)},
\qquad
R_{T,L}:=F_T^{(z)}-\E\left[F_T^{(z)}\right]-S_{T,L}.
\]
Since $|C|\le1$, for $q\ge L+1\ge m$ we have
$\|g_q\|_\infty\le\|C^m\|_{L^1(E)}<\infty$.
Thus, by \eqref{eq:phiexp}, \eqref{ortovar} and \eqref{varasymp}
\begin{equation*}
\limsup_{T\rightarrow\infty}\sigma_T^{-2}\Var(R_{T,L})\lesssim \limsup_{T\rightarrow\infty}T^{-d}\Var(R_{T,L})
\lesssim\|C^m\|_{L^1(E)}
\sum_{q>L}q!c_q^2
\rightarrow0,
\qquad L\to\infty.
\end{equation*}
By reasoning as in the comparison between $F_T^{(z)}$ and $H_T^{(z)}$, for every $T,L$ large enough we have 
\[
\E\left[\left|\widetilde F_T^{(z)}-\widetilde S_{T,L}\right|^2\right]\lesssim \frac{\Var\left(F_T^{(z)}-S_{T,L}\right)}{\sigma_T^2}=\frac{\Var(R_{T,L})}{\sigma_T^2}\lesssim\sum_{q>L}q!c_q^2\rightarrow 0\qquad\text{ as }L\rightarrow\infty\,.
\]
Thus, $\Var(S_{T,L})\gtrsim \sigma_T^2\gtrsim T^d$ for $T$ and $L$ large enough. Moreover, the proof is concluded if we prove \footnote{Indeed, in this case taking first the limit in $T$ and then in $L$ allows to conclude.} $\widetilde S_{T,L}:=S_{T,L}/\sqrt{\Var(S_{T,L})}\overset{\rm law}{\rightarrow}N_z$ for every $L$ large enough. Fix $L$ and write
\[
\widetilde S_{T,L}=\sum_{q=1}^L a_{T,q}\widetilde F_{T,q}^{(z)}\,,\qquad\quad a_{T,q}:=c_q\sqrt{\frac{\Var\left( F_{T,q}^{(z)}\right)}{\Var(S_{T,L})}}\,.
\]
By Peccati-Tudor \eqref{eq:peccatitudor}, applied on a vector of real and imaginary parts if $z\neq -z$, the convergence of $\widetilde F_{T,q}^{(z)}\overset{\rm law}{\rightarrow}N_z$ for every $q\le L$ with $c_q\neq0$ allows to conclude the joint convergence to a random vector with a finite number of i.i.d. copies of $N_z$. Let us now take a sequence $T_n$, $n\in \N$, $n\rightarrow\infty$ and consider $\widetilde S_{T_n,L}$. Since $a_{T_n}:=(a_{T_n,q})_{q\le L}\in\mathbb S^{L-1}\subseteq \R^L$ is a bounded sequence, we can extract a convergent subsequence from $a_{T_n}$, say $a_{T_{n_k}}\rightarrow a$ for some $a\in\mathbb S^{L-1}$. This, combined with the joint convergence to the vector of i.i.d. copies of $N_z$ mentioned before, allows to conclude that the subsequence $\sum_{q\le L}a_{T_{n_k},q}\widetilde F_{T_{n_k},q}^{(z)}$ of $\widetilde S_{T_n,L}$ is converging to $N_z$. Since this can be done for every choice of $T_n$, we conclude that $\widetilde S_{T,L}\overset{\rm law}{\rightarrow}N_z$.
\end{proof}

\begin{proof}[Proof of \Cref{cor:Ahlfors general}]
Since $c_2\neq0$, $z=0$, $C\notin L^2$ and $C\in L^m$ for $m$ large enough, by reasoning as in
\cite[Proposition 2.2]{MainiNourdinSpectralCLT},
the contribution of the chaoses of order $q\ge3$
is negligible, in the sense that $\Var(F_T^{(0)})=\Var(S_T)(1+o(1))$ and $\widetilde F_T^{(0)}-S_T/\sqrt{\Var(S_T)}\overset{L^2(\Omega)}{\rightarrow}0$, where $S_T:=c_1F_{T,1}^{(0)}+c_2F_{T,2}^{(0)}$. In particular, the variance asymptotics follow by \eqref{ortovar} and \Cref{prop:second-chaos-spectral-set}.
Set $\sigma_T^2:=\Var(S_T)$ and
\[
v_{T,i}:=\frac{c_i^2\Var(F_{T,i}^{(0)})}{\sigma_T^2},
\qquad i=1,2.
\]
By orthogonality \eqref{ortovar}, $v_{T,1}+v_{T,2}=1$.
Along any sequence $T_n\to\infty$, we can extract a
subsequence such that $v_{T_{n_k},i}\to v_i$, $i=1,2$.
Since the first chaos is Gaussian, possibly degenerate,
and the second satisfies a CLT by
\Cref{prop:second-chaos-spectral-set},
Peccati--Tudor \eqref{eq:peccatitudor} yields the joint convergence to a standard Gaussian vector with variances $v_1,v_2$, yielding
\[
\frac{S_{T_{n_k}}}{\sigma_{T_{n_k}}}
\overset{\mathrm{law}}{\longrightarrow}
\mathcal N(0,v_1+v_2)=N.
\]
Since $T_n\to\infty$ is arbitrary, we have $S_T/\sqrt{\Var(S_T)}\overset{\mathrm{law}}{\longrightarrow} N$ and a CLT for $F_T^{(0)}$.
\end{proof}

\section{Coarea formula, Ehresmann's theorem and proof of \Cref{lem:twofold-convolution}}\label{sec:ehresmann}

We recall some geometric tools used in the
proof of \Cref{lem:twofold-convolution}.
We formulate them for smooth submanifolds
of Euclidean spaces, which is the setting
needed below.

Let $\mathcal M\subset\R^a$ and
$\mathcal N\subset\R^b$ be smooth embedded
submanifolds of dimensions $m\ge n\ge1$. 
Let
$T_x\mathcal M\subseteq \R^a$ be the tangent space of $\mathcal{M}$ in $x$.
The Riemannian metric induced by $\R^a$ on $\mathcal{M}$ is simply the restriction of the Euclidean inner product
to each tangent space, namely
\[
\langle u,v\rangle_{T_x\mathcal M}
=\langle u,v\rangle_{\R^a},
\qquad u,v\in T_x\mathcal M,
\]
and analogously for $\mathcal N$. Let $f:\mathcal M\to\mathcal N$ be
smooth and let $Df_x:T_x\mathcal M\to T_{f(x)}\mathcal N$ be its differential at $x$.
We write $Df_x^*:T_{f(x)}\mathcal N\to T_x\mathcal M$ for its adjoint,
satisfying
\[
\langle Df_xu,v\rangle
=\langle u,Df_x^*v\rangle,
\qquad
u\in T_x\mathcal M,\quad
v\in T_{f(x)}\mathcal N,
\]
The tangential Jacobian of $f$ is
defined by
\[
J_f(x):=
\sqrt{\det(Df_xDf_x^*)},
\]
In particular, $J_f(x)>0$ exactly
when $Df_x$ is surjective. 

A crucial tool will be the co-area formula
\cite[Theorem~3.2.22]{FedererGMT}: for every nonnegative
measurable $h$
\begin{equation}\label{eq:coarea-prelim}
\int_{\mathcal M}h(x)J_f(x)
\,d\mathcal H^m(x)
=
\int_{\mathcal N}
\int_{f^{-1}(z)}
h(x)\,d\mathcal H^{m-n}(x)
\,d\mathcal H^n(z)\,,
\end{equation}
were $\mathcal H^k$ is the Hausdorff measure.
\eqref{eq:coarea-prelim}
also holds when $Df_x$ fails to be
surjective at some points\footnote{By Sard's theorem
\cite[Theorem~6.10]{LeeSmoothManifolds}, $J_f>0$ on $f^{-1}(z)$ for $\mathcal{H}^n$ a.e. $z\in\mathcal{N}$. For such $z$, the regular
value theorem ensures that
$f^{-1}(z)$ is a smooth
submanifold of dimension $m-n$.}.

We shall also use the area formula,
corresponding to $m=n$.
If $f:\mathcal M\to\mathcal N$
is a diffeomorphism between
$k$-dimensional embedded submanifolds,
then, for every nonnegative
measurable function $h$
\begin{equation}\label{eq:area-prelim}
\int_{\mathcal N}h(y)
\,d\mathcal H^k(y)
=
\int_{\mathcal M}
h(f(x))J_f(x)
\,d\mathcal H^k(x),
\end{equation}
In dimension $k=0$,
$\mathcal H^0$ is the counting measure
and $J_f=1$.

Finally, we recall Ehresmann's fibration
theorem \cite[Theorem~8.5.10]
{DundasDifferentialTopology} and its consequences.
A smooth map $f:\mathcal M\to\mathcal N$
is called a submersion if $Df_x$ is
surjective at every point, equivalently
if $J_f(x)>0$ for every $x\in\mathcal M$.
It is called proper if the preimage of
every compact set is compact.

\begin{lemma}[Smooth variation of fibres]
\label{lem:smooth-fibres}
Let $\mathcal M\subset\R^a$ and
$\mathcal N\subset\R^b$ be smooth embedded
submanifolds of dimensions $m\ge n$, and
let $f:\mathcal M\to\mathcal N$ be a proper
smooth submersion. Then the following hold.
\begin{enumerate}
\item Every nonempty fibre $f^{-1}(z)$ is
a smooth embedded submanifold of dimension
$m-n$.

\item For every $z_0\in f(\mathcal M)$,
there exist an open neighbourhood $V$ of
$z_0$ and a diffeomorphism
\[
\Psi:f^{-1}(V)\to V\times f^{-1}(z_0),
\qquad \pi_1\circ\Psi=f.
\]
Writing $\Phi:=\Psi^{-1}$, the following restrictions are diffeomorphisms
\[
\Phi_z:f^{-1}(z_0)\to f^{-1}(z),
\qquad \Phi_z(u):=\Phi(z,u),
\quad z\in V,
\]

\item The map $(z,u)\mapsto\Phi_z(u)$
is smooth on $V\times f^{-1}(z_0)$.
Moreover, the tangential Jacobians
$J_{\Phi_z}(u)$ are strictly positive
and continuous in $(z,u)\in V\times f^{-1}(z_0)$.

\item If $\mathcal N$ is connected and
$f$ is surjective, then all fibres
$f^{-1}(z)$, $z\in\mathcal N$, are
diffeomorphic.
\end{enumerate}
\end{lemma}

\begin{proof}
Since $f$ is a submersion, every
$z\in\mathcal N$ is a regular value.
The regular value theorem therefore
implies that each nonempty fibre is
a smooth embedded submanifold of
dimension $m-n$, proving (1).

The first assertion in (2) is Ehresmann's theorem.
Since $\pi_1\circ\Psi=f$, for every
$z\in V$ $\Phi_z$ is a bijection,
with inverse
$\pi_2\circ\Psi|_{f^{-1}(z)}$.
Both maps are smooth, being restrictions
or compositions of smooth maps.
This proves (2). Since $\Phi_z(u)=\Phi(z,u)$ and
$\Phi$ is smooth, the map
$(z,u)\mapsto\Phi_z(u)$ is smooth.
Thus the differentials $D\Phi_z(u)$
depend continuously on $(z,u)$ when
expressed in local coordinate and $(z,u)\mapsto J_{\Phi_z}(u)$ is continuous.
Moreover,
since $\Phi_z$ is a diffeomorphism $J_{\Phi_z}>0$, which concludes (3).
Finally, if $\mathcal N$ is connected,
any two of its points can be joined
by a path. Covering this path by
finitely many neighbourhoods as above gives
diffeomorphisms between all the fibres,
proving (4).
\end{proof}

We also need some properties on the Gauss map $n:M\rightarrow\mathbb S^{d-1}$, defined as follows: $\forall x\in M$, $n(x)$ is the only normal vector to $T_xM$ such that $x+t\,n(x )\notin K$ $\forall t>0$.

\begin{lemma}[Gauss map]
\label{lem:gauss-map}
Let $K\subset\R^d$, $d\ge2$, be a
compact convex set with nonempty interior
and smooth boundary $M=\partial K$.
Let $n:M\to\mathbb S^{d-1}$ be the
outward Gauss map.

If $K$ is strictly convex, then $n$ is bijective.
If additionally $M$ has non-vanishing
Gaussian curvature, then $n$ is a
diffeomorphism. If $K=-K$, then
$n(-x)=-n(x)$ for all $x\in M$.
\end{lemma}

\begin{proof}
Fix $u\in\mathbb S^{d-1}$.
By compactness, the function
$x\mapsto\langle x,u\rangle$
attains its maximum on $K$ at some
$x_u\in M$\footnote{Indeed, an interior maximum
would lead to a contradiction $\langle x+tu,u \rangle>\langle x,u\rangle\in K$ for $t$ small enough.}.
Since $x_u$ is also a maximizer on $M$,
$\langle v,u\rangle=0$
for every $v\in T_{x_u}M$.
Moreover,
\[
\langle x-x_u,u\rangle\le 0
\qquad\forall\,x\in K\quad \implies \quad x_u+t\,u\not\in K\,\quad \forall t>0
\]
Thus $n(x_u)=u$, proving that $n$ is
surjective. Moreover, by convexity and smoothness, $\forall$ $x\in M$,
\[
\langle y-x,n(x)\rangle\le0,
\qquad \forall \,y\in K.
\]
Indeed, otherwise the segment
$x+t(y-x)\in K$ would be out of $K$ for $t$ small enough.
Consequently, if $n(x)=n(y)=u$,
both $x$ and $y$ maximize
$\langle\cdot,u\rangle$ on $K$.
If $x\ne y$, their midpoint also
maximizes this function, contradicting
strict convexity, since an interior
point cannot be a maximizer.
Hence $n$ is injective.

Since $M$ is a smooth hypersurface,
the outward Gauss map $n$ is smooth. Moreover, by assumption the Gaussian
curvature is nowhere zero, so $Dn_x$
is invertible at every $x\in M$. Thus, by the inverse function theorem, see e.g.
\cite[Theorem~4.5]{LeeSmoothManifolds},
$n$ is a local diffeomorphism.
Since it is bijective, it is also a
(global) diffeomorphism, see e.g. \cite[Proposition~4.6(f)]{LeeSmoothManifolds}.
Finally, if $K=-K$, then $T_{-x}\mathcal{M}=T_x\mathcal{M}$ and so $n(-x)\in\{\pm n(x)\}$. If $n(x)$ is the outward normal, then $x+tn(x)\notin K$ for $t>0$, that implies $-x+t(-n(x))\notin K$ and $-n(x)=n(-x)$.
\end{proof}

We now apply these tools to the
addition map on a hypersurface $M$ as in \Cref{lem:twofold-convolution}.

\begin{lemma}
\label{lem:addition-map}
Fix $d\ge2$ and let $K\subseteq\R^d$,
$M=\partial K$ be as in
\Cref{lem:twofold-convolution}.
Consider $F:M\times M\to\R^d$,
$F(x,y)=x+y$, and let
$n:M\to\mathbb S^{d-1}$ be the
outward Gauss map.
Then the tangential Jacobian of $F$ is
\begin{equation}\label{JF}
J_F(x,y)
=2^{(d-2)/2}
\bigl(1-\langle n(x),n(y)\rangle^2
\bigr)^{1/2}.
\end{equation}
Moreover, the critical set of $F$ is
\begin{equation}\label{2.2}
\Sigma:=\{(x,y)\in M\times M:
J_F(x,y)=0\}
=\{(x,x):x\in M\}
\cup\{(x,-x):x\in M\}.
\end{equation}
Finally, if
$O:=\operatorname{int}(2K)\setminus\{0\}$, then $F^{-1}(O)=M\times M\setminus\Sigma$ and 
$F:F^{-1}(O)\to O$ is a proper surjective smooth submersion.
\end{lemma}

\begin{proof}
Since $M$ is a smooth hypersurface, the tangent space to the product
$M\times M$ at $(x,y)$ is
$T_xM\times T_yM =n(x)^\perp \times n(y)^\perp$,
Since $F(x,y)=x+y$, its differential is $DF_{(x,y)}:n(x)^\perp \times n(y)^\perp\rightarrow \R^d$
\[
DF_{(x,y)}(v,w)=v+w,
\qquad
v\in n(x)^\perp,\quad w\in n(y)^\perp.
\]
Moreover, standard computations show that the adjoint $DF^*_{(x,y)}:\R^d\rightarrow n(x)^\perp \times n(y)^\perp$ is given by $DF^*_{(x,y)}(w)=(P_xw,P_yw)$, where $P_x:=I-n(x)n(x)^t$ is the projection map on $n(x)^\perp$.
Thus, we get
\[
DF_{(x,y)}DF_{(x,y)}^*
=P_x+P_y=2I-n(x)n(x)^t-n(y)n(y)^t.
\]
Standard arguments about eigenvalues allow to compute the determinant of $J_F$, proving \eqref{JF}.
By \eqref{JF}, $J_F(x,y)=0$ $\iff$
$n(y)=\pm n(x)$. Moreover, by \Cref{lem:gauss-map} $n$ is bijective and $n(-x)=-n(x)$, which proves \eqref{2.2}.
By convexity,
$F(M\times M)\subseteq2K$.
Moreover, strict convexity implies
that if $x,y\in M$, $x\neq y$,
then $(x+y)\in\operatorname{int}(2K)$, so $F^{-1}(\partial(2K))=F^{-1}(2M)=\{(x,y)\in M\times M\,;\,x=y\}$.
Since $F^{-1}(0)=\{(x,y)\in M\times M\,;\,x=-y\}$, we conclude that $F^{-1}(O)=M\times M\setminus\Sigma$.
It remains to prove that $F:M\times M\setminus \Sigma\rightarrow O$ is a proper surjective submersion. By \eqref{JF}, $J_F>0$ on $F^{-1}(O)$, so the
restriction $F:F^{-1}(O)\to O$ is a smooth
submersion. Moreover, its image $F(F^{-1}(O))$ is open in $O$, since smooth
submersions are open maps
(see e.g. \cite[Proposition~4.28]{LeeSmoothManifolds}),
and closed in $O$, since
$F(F^{-1}(O))=F(M\times M)\cap O$ and
$F(M\times M)$ is compact by continuity of $F$.
Moreover, the image is nonempty and $O$ is
connected, so we have $F(F^{-1}(O))=O$ and $F$ surjective.
Finally, for $B\subset O$ compact, $F^{-1}(B)$ is compact, since it is closed in $M\times M$
 compact. Thus the restriction is proper.
\end{proof}

We are finally ready to prove
\Cref{lem:twofold-convolution}.

\begin{proof}[Proof of
\Cref{lem:twofold-convolution}]

Since $M$ is a compact smooth hypersurface,
it admits finitely many local graph charts
with uniformly bounded derivatives.
The area formula then gives $\mathcal H^{d-1}(M\cap B_r(x))
\asymp r^{d-1}$
uniformly for $x\in M$ and
$0<r<r_0$, proving the Ahlfors bounds.

Let $F$, $J_F$, $\Sigma$ and $O$
be as in \Cref{lem:addition-map}.
Since $F^{-1}(O)=(M\times M)\setminus\Sigma$
and $\mathcal H^{2d-2}(\Sigma)=0$, we have
\[
G^{*2}(\R^d\setminus O)
=(G\otimes G)(\Sigma)
=\frac{\mathcal H^{2d-2}(\Sigma)}
{\mathcal H^{d-1}(M)^2}=0.
\]
Applying the coarea formula
\eqref{eq:coarea-prelim}
to the restriction
$F:F^{-1}(O)\to O$,
with $h=\mathbf1_{F^{-1}(A)}/J_F$,
we obtain for every Borel set
$A\subset O$ (recall the choice \eqref{eq:muhat})
\begin{equation}\label{g2density}
G^{*2}(A)=\int_Ag_2(z)\,\wmu(dz),
\qquad
g_2(z):=
\frac{(2\pi)^d}{\mathcal H^{d-1}(M)^2}
\int_{F^{-1}(z)}
\frac{d\mathcal H^{d-2}}{J_F},
\quad\text{for a.e. }z\in O.
\end{equation}
We next show that $g_2$ admits a continuous
and strictly positive version on $O$.
By \Cref{lem:addition-map,lem:smooth-fibres}, the fibres
$F^{-1}(z)$, $z\in O$, are compact and
nonempty, and $J_F>0$ on $F^{-1}(O)$.
We may therefore define $g_2(z)$ for every
$z\in O$ by
\eqref{g2density}. Fix $z_0\in O$. By
\Cref{lem:addition-map,lem:smooth-fibres},
there exist a neighbourhood $V$ of $z_0$
and a smooth map
\[
\Phi:V\times F^{-1}(z_0)\to F^{-1}(V)
\]
such that, for every $z\in V$,
$\Phi_z(u):=\Phi(z,u)$ defines a
diffeomorphism from $F^{-1}(z_0)$
onto $F^{-1}(z)$.

By the area formula \eqref{eq:area-prelim},
for every $z\in V$ we have
\[
g_2(z)=
\frac{(2\pi)^d}{\mathcal H^{d-1}(M)^2}
\int_{F^{-1}(z_0)}
\frac{J_{\Phi_z}(u)}
{J_F(\Phi_z(u))}
\,d\mathcal H^{d-2}(u).
\]
By \Cref{lem:addition-map,lem:smooth-fibres}, both the numerator and denominator
are continuous and strictly positive
in $(z,u)$. Since $F^{-1}(z_0)$ is
compact and nonempty, $g_2$ is
continuous and strictly positive on $V$.
Since $z_0\in O$ was arbitrary, $g_2$ is continuous and positive on $O$.

Assume now that $d=2$. By
\Cref{lem:smooth-fibres}, the fibres
$F^{-1}(z)$, $z\in O$, are all diffeomorphic manifolds of dimension $0$, hence they have the same
finite positive cardinality.
Since $n$ is a diffeomorphism by \Cref{lem:gauss-map}
and $M$ is compact, $n$ and its inverse are both (global) Lipschitz functions on $M$. Using $n(-y)=-n(y)$
and \eqref{JF}, we get for some
$\delta>0$ and uniformly for $x,y\in M$ with
$0<\|x+y\|<\delta$,
\[
J_F(x,y)\asymp\|n(x)+n(y)\|
\asymp\|x+y\|,
\]
Since $\mathcal H^0$ is the counting
measure, \eqref{g2density} yields
$g_2(z)\asymp\|z\|^{-1}$ for $0<\|z\|<\delta$.
\end{proof}

\makeatletter
\let\@tocwrite\@gobbletwo
\makeatother

	\section*{Acknowledgements}
    The author is grateful to E. Di Bernardino, L. Gass, F. Grotto, I. Nourdin, N. Leonenko, D. Marinucci, G. Peccati, M. Rossi, R. Shevchenko, M. Stecconi and A. P. Todino for the useful discussions.
    
    The author is associated with INdAM (Istituto Nazionale di Alta Matematica “Francesco Severi”) and the GNAMPA group.
	This work was partially supported by the MUR Excellence Department Project MatMod@TOV, awarded to the Department of Mathematics, University of Rome Tor Vergata (CUP E83C18000100006). We also acknowledge financial support from the MUR 2022 PRIN project GRAFIA (project code 202284Z9E4).

\section*{Declaration of AI use}
The author used ChatGPT 5.6 and 6 as a supporting tool for: language editing, generation of \Cref{fig:spectral-delocalization}, bibliographic exploration of related/needed results, refinement and discussion of mathematical results. 
The author takes full responsibility for the mathematical arguments, their correctness, and the final content of the manuscript.

	\bibliographystyle{alpha}
		
\end{document}